\documentclass[12pt]{article}

\usepackage{amssymb,mathtools,amsthm, amsmath}
\usepackage{fullpage}
\usepackage[small, labelfont = bf]{caption}
\usepackage{tikz-cd}
\usepackage[numbers,sort]{natbib}

\newcommand{\blocktheorem}[1]{%
  \csletcs{old#1}{#1}% Store \begin
  \csletcs{endold#1}{end#1}% Store \end
  \RenewDocumentEnvironment{#1}{o}
    {\par\addvspace{1.5ex}
     \noindent\begin{minipage}{\textwidth}
     \IfNoValueTF{##1}
       {\csuse{old#1}}
       {\csuse{old#1}[##1]}}
    {\csuse{endold#1}
     \end{minipage}
     \par\addvspace{1.5ex}}
}

\newcommand{\BD}       {{\mathrm{BD}}}
\newcommand{\bd}       {{\mathrm{bd}}}
\newcommand{\rk}        {{\mathrm{rk}}}
\newcommand{\im}        {{\mathrm{im}}\,}
\newcommand{\Vect}      {{\mathsf{Vec}}}
\newcommand{\Hilb}      {{\mathsf{Hilb}}}
	
\newcommand{\Dgm}      {{\mathrm{Dgm}}}

\newtheoremstyle{amit}% name
{7pt}% Space above
{7pt}% Space below
{}% Body font
{0pt}% Indent amount
{\bf}% Theorem head font
{:}% Punctuation after theorem head
{.5em}% Space after theorem head
{}% Theorem head spec (can be left empty, meaning `normal')

\theoremstyle{amit}
\newtheorem{defn}{Definition}[section]
\newtheorem{prop}[defn]{Proposition}
\newtheorem{lem}[defn]{Lemma}
\newtheorem{thm}[defn]{Theorem}
\newtheorem{cor}[defn]{Corollary}
\newtheorem{rmk}[defn]{Remark}
\newtheorem{ex}[defn]{Example}

\begin{document}

\title{The M\"obius Homology Laplacian for Persistence Modules}

\author{Tatum Rask \\ \normalsize{University of Vermont, Department of Mathematics and Statistics}
}
\date{}

\maketitle

\begin{abstract}
Persistent homology provides a framework for studying the evolution of topological features across filtered spaces. While persistence diagrams give discrete invariants that record the birth and death of homological features, they do not directly encode analytic or geometric structure associated to persistence modules. This paper develops a spectral framework for persistent homology by introducing and studying a Laplace operator on the M\"obius chain complex, called the \textit{M\"obius homology Laplacian}.

M\"obius homology categorifies M\"obius inversion for persistence modules by replacing integer-valued data with vector-space-valued chain complexes. By equipping these chain complexes with inner products, we define a Laplace operator on the M\"obius chain complex whose kernel provides canonical cycle representatives for the associated M\"obius homology space. We analyze the spectral properties of this operator in two primary settings: one-parameter persistence modules and birth-death modules arising from filtrations of simplicial complexes. For one-parameter modules, we completely characterize the eigenvalues of the degree-one M\"obius homology Laplacian. For birth-death modules, we show that the trace and determinant of the degree-zero M\"obius homology Laplacian reflect the combinatorial structure of the persistence diagram. 
\end{abstract}

%%%%%%%%%%%%%%%%%%%%%%%%%%%%%%%%%%%%%%%%%%%%%%
%%%%%%%%%%%%%%%%%%%%%%%%%%%%%%%%%%%%%%%%%%%%%%
\section{Introduction}
\label{sec:introduction}
%%%%%%%%%%%%%%%%%%%%%%%%%%%%%%%%%%%%%%%%%%%%%%
%%%%%%%%%%%%%%%%%%%%%%%%%%%%%%%%%%%%%%%%%%%%%%
Persistent homology was developed as a tool for extracting topological information from data \cite{Landi97, Robins99, ELZ02}. In topological data analysis, one often constructs a filtered simplicial complex from a point cloud of data (for example, the Vietoris-Rips or \v{C}ech complexes  \cite[Section 2.2]{carlsson09}) and tracks the evolution of homological features as the filtration value varies. A key feature of persistent homology is \textit{stability}: persistence diagrams vary continuously under perturbations of the data, ensuring robustness to noise \cite{stability-persistence, chazal09proximity}.  

Although motivated by data analysis, persistent homology has developed into a rich mathematical theory. Filtered simplicial complexes give rise to \textit{persistence modules}: sequences of vector spaces with linear maps between them. In the one-parameter setting, persistence modules admit complete discrete invariants in the form of persistence diagrams or, equivalently, barcodes \cite{CaZo09}. Persistence diagrams record the birth and death of homological features, and they can be computed using algebraic and combinatorial methods such as rank invariants \cite{CaZo09}, birth-death functions \cite{McCPa20-edit-distance}, and interval decompositions \cite{ZoCa05}. These constructions are equivalent in the one-parameter case. However, all three constructions are combinatorial in nature. They record when features exist, but they do not provide canonical cycle representatives or analytic structure. This work is motivated by the following question:

\begin{quote}
\textit{Can persistence admit an intrinsic analytic structure, analogous to Hodge theory in classical topology?}
\end{quote}

In classical topology, algebraic and combinatorial invariants are enriched by analytic methods such as the Hodge Laplacian on a Riemannian manifold or its combinatorial counterpart on a simplicial complex. The kernel of the Hodge Laplacian recovers homology spaces and provides canonical harmonic representatives for each homology class, and its non-zero spectrum reflects geometric and combinatorial structure \cite{hodge1941theory, Eckmann1944, KirchhoffUeberDA, HORAK2013}. Persistence modules have not traditionally admitted such operators. This paper develops a framework in which an analytic structure naturally emerges. To do so, we utilize \textit{M\"obius homology}, a categorification of M\"obius inversion.

M\"obius inversion is a classical tool for extracting local contribution from cumulative data on finite posets \cite{Rot64}. In persistence theory, rank and birth-death functions encode cumulative information about the presence of homological features across intervals. The persistence diagram can be computed by applying M\"obius inversion to these functions.

M\"obius homology categorifies M\"obius inversion \cite{patel2023mobius}. Instead of assigning integers to intervals, it assigns vector spaces organized into a chain complex indexed by intervals. 
In the one-parameter setting, M\"obius homology recovers the same interval decomposition data as barcodes (see Section~\ref{sec:mob-homology-calculation}). In the case of multiparameter modules, where no complete discrete invariants exist, M\"obius homology can distinguish non-isomorphic modules that are indistinguishable by M\"obius inversion and persistence diagrams alone. The additional structure of chain complexes, together with the introduction of an inner product, allow for the computation of adjoint operators and, ultimately, a combinatorial Laplacian acting on the M\"obius chain spaces (see Section~\ref{sec:mob-laplacian}). 

The central idea of this work is that M\"obius homology admits a natural Laplace operator, called the \textit{M\"obius homology Laplacian}, that provides analytic structure for persistence modules.

%%%%%%%%%%%%%%%%%%%%%%%%%%%%%%%%%%%%%%%%%%%%%%
\subsection{Related Work}
\label{sec:related-work}
%%%%%%%%%%%%%%%%%%%%%%%%%%%%%%%%%%%%%%%%%%%%%%

This work presents a novel perspective on the analytic and spectral properties of persistence modules utilizing the tools of M\"obius homology \cite{patel2023mobius}. Our work builds upon and complements existing literature in combinatorial Hodge theory, harmonic persistent homology and persistent Laplacians, and monoidal M\"obius inversions for birth-death spaces. 

\paragraph{Combinatorial Hodge Laplacians}
This work follows a line of literature on the combinatorial Hodge Laplacian and its spectrum. The study started with Kirchhoff's paper on the graph Laplacian \cite{KirchhoffUeberDA} and was popularized by Fiedler, who studied properties of its smallest non-zero eigenvalue \cite{Fiedler1973}. Eckmann was the first to generalize the graph Laplacian to simplicial complexes, introducing the combinatorial Laplacian \cite{Eckmann1944}.
Hansen and Ghrist extend this work even further, defining the cellular sheaf Laplacian \cite{Hansen2019spectral}.

The spectrum of the graph Laplacian is known to reflect the combinatorial structure of the underlying graph. For example, Kirchhoff proves the matrix-tree theorem in \cite{KirchhoffUeberDA}, relating the non-zero eigenvalues to the number of spanning trees in a graph. Biggs extends the matrix-tree theorem, relating forests of various sizes to the characteristic polynomial of the graph Laplacian \cite{Biggs93AlgebraicGraphTheory}. The Cheeger inequality relates the \textit{conductance} of a graph to its second-smallest eigenvalue \cite{Cheeger1971, FanChungCheeger}. Fisher studied the \textit{normalized} graph Laplacian (as a discrete analogue to the Laplace-Beltrami operator) and its connection to random walks on lattices in \cite{FISHER1966105}. Moreover, Anderson and Morley prove important results about the spectrum of the graph Laplacian in \cite{Anderson01101985}. 

Similar results for the combinatorial Laplace operator connect the spectrum to the combinatorial and geometric structure of the underlying simplicial complex. Duval and Reiner, for one, show that \textit{shifted} simplicial complexes always have integer eigenvalues that arise from an associated integer partition \cite{duval2002shifted}. 
Furthermore, Horak and Jost explicitly compute an upper bound on the spectrum and, in the case of \textit{normalized} combinatorial Laplacians, provide a full characterization of the eigenvalues for certain types of simplicial complexes \cite{HORAK2013}. Steenbergen dives further into the spectrum of combinatorial Laplacians in his doctoral dissertation \cite{Steenbergen}, and in \cite{Bernardi16, Contreras24} the characteristic polynomial of the combinatorial Laplacian is related to higher-dimensional forests and discrete Morse theory. 

For the cellular (co)sheaf Laplacian, Hansen and Ghrist provide a first spectral analysis in \cite{hansen20phd, Hansen2019spectral}. A Cheeger-type inequality for cellular sheaves is proved in \cite{guadagnuolo2026}, among other interesting spectral results. However, because a cellular (co)sheaf provides so much additional data, the spectrum of its associated Laplace operator is, in general, much more mysterious. In this work, we study a special case of the cellular cosheaf Laplacian. In particular, we utilize the inherent structure of the M\"obius chain complex, obtained from the \textit{order cosheaf} of a persistence module, to relate the spectrum of its associated cellular cosheaf Laplacian to the underlying module.

\paragraph{Spectral Methods in Persistent Homology}
Given the success of the combinatorial Laplacian as a method for studying graphs and simplicial complexes, it should come as no surprise that such spectral methods have been extended to \textit{filtered} simplicial complexes. This study originated in a 2014 talk on persistent harmonic forms by Andr\'e Lieutier \cite{lieuter_talk}. Since then, authors such as Wang, Nguyen, and Wei \cite{wang20persistentspectralgraphs} and Mémoli, Wan, and Wang \cite{Memloi2022Laplacian} have investigated the persistent Laplacian further. The latter paper, for example, presents interesting results about the spectrum of the persistent Laplacian: in particular, the authors prove a persistent version of the Cheeger inequality. They also prove that the eigenvalues are monotonic and stable. Furthermore, Wei and Wei define the persistent \textit{sheaf} Laplacian in \cite{wei2025persistentsheaf}.

On the other hand, Basu and Cox present an analysis of harmonic persistent homology spaces in \cite{basu2022harmonic}. In particular, harmonic persistent homology spaces provide canonical cycle representatives for barcodes that maximize ``relative essential content." We will show a similar result, heavily inspired by Basu and Cox's work, in Section~\ref{subsec:harmonic-reps}. More recently, Basu, Guzm\'an-S\'aenz, and Parida show that essential simplices ``dominate" harmonic cycle representatives \cite{basu2026dominate}. 

The approaches above act on filtered chain complexes at individual parameter values whereas
the M\"obius homology Laplacian acts on the \textit{M\"obius} chain complex, whose construction depends on the combinatorial structure of the underlying indexing poset. 

\paragraph{Grassmannian Persistence Diagrams and Orthogonal M\"obius Inversion} 

 Another approach for obtaining cycle representatives for persistence diagrams and utilizing vector space-level data appears in Gülen, Mémoli, and Wan's work on the Grassmannian persistence diagram \cite{gulen2025grassmannian, gulen2025grassmannian1d}. In their work, the authors define monoidal M\"obius inversion and apply it to both the birth-death spaces and the kernel of the persistent Laplacian. The result is a vector space of cycle representatives for the persistence diagrams. In Appendix~\ref{appendix:lit-review}, we look into the relationship between our work and the Grassmannian persistence diagram. 

%%%%%%%%%%%%%%%%%%%%%%%%%%%%%%%%%%%%%%%%%%%%%%
\subsection{Contributions}
%%%%%%%%%%%%%%%%%%%%%%%%%%%%%%%%%%%%%%%%%%%%%%
The overarching contribution of this work is establishing foundational elements of a spectral theory for persistence modules. All contributions are also present in the author's doctoral dissertation \cite{rask26phd}.

More specifically, we first define the M\"obius homology Laplacian for persistence modules indexed by finite posets. 
As to be expected, the kernel of this operator is canonically isomorphic 
to M\"obius homology (see Theorem~\ref{prop:harmonic_cosheaf_reps}). 

Second, we provide a complete characterization of the eigenvalues of the degree-one M\"obius homology Laplacian in the one-parameter setting. The eigenvalues of this Laplace operator are determined by the barcode of the underlying persistence module (see Theorem~\ref{thm:spectrum_of_L1b}).

Third, when persistence modules arise from simplicial filtrations, the analytic viewpoint reflects the geometry of the underlying simplicial complex and the combinatorics of the persistence diagram. This is discussed in detail in Section~\ref{sec:birth-death-laplacian}. The degree-zero Laplacian acts on a birth-death module constructed directly from 
the filtration, and harmonic representatives emphasize ``newly added simplices" (see Theorem~\ref{prop:relative-new-content}). Moreover, spectral 
invariants such as the trace and determinant admit combinatorial descriptions in terms of the birth-death function (see Theorems~\ref{trace-of-L0ab} and~\ref{prop:lower-bound-L0ab}).

In all, this work introduces the M\"obius homology Laplace operator, establishes its fundamental properties, and connects it to 
existing perspectives in persistence theory. Many deeper analytic and geometric consequences 
remain to be explored.

\subsection{Outline}

This paper is organized as follows. In Section~\ref{sec:background}, we review the relevant prerequisite materials. In Section~\ref{sec:modules-mob-homology}, we introduce the central objects and constructions in this work: persistence modules and M\"obius homology.

In Section~\ref{sec:mob-laplacian} we define the M\"obius homology Laplacian in general. 

In Section~\ref{sec:the1dcase}, we focus on one-parameter persistence modules. We compute M\"obius homology spaces for such modules and completely characterize the spectrum of their associated degree-one M\"obius homology Laplacian.

In Section~\ref{sec:birth-death-laplacian}, we study the birth-death module for multiparameter filtrations. We provide a first analysis of both the kernel and the non-zero spectrum of the associated degree-zero M\"obius homology Laplacian.

%%%%%%%%%%%%%%%%%%%%%%%%%%%%%%%%%%%%%%%%%%%%%%
\subsection{Acknowledgements}
%%%%%%%%%%%%%%%%%%%%%%%%%%%%%%%%%%%%%%%%%%%%%%
This work is the result of the author's PhD dissertation, completed at Colorado State University under the supervision of Amit Patel. This work was greatly improved by conversations with and suggestions from Amit.

\section{Background}
\label{sec:background}
We begin with a brief discussion of mathematical constructions and results that play a critical role in this work. M\"obius inversion (and its role in persistent homology) will be a consistent theme throughout this work. Simplicial homology and, in particular, simplicial cosheaf homology will provide the algebraic basis needed in this work, to which we will add additional structure via inner products. The inner product structure allows for the construction of combinatorial Laplace operators: square matrices that act the chain spaces of graphs, simplicial complexes, and simplicial cosheaves.

\subsection{Posets and M\"obius Inversion}\label{sec:mob-inverison-approach}

Recall that a (finite) poset $P$ is a (finite) set together with a partial order $\leq$ that is reflexive, antisymmetric, and transitive.

\begin{defn}\label{def:mobius-inversion-form}
Let $P$ be a finite poset, and let $f: P \to \mathbb{Z}$ be an integer-valued function on $P$. Then, the \textbf{M\"obius inverse} (or, inversion) of $f$ is defined as the unique function $\partial f: P \to \mathbb{Z}$ satisfying \[
f(b) = \sum_{a \leq b} \partial f(a)
\]
for all $b \in P.$
\end{defn}

The M\"obius inverse of $f$ acts as a sort of ``combinatorial derivative." Indeed, the formula above looks like a discrete analog of the fundamental theorem of calculus. With this interpretation, if $f$ captures cumulative information on $P$, then its M\"obius inverse $\partial f$ describes local contribution.

\subsection{Hilbert Spaces and Operators}
\label{sec:hilb}

Throughout this work, we will consider finite-dimensional vector spaces over $\mathbb{R}$. We denote by $\Vect_\mathbb{R}$ the category of finite-dimensional vector spaces with morphisms given by linear maps. When equipped with an inner product, such vector spaces are \textit{Hilbert spaces}. Denote by $\mathrm{Hilb}_\mathbb{R}$ the category of finite-dimensional Hilbert spaces over $\mathbb{R}$ with morphisms given by (bounded) linear maps. \textit{From here on out, we will assume all vector spaces and Hilbert spaces are finite dimensional.}

Given a bounded linear map $f: V \to W$ between two Hilbert spaces, its \textit{adjoint} is the unique bounded linear map $f^*: W \to V$ so that $\langle w, f(v)\rangle_W = \langle f^*(w), v\rangle_V$. We call a bounded linear map $f: V \to V$ from a Hilbert space to itself an \textit{operator} on $V$. An operator satisfying $f = f^*$ is called \textit{self-adjoint}. An operator satisfying $\langle f(v), v\rangle_V \geq 0$ for all $ v\in V$ is called \textit{positive semidefinite}. Equivalently, since we are considering finite-dimensional Hilbert spaces, a self-adjoint operator is positive semidefinite if and only if all of its  eigenvalues are non-negative.

\subsection{Simplicial Homology}

A (finite) \textit{simplicial complex} $K$ over a finite vertex set $V$ is a collection of subsets of $V$, called \textit{simplices}, that is closed downwards. That is, if $\sigma \in K$ and $\tau \subseteq \sigma$, then $\tau \in K$. We call $\sigma$ an \textit{$i$-simplex} if~$|\sigma| = i+1$. Let $K^i$ be the \textit{$i$-skeleton} of $K$ (that is, all $\sigma$ such that $|\sigma| \leq i+1$), and denote the set of all $i$-simplices of $K$ by $K^{[i]} = K^i - K^{i-1}$. An orientation on a $i$-simplex $\sigma$ is given by ordering its vertices, $\sigma = [v_0 v_1 \ldots v_i]$.

For $i \geq 0$, the \textit{$i$-th chain space} $C_i(K)$ is the real vector space with basis given by the set of oriented $i$-simplices $K^{[i]}$. Define the \textit{boundary map} $\partial_i\colon C_i(K)\to C_{i-1}(K)$ on $i$-simplices as follows:
\[
\partial_i(\sigma) = \sum_{j = 0}^i (-1)^j [v_0 \ldots \hat{v}_j \ldots v_i],
\]
and extend linearly to obtain a map on all of $C_i(K)$. Here, $\hat{v}_j$ means that the $j$th vertex has been removed. 

For $i\geq 0$, define the \textit{cycle} and \textit{boundary} subspaces as
\[
Z_i(K) = \ker\partial_i \qquad\text{and}\qquad B_i(K) = \im \partial_{i+1}
\]
respectively. The \textit{$i$-th homology space} is given by the quotient $H_i(K) = Z_i(K)/B_i(K)$.

Now, consider a \textit{subcomplex} $A$ of the simplicial complex $K$. That is, $A$ is a subset of $K$ that is itself a simplicial complex. We wish to describe homology of $K$ \textit{relative} to the subcomplex $A$. For $A \subseteq K$, the \textit{relative $i$-chain space} is defined as the quotient 
\[
C_i(K, A) = \frac{C_i(K)}{C_i(A)}. 
\]
A quick check shows that the boundary maps on $K$ restrict to boundary maps on $A$ and, thus, respect the quotient. That is, we obtain maps $\partial_i: C_i(K,A) \to C_{i-1}(K, A)$.

For $i \geq 0$, define the $i$-th \textit{relative cycle and boundary spaces} as 
\[
Z_i(K, A) = \ker \partial_i \quad \text{and} \quad B_i(K, A) = \im \partial_{i+1},
\]
respectively. Then, the $i$-th \textit{relative homology space} is $H_i(K, A) = Z_i(K, A) / B_i(K, A)$.

\subsection{Simplicial Cosheaves and their Homology}
\label{subsec:sheafy}

Let $K$ be a simplicial complex. Any simplicial complex is a category, where the objects are the simplcies and there is a morphism $\sigma \to \tau$ whenever $\tau \subseteq \sigma$ (we say $\tau$ is a \textit{face} of $\sigma$). 

A \textit{simplicial cosheaf} valued in $\Vect_\mathbb{R}$ is a functor $\mathcal{F}: K \to \Vect_{\mathbb{R}}$. Intuitively, simplicial cosheaves assign vector spaces to each simplex, and the linear maps propagate information \textit{down} in dimension. That is, if $\tau$ is a face of $\sigma$, then there is a linear map $\mathcal{F}(\sigma \to \tau)$. Morphisms of simplicial cosheaves are natural transformations. 

Given a simplicial cosheaf $\mathcal{F}$, define the \textit{$i$-th chain space} as the direct sum: \[
C_i(K;\mathcal{F}) = \bigoplus_{\sigma \in K^{[i]}} \mathcal{F}(\sigma).
\]
Define the \textit{$i$-th boundary map} as the direct sum:
\[
\partial_i = \bigoplus_{\sigma \in K^{[i]}}\left( \sum_{\sigma < \tau} [\tau:\sigma]\mathcal{F} (\sigma \to \tau)  \right),
\]
where $[\tau:\sigma] = (-1)^{j}$ for $\sigma = [v_0 v_1 \cdots v_i]$ and $\tau = [v_0 \cdots \hat{v_j} \cdots v_i]$. Here, $\hat{v_j}$ means we remove the $j$th vertex.

The \textit{cosheaf homology} spaces, then, are given by
\[
H_i(K;\mathcal{F}) = \ker\partial_i/\im\partial_{i+1}.
\]

This construction generalizes simplicial homology by allowing each simplex to carry the additional structure of a vector space, where the boundary maps describe how information is related across faces. See \cite{curry2014sheaves} for more on simplicial cosheaves and cosheaf homology. 

\subsection{Combinatorial Hodge Laplacians} 
\label{subsec:combo-lap}

Combinatorial Hodge Laplacians are self-adjoint, positive semidefinite operators that act on the chain spaces of graphs, simplicial complexes, and simplicial cosheaves. 

\paragraph{Simplicial Hodge Laplacians}

A common inner product structure for $C_i(K)$ is given by declaring that the set of $i$-simplices of $K$ forms an orthonormal basis. That is, for all $\sigma, \sigma' \in K^{[i]}$, we define
\[
\langle \sigma, \sigma'\rangle_{C_i(K)} = \begin{cases}
1 & \text{if } \sigma = \sigma'\\
0 & \text{else}
\end{cases}
\]
and extended linearly. 
Thus, $C_i(K)$ is a Hilbert space and so are the cycle and boundary subspaces, $Z_i(K)$ and $B_i(K)$. Since each space is finitely generated, the boundary maps $\partial_i$ are bounded.

\begin{defn}[\cite{Eckmann1944}]
Let $K$ be a (finite) simplicial complex, and endow each chain space $C_i(K)$ with an inner product. The $i$-th \textbf{combinatorial Laplacian} is defined as:
\[
L_i = \partial_{i+1} \partial_{i+1}^* + \partial_i^* \partial_i: C_i(K) \to C_i(K).
\]
\end{defn}

Furthermore, because we endowed each chain space $C_i(K)$ with the inner product
for which the $i$-simplices form an orthonormal basis, adjoints are given by transposes. That is, the combinatorial Laplacian is written as $L_i = \partial_{i+1} \partial_{i+1}^T + \partial_i^T \partial_i$. When the simplicial complex needs emphasis, we will write $L_i^K$.

In particular, a graph $G$ can be viewed a 1-dimensional simplicial complex. When considering a graph, the $0$-th Laplacian $L_0 = \partial_1 \partial_1^T$ is called the \textit{graph Laplacian}. The graph Laplacian may also be written as the ``degree matrix" minus the``adjacency matrix." That is, we write $L_0 = D - A$ where $D$ is the diagonal matrix with entries $D(i,i)$ given by the degree of vertex $i$. Define the adjacency matrix $A$ by setting $A(i,j) = 1$ if vertex $i$ and vertex $j$ are adjacent and $A(i,j) = 0$ otherwise.

It is known that the spectrum of the combinatorial Laplacian for graphs and simplicial complexes reflect the combinatorial structure of the underlying object. See Section~\ref{sec:related-work} for a brief overview of the literature in this area. In particular, a result by Anderson and Morley states that the eigenvalues of the graph Laplacian for a graph with $n$ vertices will always lie in the (closed) interval $[0,n]$ \cite[Theorem 1]{Anderson01101985}. The proof of that result includes the following proposition, which will be of use later in this work.

\begin{prop}\label{prop:spectrum_of_L0Kn}
    Let $K_n$ be the complete graph on $n$ vertices. Then, the eigenvalues of $L_0^{K_n}$ are~0 with multiplicity 1 and $n$ with multiplicity $n-1$.
    \begin{proof}
    Since each vertex in $K_n$ has degree $n-1$, the adjacency matrix for $K_n$ is given by~$A = J - I_n$, where $J$ is the all-ones matrix and $I_n$ is the $n\times n$ identity matrix. Thus,
    \[
    L_0^{K_n} = D - A = (n-1)I_n - (J - I_n) = nI_n - J.
    \]
    Let $\vec{1}$ denote the all-ones vector. Then, $J$ 
    satisfies $J \vec{1} = n\vec{1}$ and $
    L_0^{K_n}\vec{1} = \vec{0}.$
    Moreover, $Jx = \vec{0}$ for all $x$ orthogonal to $\vec{1}$. Thus, $
    L_0^{K_n} x = n x.$ 
    Therefore, the eigenvalues of $L_0^{K_n}$ are 0 with multiplicity~1 and $n$ with multiplicity $n-1$.
    \end{proof}
\end{prop}

    \begin{rmk}\label{rmk:entries-sum-to0}
    In the proof of Proposition~\ref{prop:spectrum_of_L0Kn}, we showed that the eigenvectors $x$ of $L_0^{K_n}$ with non-zero eigenvalue are orthogonal to $\vec{1}$. Indeed, since the Laplacian is self-adjoint, eigenvectors with distinct eigenvalues must be orthogonal. In this case, $ 0 = \langle x, \vec{1} \rangle = \sum_{i = 1}^n x_i$. That is, the entries of any eigenvector $x$ for $L_0^{K_n}$ with eigenvalue $n$ must sum to 0.
    \end{rmk}

    \paragraph{Simplicial Cosheaf Laplacians}
    
    Hansen and Ghrist extend the simplicial Hodge Laplacian to include the additional data of a \textit{cosheaf} over a simplicial complex \cite{Hansen2019spectral, hansen20phd}.
    
    Let $\mathcal{F}$ be a simplicial cosheaf. Suppose, now, that each $\mathcal{F}(\sigma)$ carries the additional structure of a Hilbert space. Since we are restricting to finite simplicial complexes and finite dimensional Hilbert spaces, all chain spaces are necessarily finitely generated and all maps $\mathcal{F}(\sigma \to \tau)$ are bounded linear maps. Thus, we will now consider cosheaves valued over the category of Hilbert spaces: $\mathcal{F}: K \to \mathrm{Hilb}_{\mathbb{R}}$. In \cite{Hansen2019spectral}, the authors call these \textit{weighted simplicial cosheaves}. 

With an inner product on each costalk $\mathcal{F}(\sigma)$, we obtain an inner product on cosheaf chain spaces as follows. Each summand in $C_i(K; \mathcal{F}) = \bigoplus_{\sigma, \; \dim \sigma = i} \mathcal{F}(\sigma)$ retains its inner product structure, and all summands are mutually orthogonal. That is, $\langle\cdot, \cdot\rangle_{C_i(K; \mathcal{F})}$ is given by:
\[
\langle x, y \rangle_{C_i(K; \mathcal{F})} = \begin{cases}\langle x, y \rangle_{\mathcal{F}(\sigma)} & \text{if } x, y \in \mathcal{F}(\sigma), \\
0 & \text{if } x \in \mathcal{F}(\sigma), y \in \mathcal{F}(\sigma') \text{ for } \sigma \neq \sigma'.
\end{cases}
\]

Whereas cosheaf homology itself is an algebraic operation, the addition structure of inner products on each costalk $\mathcal{F}(\sigma)$ and, thus, each chain space $C_i(K; \mathcal{F})$ will allow us to investigate \textit{analytic} properties of these cosheaves. We do so by extending the definition of the combinatorial Hodge Laplacian to simplicial cosheaves.

\begin{defn}[{\cite{Hansen2019spectral}}]\label{def:cosheaf-Laplacian}
Given a simplicial cosheaf $\mathcal{F}: K \to \mathrm{Hilb}_{\mathbb{R}}$, the \textbf{$i$-th cosheaf Laplacian} is the operator on $C_i(K; \mathcal{F})$ defined as follows:
\[
L_i^\mathcal{F} = \partial_{i+1} \partial_{i+1}^* + \partial_i^* \partial_i.
\]
\end{defn}
Although the definition above looks identical to that of the combinatorial Laplacian, it is important to remember that we have the extra information of a cosheaf $\mathcal{F}$ over $K$.

%%%%%%%%%%%%%%%%%%%%%%%%%%%%%%%%%%%%%%%%%%%%%%
%%%%%%%%%%%%%%%%%%%%%%%%%%%%%%%%%%%%%%%%%%%%%%
\section{Persistence Modules and M\"obius Homology}
\label{sec:modules-mob-homology}
%%%%%%%%%%%%%%%%%%%%%%%%%%%%%%%%%%%%%%%%%%%%%%
%%%%%%%%%%%%%%%%%%%%%%%%%%%%%%%%%%%%%%%%%%%%%%

In this work, we will use a cellular cosheaf Laplace operator in order to obtain analytical information for \textit{persistence modules} using the theory of \textit{M\"obius homology}.

\subsection{Persistence Modules}
\label{subsec:modules}

We now define persistence modules, the central objects of study in this work. Indeed, these are the mathematical objects for which we will develop a spectral theory.

\begin{defn}\label{def:p-mod}
Let $P$ be a poset. A \textbf{persistence module} is a functor $M: P \to \Vect_\mathbb{R}$.
\end{defn}

That is, persistence modules map poset elements $p \in P$ to vector spaces $M(p)$ and poset relations $p \leq q$ to linear maps $M(p \leq q)$ so that everything commutes. Morphisms between two persistence modules $M, N: P \to \Vect_\mathbb{R}$, then, are natural transformations $\varphi: M \Rightarrow N$. Two modules are isomorphic if there is a natural isomorphism between them.

Given two persistence modules $M: P \to \Vect_\mathbb{R}$ and $N: P \to \Vect_\mathbb{R}$, their direct sum is another persistence module $M \oplus N: P \to \Vect_\mathbb{R}$ defined as follows: for poset elements $a \in P$, $(M \oplus N )(a) = M(a) \oplus N(a)$. For poset relations, $(M \oplus N )(a \leq b) = M(a \leq b) \oplus N (a \leq b)$. 

For the remainder of this paper, we require the following finiteness conditions: first, we will restrict to \textit{finite} posets $P$. As mentioned previously, we also require that each $M(p)$ be finite-dimensional.

In persistence theory, we aim to answer the following question: when do features first appear (are ``born") in a module $M$, and when are they mapped to zero (``die")? 

Observe that features may persist indefinitely rather than disappearing at a finite parameter value. In order to capture such features alongside those that die at a finite index, it is standard to adjoin a formal top element to the finite poset $P$ representing ``infinity." 

\begin{defn}\label{def:extended-poset}
Let $P$ be a finite poset. We define $$\overline{P} = P \cup \{ \infty\}$$
to be the poset obtained by adjoining a top element $\infty$ with the convention that $a < \infty$ for all elements $a \in P$.
\end{defn}

Thus, we will use elements of $P$ to identify the birth of a feature in $M: P \to \Vect_\mathbb{R}$, and we will use elements in $\overline{P}$ to identify its death. We represent the birth and the death of a feature together using \textit{symbolic persistence intervals.} 

\begin{defn}\label{def:halfopen-ints} 
Let $P$ be a finite poset and let $\overline{P} = P \cup \{ \infty\}$. A half-open \textbf{persistence interval} is a formal symbol $[a,b)$ with $a \in P$ and $b \in \overline{P}$ such that $a \leq b$. Denote by $\Dgm(P)$ the set of all persistence intervals. The subset $\{[a,a) \; | \; a \in P\}$ is called the \textbf{diagonal}.
\end{defn}

We emphasize that $[a, b)$ is not interpreted as a subset of $P$, but rather as a symbolic object recording a birth index $a$ and a death index $b$. 

For modules indexed over a \textit{totally ordered} set $T$, which we call one-parameter modules, the appearance and disappearance of features throughout the module correspond with its algebraic decomposition into ``interval modules." This decomposition is unique up to isomorphism. We defer the statement of this structural theorem to Section~\ref{sec:the1dcase}. When considering modules over a finite poset $P$ (not necessarily totally ordered), there exists a similar decomposition into indecomposable summands, but finding such a decomposition is often considered infeasible \cite{bauer2025multiparameter}. 

M\"obius inversion provides another approach for describing the appearance and disappearance of features in a module. Although the M\"obius inversion approach to persistence was present in original papers such as \cite{ELZ02}, Patel was the first to make concrete connections between persistent homology and Rota's work on algebraic combinatorics \cite{Patel2018, Rot64}. Loosely speaking, the pipeline is as follows. Given a persistence module $M: P \to \Vect_\mathbb{R}$, define a suitable function $f: \Dgm (P) \to \mathbb{Z}$ obtained from the data of the module $M$. The persistence diagram of $M$, then, is the M\"obius inversion $\partial f$. 

Two common functions with which to study persistence are the rank function \cite{kim2018generalized, clause2024generalizedrankinvariantmobius, carlsson09} and the birth-death function \cite{McCPa20-edit-distance, GulenMcCleary}. In the one-parameter case, the persistence diagrams obtained from the rank and the birth-death functions are equivalent to one another, and they are also equivalent to the interval decomposition approach. However, no complete discrete invariant exists for multiparameter modules $M: P \to \Vect_\mathbb{R}$ indexed over a finite poset $P$ (not necessarily totally ordered) \cite{CaZo09}. Furthermore, the rank and the birth-death functions of the same multiparameter module $M$ may result in different persistence diagrams.  Nonetheless, M\"obius inversion is commonly used in the literature to describe the persistence of features in $M$ \cite{Patel2018, McCPa20-edit-distance, kim2018generalized}. We provide more details on the M\"obius inversion approach to persistence in Section~\ref{sec:filtrations}. 

M\"obius inversion and persistence diagrams provide integer-valued data. In order to develop a spectral theory for persistence modules, we need higher order structure with which to analyze a module. Patel and Skraba's M\"obius homology provides such an approach \cite{patel2023mobius}.

\subsection{M\"obius Homology}
\label{sec:mob-homology-defn}

M\"obius homology categorifies M\"obius inversion. It assigns a chain complex to each poset element $c \in P$ whereas M\"obius inversion assigns only an integer. Loosely speaking, given a module $M: P \to \Vect_\mathbb{R}$, M\"obius homology is the homology of a cosheaf on the order complex of $P$ that is constructed using the data of the persistence module $M$.

First, we define the order complex of a poset $P$. This is a method for obtaining a simplicial complex associated to a poset. With the structure of a simplicial complex, we can use tools from algebraic topology (namely, homology and simplicial cosheaves) to study the structure of $P$. 

\begin{defn}
Let $P$ be a poset. Its \textbf{order complex}, denoted $\Delta P$, is the simplicial complex constructed by including an $i$-simplex for all $i$-chains $a_0 < a_1 < \cdots < a_i$ in $P$. Since subchains are necessarily chains themselves, $\Delta P$ is closed downwards and is, indeed, a simplicial complex.
\end{defn}

\begin{figure}
% !TEX root = main.tex

\centering

\tikzset{every picture/.style={line width=0.75pt}} %set default line width to 0.75pt        

\begin{tikzpicture}[x=0.75pt,y=0.75pt,yscale=-1,xscale=1]
%uncomment if require: \path (0,132); %set diagram left start at 0, and has height of 132

%Straight Lines [id:da7601651934470641] 
\draw    (233,31.44) -- (258,59.44) ;
%Straight Lines [id:da4238914396172333] 
\draw    (198,59.44) -- (221,31.44) ;
%Straight Lines [id:da2620891244261778] 
\draw    (198,77.44) -- (223,105.44) ;
%Straight Lines [id:da32617217121231157] 
\draw    (234,105.44) -- (257,77.44) ;
%Shape: Right Triangle [id:dp6066332272388044] 
\draw  [fill={rgb, 255:red, 155; green, 155; blue, 155 }  ,fill opacity=1 ] (387,43.59) -- (457,83.59) -- (387,83.59) -- cycle ;
%Shape: Right Triangle [id:dp13293097501494877] 
\draw  [fill={rgb, 255:red, 155; green, 155; blue, 155 }  ,fill opacity=1 ] (457,83.59) -- (387,43.59) -- (457,43.59) -- cycle ;
%Shape: Circle [id:dp5112870017337332] 
\draw  [fill={rgb, 255:red, 0; green, 0; blue, 0 }  ,fill opacity=1 ] (453.59,83.59) .. controls (453.59,81.71) and (455.12,80.19) .. (457,80.19) .. controls (458.88,80.19) and (460.41,81.71) .. (460.41,83.59) .. controls (460.41,85.47) and (458.88,87) .. (457,87) .. controls (455.12,87) and (453.59,85.47) .. (453.59,83.59) -- cycle ;
%Shape: Circle [id:dp855889139696355] 
\draw  [fill={rgb, 255:red, 0; green, 0; blue, 0 }  ,fill opacity=1 ] (453.59,43.59) .. controls (453.59,41.71) and (455.12,40.19) .. (457,40.19) .. controls (458.88,40.19) and (460.41,41.71) .. (460.41,43.59) .. controls (460.41,45.47) and (458.88,47) .. (457,47) .. controls (455.12,47) and (453.59,45.47) .. (453.59,43.59) -- cycle ;
%Shape: Circle [id:dp6722359610425004] 
\draw  [fill={rgb, 255:red, 0; green, 0; blue, 0 }  ,fill opacity=1 ] (383.59,43.59) .. controls (383.59,41.71) and (385.12,40.19) .. (387,40.19) .. controls (388.88,40.19) and (390.41,41.71) .. (390.41,43.59) .. controls (390.41,45.47) and (388.88,47) .. (387,47) .. controls (385.12,47) and (383.59,45.47) .. (383.59,43.59) -- cycle ;
%Shape: Circle [id:dp10307412475049493] 
\draw  [fill={rgb, 255:red, 0; green, 0; blue, 0 }  ,fill opacity=1 ] (383.59,83.59) .. controls (383.59,81.71) and (385.12,80.19) .. (387,80.19) .. controls (388.88,80.19) and (390.41,81.71) .. (390.41,83.59) .. controls (390.41,85.47) and (388.88,87) .. (387,87) .. controls (385.12,87) and (383.59,85.47) .. (383.59,83.59) -- cycle ;

% Text Node
\draw (223,106.84) node [anchor=north west][inner sep=0.75pt]    {$a$};
% Text Node
\draw (187,59.84) node [anchor=north west][inner sep=0.75pt]    {$b$};
% Text Node
\draw (257,62.84) node [anchor=north west][inner sep=0.75pt]    {$c$};
% Text Node
\draw (224,14.84) node [anchor=north west][inner sep=0.75pt]    {$d$};
% Text Node
\draw (373,26.84) node [anchor=north west][inner sep=0.75pt]    {$a$};
% Text Node
\draw (462.41,23.84) node [anchor=north west][inner sep=0.75pt]    {$b$};
% Text Node
\draw (462.41,82.99) node [anchor=north west][inner sep=0.75pt]    {$d$};
% Text Node
\draw (368,86.4) node [anchor=north west][inner sep=0.75pt]    {$c$};

\end{tikzpicture}
\vspace{-0.2in}
\caption{The Hasse diagram of a poset $P$ (left) and its order complex $\Delta P$ (right).}
\label{fig:ordercomplex}
\end{figure}
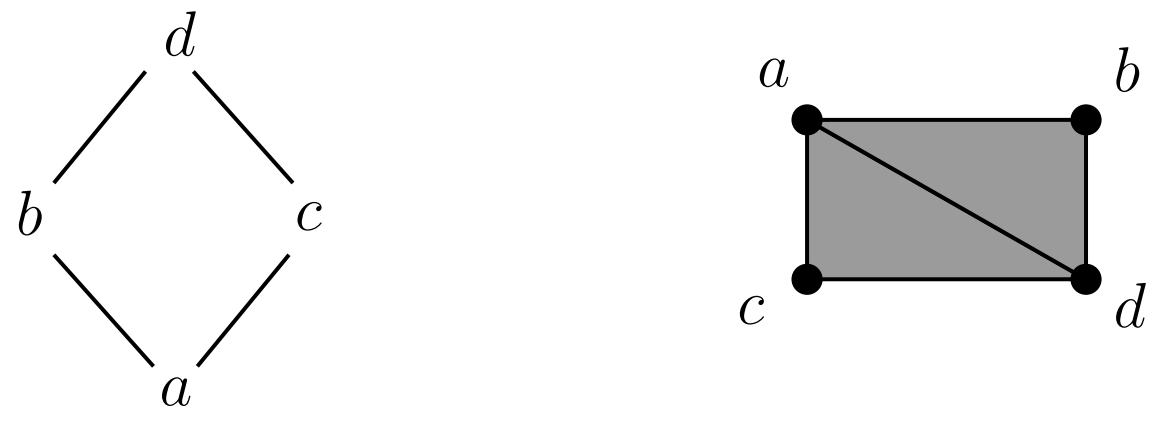

Intuitively, the order complex captures the combinatorial ``shape" of a poset by encoding the elements of $P$ along with comparability and sequences of comparable elements (that is, longer chains). See Figure~\ref{fig:ordercomplex} for an example of an order complex.

Observe that, when $P$ is a finite poset, $\Delta P$ is a finite simplicial complex. Now, given a persistence module $M: P \to \Vect_\mathbb{R}$ over such a finite poset $P$, we will construct a cosheaf on the order complex $\Delta P$ using the data of $M$.

\begin{defn}[{\cite[Definition 3.6]{patel2023mobius}}]\label{def:order_complex}
Let $M: P \to \Vect_\mathbb{R}$ be a module over a finite poset $P$. The \textbf{order cosheaf} of $M$ is the simplicial cosheaf $\underline{M}: \Delta P \to \Vect_\mathbb{R}$ defined as follows. Given a simplex $\sigma \in \Delta P$, let $\underline{M}(\sigma) = M(\min \sigma)$. For $\tau \subseteq \sigma$, the morphism $\underline{M}(\sigma \to \tau)$ is given by the linear map $M(\min\sigma \leq \min \tau)$.
\end{defn}

We will now focus on a specific element $c$ in our poset $P$ by considering its downset: that is, all elements $x$ such that $x \leq c$. In the order complex, a similar construction gives way to a subcomplex of $\Delta P$.

\begin{defn}[{\cite[Definition 3.7]{patel2023mobius}}]\label{def:lower_complex}
The \textbf{lower complex} of $c \in P$ is the subcomplex of $\Delta P$ given by $\Delta P_{\leq c} = \{\sigma \in \Delta P \; | \; \max \sigma \leq c\}$. The \textbf{strict lower complex} of $c \in P$ is the subcomplex given by $\Delta P_{< c} = \{\sigma \in \Delta P \; | \; \max \sigma < c\}$.
\end{defn}

That is, the lower complex of $\Delta P$ at $c$ contains only the chains that end at or before $c$. The \textit{strict} lower complex of $\Delta P$ at $c$ contains only the chains that end (strictly) before $c$. We are specifically interested in isolating chains that end \textit{exactly} at $c$: thus, we compute relative chain spaces.

\begin{defn}[{\cite[Definition 3.8]{patel2023mobius}}]\label{def:mobius_chain_complex}
Let $P$ be a finite poset, and let $M: P \to \Vect_\mathbb{R}$ be a module. Fix $c \in P$. The \textbf{Möbius chain complex} of $M$ at $c$ is the relative simplicial cosheaf chain complex given by
\[
C^\downarrow_\bullet M(c) = C_\bullet(\Delta P_{\le c}, \Delta P_{< c}; M).
\]
Denote by 
\[
\partial_i^\downarrow: C_i^\downarrow M(c) \to C_{i-1}^\downarrow M(c)
\]
the boundary maps in $C^\downarrow_\bullet M(c)$.

Applying homology in dimension $i$ yields the \textbf{M\"obius homology space} of $M$ at $c$: 
\[H_i^\downarrow M(c) = \ker \partial_i^\downarrow / \im \partial_{i+1}^\downarrow. \]
\end{defn}

That is, M\"obius homology assigns a \textit{vector space} to every $i \geq 0$ and element $c \in P$, retaining chain-level structure that is lost when considering integer-valued invariants. Nonetheless, when ``decategorifying," we relate M\"obius homology to M\"obius inversion via the following proposition.

\begin{prop}[{\cite[Theorem 3.13]{patel2023mobius}}]\label{prop:mob-hom-categorify}
Let $P$ be a finite poset, let $M: P \to \Vect_\mathbb{R}$ be a module, and let $m: P \to \mathbb{Z}$ be the dimension function of $M$. That is, $m(c) = \dim M(c)$. Then, for all~$c \in P$,
\begin{equation}\label{eq:mob-eulerchar}
\partial m(c) = \sum_{i \geq 0} (-1)^i \dim H_i^\downarrow M(c).
\end{equation}
\end{prop}

The proposition above states precisely how M\"obius homology is a categorification of M\"obius inversion. Indeed, the M\"obius inversion of the dimension function  is the \textit{Euler characteristic} of the order cosheaf $\underline{M}$. See \cite{patel2023mobius} for examples and further discussion of using M\"obius homology to distinguish modules that are indistinguishable by integer-valued invariants (such as dimension functions or persistence diagrams) alone. 

When the vector spaces mapped to via a module $M$ carry the additional structure of an inner product, we obtain an inner product structure on the M\"obius chain spaces associated to $M$. In particular, the boundary operators in the M\"obius chain complex admit adjoint operators. This, ultimately, allows for the construction of a cellular cosheaf Laplace operator on the M\"obius chain complex, which we call the \textit{M\"obius homology Laplacian}.

%%%%%%%%%%%%%%%%%%%%%%%%%%%%%%%%%%%%%%%%%%%%%%
%%%%%%%%%%%%%%%%%%%%%%%%%%%%%%%%%%%%%%%%%%%%%%
\section{The M\"obius Homology Laplacian}
\label{sec:mob-laplacian}
%%%%%%%%%%%%%%%%%%%%%%%%%%%%%%%%%%%%%%%%%%%%%%
%%%%%%%%%%%%%%%%%%%%%%%%%%%%%%%%%%%%%%%%%%%%%%

For the remainder of this paper, we will consider persistence modules valued in the category of Hilbert spaces: $M: P \to \Hilb_\mathbb{R}$. We again require that each $M(c)$ is finite dimensional. 

Recall from Section~\ref{subsec:combo-lap} that the M\"obius chain space $$C_i^\downarrow M(c) = \bigoplus_{i\text{-chains } \sigma, \; \max(\sigma) = c} M(\min \sigma)$$ inherits an inner product structure as follows. In $C_i^\downarrow M(c)$, each summand $M(\min \sigma)$ retains its inner product structure and is orthogonal to every $M(\min \sigma')$ for distinct $i$-chains $\sigma' \neq \sigma$. With an inner product structure on each M\"obius chain space, we can define adjoint operators for each boundary map $\partial_i^\downarrow$. Thus, assuming that all of the Hilbert spaces $M(c)$ are finite-dimensional, we construct the cellular cosheaf Laplacian on the M\"obius chain space by way of \cite{Hansen2019spectral}.

\begin{defn}\label{def:mob-laplacian}
Let $M: P \to \Hilb_\mathbb{R}$ be a module over a finite poset $P$ whose stalks are finite-dimensional. Then, the degree-$i$ \textbf{M\"obius homology Laplacian}
of $M$ at $c \in P$ is the operator
\begin{align*}
    L_i^\downarrow M(c) =   \partial_{i+1}^\downarrow(\partial_{i+1}^\downarrow)^* +  (\partial_i^\downarrow)^*\partial_i^\downarrow 
\end{align*}
acting on the M\"obius chain space $C_i^\downarrow M(c)$.
Recall that $\partial_i^\downarrow$ denotes the $i$-th boundary operator in the M\"obius chain complex $C_\bullet^\downarrow M(c)$.
\end{defn}

It is worth highlighting that the M\"obius homology Laplacian is defined point-wise at a fixed index $c \in P$. Thus, provided a module $M: P \to \Hilb_\mathbb{R}$, we may construct a Laplacian $L_i^\downarrow M(c)$ for \textit{any} poset element $c \in P$ and for \textit{any} dimension $i \geq 0$. 

For every $c \in P$ and $i \geq 0$, a quick check shows that $L_i^\downarrow M(c)$ is self-adjoint and positive semidefinite. We state other properties of $L_i^\downarrow M(c)$ below. Although these properties are important enough to state clearly, they are standard in discrete Hodge theory. Thus, we defer the proofs to Appendix~\ref{appendix:proofs} or provide relevant citations.

First, we immediately obtain the following result by applying \cite[Theorem 3.1]{Hansen2019spectral} to the M\"obius homology construction.

\begin{thm}\label{prop:harmonic_cosheaf_reps} Let $P$ be a finite poset, let $c \in P$, and let $M:P \to \Hilb_\mathbb{R}$ be a module. Then, 
\[\ker L_i^\downarrow M(c) \cong H_i^\downarrow M(c).\]
\begin{proof}
    See \cite[Theorem 3.1]{Hansen2019spectral}.
\end{proof}
\end{thm}

Indeed, in order to prove the statement above, Hansen and Ghrist prove the equality $\ker L_i^\downarrow M(c) = \ker \partial_i^\downarrow \cap \im (\partial_{i+1}^\downarrow)^\perp.$ Then, the orthogonal projection map $\mathrm{proj}_{(\im \partial_{i+1}^\downarrow)^\perp}(-)$  yields the desired isomorphism.

Furthermore, as in simplicial harmonic homology, $\ker L_i^\downarrow M(c)$ provides \textit{energy minimizing cycle representatives} for $H_i^\downarrow M(c)$, which on its own provides only an equivalence class of cycles. 

\begin{defn}\label{def:space-of-reps} Let $P$ be a finite poset, let $c \in P$, and let $M:P \to \Hilb_\mathbb{R}$ be a module.  
Denote by \[\varphi: \ker \partial_i^\downarrow \to H_i^\downarrow M(c)\] the natural surjection obtained by taking the quotient $\ker \partial_i^\downarrow / \im \partial_{i+1}^\downarrow$. The space of \textbf{cycle representatives} for a homology class $[x] \in H_i^\downarrow M(c)$ is given by $\varphi^{-1}([x])$.
\end{defn}

Obtaining cycle representatives for $H_i^\downarrow M(c)$ is akin to finding a suitable section of $\varphi$. In the following theorem, we show that $\ker L_i^\downarrow M(c)$ provides a method for obtaining such cycle representatives.

\begin{thm}\label{prop:projection-section} Let $P$ be a finite poset, let $c \in P$, and let $M:P \to \Hilb_\mathbb{R}$ be a persistence module.
The inclusion $s: \ker L_i^\downarrow M(c) \hookrightarrow \ker \partial_i^\downarrow$ induces a section of the quotient map $\varphi: \ker \partial_i^\downarrow \to H_i^\downarrow M(c)$.
\begin{proof}
See Appendix~\ref{appendix:firstproof}.
\end{proof}
\end{thm}

In other words, Theorem~\ref{prop:projection-section} states that for each homology class $[x] \in H_i^\downarrow M(c)$, the element $\mathrm{proj}_{(\mathrm{im}\partial_{i+1}^\downarrow)^\perp}(x)$ is the \textit{unique representative} of $[x]$ lying in $\ker L_i^\downarrow M(c)$. Moreover, these representatives have nice geometric properties. 

\begin{defn}
Let $x$ be a vector in an inner product space $V$. The \textbf{energy} of $x$ is given by its norm squared $\|x\|_V^2 = \langle x,x\rangle_V$.
\end{defn}

\begin{thm}\label{prop:minimal-energy}
Let $P$ be a finite poset, let $c \in P$, and let $M:P \to \Hilb_\mathbb{R}$ be a module. Let $x \in \ker L_i^\downarrow M(c)$ and let $y = x + \partial_{i+1}^\downarrow z$ be another cycle representative of the same homology class. Then,
\[
\| x\|^2 \leq \|y \|^2.
\]
That is, $x \in \ker L_i^\downarrow M(c)$ is the unique representative with minimal energy.
\begin{proof}
    See Appendix~\ref{appendix:secondproof}
\end{proof}
\end{thm}

That is, the kernel of the M\"obius homology Laplacian provides canonical and energy minimizing representatives for M\"obius homology spaces. In the sections to follow, we will utilize these results for modules arising from persistent homology. Furthermore, for these modules specifically, the inherent structure of M\"obius homology makes its associated cosheaf Laplacian amenable to spectral analysis. 

To illustrate the claim above, we will investigate the M\"obius homology Laplacian in-depth for persistence modules with well-known structure. Indeed, in the next section, the well-known interval decomposition of one-parameter modules will allow for not only an explicit computation of M\"obius homology, but also a complete characterization of the spectrum of its associated degree-one M\"obius homology Laplacian. In later sections, will we focus on modules arising from simplicial filtrations. For these modules, the structure of M\"obius homology also makes the associated M\"obius homology Laplacian amenable to spectral analysis.

\section{One-Parameter Persistence Modules}
\label{sec:the1dcase}

We begin our study of the M\"obius homology Laplacian by focusing on one-parameter modules: that is, modules $M: T \to \Vect_\mathbb{R}$ for $T$ a finite totally ordered set. One-parameter modules have a unique (up to isomorphism) decomposition into interval modules. The goal of this section is to show that both M\"obius homology and the associated M\"obius homology Laplacian identify the already-known structure of one-parameter modules. Thus, we begin by reviewing the structural theorem for one-parameter persistence modules.

\begin{defn}\label{def:int-mod}
Let $T$ be a finite totally ordered set and let $[a,b) \in \Dgm (T)$ with $a < b$. The \textbf{interval module} $I_{[a,b)}: T \to \Vect_\mathbb{R}$ is defined as follows:
\[
I_{[a,b)}(c) = \begin{cases}
\mathbb{R}, & \text{if } a \leq c < b \\
0, & \text{otherwise.}
\end{cases}
\]
For $c \leq d$ in $T$, the morphism
\[
I_{[a,b)}(c \leq d): I_{[a,b)}(c) \to I_{[a,b)}(d)
\]
is the identity whenever $I_{[a,b)}(c)= I_{[a,b)}(d) = \mathbb{R}$ and is the zero map otherwise.
\end{defn}

The following proposition utilizes results from representation theory and the structural theorem of finitely generated modules over principal ideal domains. We state it using the notation from persistent homology originating with Zomorodian and Carlsson's work.

\begin{prop}[{\cite[Section 3.2]{ZoCa05}}]\label{def:interval-decomp}
Let $T$ be a finite totally ordered set and $M: T \to \Vect_\mathbb{R}$ be a persistence module such that each $M(c)$ is finitely generated. Then, there exists a function
\[
m_\mathcal{B}: \Dgm (T) \to \mathbb{Z}_{\geq 0}
\]
such that 
\[
M \cong \bigoplus_{[a,b) \in \Dgm (T)} I_{[a,b)}^{\oplus m_\mathcal{B}([a,b))}.
\]
The function $m_\mathcal{B}$ (equivalently, the associated multiset of intervals) is uniquely determined by $M$. We call $\bigoplus_{[a,b) \in \Dgm (T)} I_{[a,b)}^{\oplus m_\mathcal{B}([a,b))}$ the \textbf{interval decomposition of $M$}.
\end{prop}

In particular, because the function $m_\mathcal{B}$ is uniquely determined by $M$, the interval decomposition of a one-parameter module $M: T \to \Vect_\mathbb{R}$ determines its isomorphism type.

\begin{rmk}\label{remark:nonunique-iso}
Although the decomposition of a module $M$ given in Proposition~\ref{def:interval-decomp} is unique, the \textit{isomorphism} $\varphi: M \Rightarrow\bigoplus_{[a,b) \in \Dgm (T)} I_{[a,b)}^{\oplus m_\mathcal{B}([a,b))}$ is not.
\end{rmk}

\begin{defn}\label{def:bar} Let $T$ be a finite totally ordered set and let $M: T \to \Vect_\mathbb{R}$ be a persistence module. The \textbf{barcode} of $M$ is the finite-support function 
\[
m_\mathcal{B}: \Dgm (T) \to \mathbb{Z}_{\geq 0}
\]
such that 
\[M \cong \bigoplus_{[a,b) \in \Dgm (T)} I_{[a,b)}^{\oplus m_\mathcal{B}([a,b))}\]
as in Proposition~\ref{def:interval-decomp}. We call an interval $[a,b)$ a \textbf{bar} in the barcode if $m_\mathcal{B}([a,b)) \neq 0$, and the value $m_\mathcal{B}([a,b))$ is called the \textbf{multiplicity} of the bar $[a,b)$.
\end{defn}

\begin{rmk}
Observe that we excluded intervals of the form $[a,a)$ from the definition of an interval module (Definition~\ref{def:int-mod}). Thus, diagonal intervals $[a,a)$ necessarily have multiplicity zero in the barcode for any module $M$.
\end{rmk}

\begin{figure}
% !TEX root = main.tex

\centering

\tikzset{every picture/.style={line width=0.75pt}} %set default line width to 0.75pt        

\begin{tikzpicture}[x=0.75pt,y=0.75pt,yscale=-1,xscale=1]
%uncomment if require: \path (0,81); %set diagram left start at 0, and has height of 81

%Straight Lines [id:da5362083836831932] 
\draw    (225,46.6) -- (434.5,47.43) ;
\draw [shift={(436.5,47.44)}, rotate = 180.23] [color={rgb, 255:red, 0; green, 0; blue, 0 }  ][line width=0.75]    (10.93,-3.29) .. controls (6.95,-1.4) and (3.31,-0.3) .. (0,0) .. controls (3.31,0.3) and (6.95,1.4) .. (10.93,3.29)   ;
%Straight Lines [id:da37257272185897805] 
\draw [line width=4.5]    (257,34.6) -- (351.5,34.44) ;
%Straight Lines [id:da40219505914538956] 
\draw [line width=4.5]    (305,20.6) -- (399.5,20.44) ;
%Straight Lines [id:da1235143200103348] 
\draw [line width=1.5]    (257,34.6) -- (351.5,34.44) ;
\draw [shift={(351.5,34.44)}, rotate = 179.9] [color={rgb, 255:red, 0; green, 0; blue, 0 }  ][line width=1.5]      (6.71,-6.71) .. controls (3.01,-6.71) and (0,-3.7) .. (0,0) .. controls (0,3.7) and (3.01,6.71) .. (6.71,6.71) ;
\draw [shift={(257,34.6)}, rotate = 359.9] [color={rgb, 255:red, 0; green, 0; blue, 0 }  ][fill={rgb, 255:red, 0; green, 0; blue, 0 }  ][line width=1.5]      (0, 0) circle [x radius= 4.36, y radius= 4.36]   ;
%Straight Lines [id:da7112483435307635] 
\draw [line width=1.5]    (305,20.6) -- (399.5,20.44) ;
\draw [shift={(399.5,20.44)}, rotate = 179.9] [color={rgb, 255:red, 0; green, 0; blue, 0 }  ][line width=1.5]      (6.71,-6.71) .. controls (3.01,-6.71) and (0,-3.7) .. (0,0) .. controls (0,3.7) and (3.01,6.71) .. (6.71,6.71) ;
\draw [shift={(305,20.6)}, rotate = 359.9] [color={rgb, 255:red, 0; green, 0; blue, 0 }  ][fill={rgb, 255:red, 0; green, 0; blue, 0 }  ][line width=1.5]      (0, 0) circle [x radius= 4.36, y radius= 4.36]   ;

% Text Node
\draw (251,50.4) node [anchor=north west][inner sep=0.75pt]    {$1$};
% Text Node
\draw (300,50.4) node [anchor=north west][inner sep=0.75pt]    {$2$};
% Text Node
\draw (348,50.4) node [anchor=north west][inner sep=0.75pt]    {$3$};
% Text Node
\draw (395,50.4) node [anchor=north west][inner sep=0.75pt]    {$4$};

\end{tikzpicture}
\vspace{-0.2in}
\caption{The barcode $m_\mathcal{B}$ for a module with interval decomposition $I_{[1,3)} \oplus I_{[2,4)}$. We represent bars $[a,b)$ as half-open line segments from $a$ to $b$.}
\label{fig:simple-filt-barcode}
\end{figure}
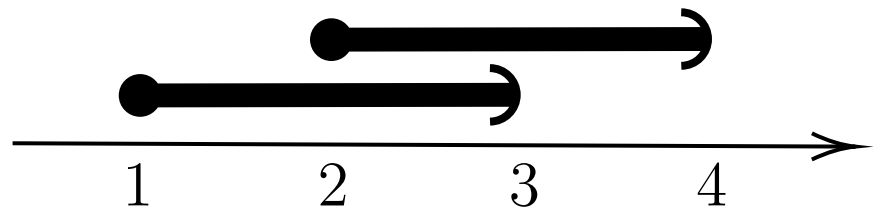

\begin{ex}\label{ex:simple-fil-barcode}
Let $M: T \to \Vect_\mathbb{R}$ be the module over the totally ordered set $1 < 2 < 3 <4 $ given by:
$$\mathbb{R} \xrightarrow{\begin{pmatrix}1 \\0\end{pmatrix}} \mathbb{R}^2 \xrightarrow{\begin{pmatrix}0 & 1\end{pmatrix}} \mathbb{R} \xrightarrow{\begin{pmatrix}0 \end{pmatrix}} 0.$$
The interval decomposition of $M$ is given by $M \cong I_{[1,3)} \oplus I_{[2,4)}$. See Figure~\ref{fig:simple-filt-barcode} for a visualization of the barcode. 
\end{ex}

We expect that M\"obius homology and the M\"obius homology Laplacian will reflect the well-known structure of one-parameter modules $M: T \to \Vect_\mathbb{R}$. Indeed, in the sections to follow, we will see that M\"obius homology identifies the birth indices and (finite) death indices of a one-parameter modules. Moreover, the spectrum of the M\"obius homology Laplacian reflects the ``length" and ``location" of bars in the barcode of $M$.

\subsection{M\"obius Homology of One-Parameter Modules}
\label{sec:mob-homology-calculation}

In this section, we will compute M\"obius homology of a module $M: T \to \Vect_\mathbb{R}$ for every element in $c \in T$. Motivated by the structure theorem for one-parameter modules, we start by computing M\"obius homology for a single interval module. 

\begin{lem}\label{lem:mobius_interval_module}
Let $I_{[a,b)}: T \to \Vect_{\mathbb{R}}$ be an interval module (see Definition~\ref{def:int-mod}). Let $c \in T$. Then,
\begin{align*}
H_0^\downarrow I_{[a,b)}(c) &\cong \begin{cases}
\mathbb{R} & \text{if } c = a,\\
0 & \text{else}.
\end{cases}\\
H_1^\downarrow I_{[a,b)}(c) &\cong \begin{cases}
\mathbb{R} & \text{if } c = b,\\
0 & \text{else}.
\end{cases}\\
H_i^\downarrow I_{[a,b)}(c) &= 0 \text{ for all } i\geq 0.
\end{align*}
\begin{proof}
Denote by $[a,b^-]$ the subposet of $T$ consisting of all poset elements $x$ such that $a \leq x < b$. That is, $[a,b^-]= \{x \in T \; | \; a \leq x < b\}$. Recall that interval modules are not defined for diagonal persistence intervals $[a,a)$. 

Fix $c \in T$. We will consider three cases: $c = a$, $c = b$, and $c \neq a, b$.

First, let $c = a$. As there are no elements $z < a \in T$ where $I_{[a,b)}(z) \neq 0$, the M\"obius chain complex $C_\bullet^\downarrow I_{[a,b)}(a)$ is given by
\[
0 \xrightarrow{} \mathbb{R} \xrightarrow{} 0.
\]
Thus, $H_0^\downarrow I_{[a,b)}(a) \cong \mathbb{R}$ and $H_i^\downarrow I_{[a,b)}(a) = 0$ for all $i \geq 1$. 

Next, let $c = b$. Since $c \in T$, $b \neq \infty$ here. Suppose $n = \# \{x \in T \; | \; a \leq x < b\}$. First, recall that~$I_{[a,b)}(x) = \mathbb{R}$ for all $a \leq x < b$. Then, there are ${n}\choose{k}$ $k$-chains in $\Delta T^\downarrow_\leq(b) \backslash \Delta T^\downarrow_<(b)$ that start at or after $a$. That is, the M\"obius chain complex $C_\bullet^\downarrow I_{[a,b)}(b)$ is given by
\[
0 \xrightarrow{} \mathbb{R}^{{n}\choose{n}}\xrightarrow{\partial_{n}^\downarrow} \cdots \xrightarrow{} \mathbb{R}^{{n}\choose{2}}\xrightarrow{\partial_{2}^\downarrow} \mathbb{R}^{{n}\choose{1}}\xrightarrow{\partial_{1}^\downarrow} 0 \xrightarrow{} 0.
\]
First, it is clear that $H_0^\downarrow I_{[a,b)}(b) = 0$. For $i \geq 0$,
we claim that $H_{i+1}^\downarrow I_{[a,b)}(b) \cong H_i (\Delta^{n-1})$ where~$\Delta^{n-1} $ denotes the standard $(n-1)$-simplex. Indeed, view the elements $x \in [a,b^-]$ as vertices with an ordering inherited from $T$, and consider all chains in $[a,b^-]$. This forms the standard $(n-1)$-simplex, and the faces of dimension $i$ are in one-to-one correspondence with the $i+1$ chains in $\Delta T^\downarrow_\leq(b) \backslash \Delta T^\downarrow_<(b)$. The boundary map $\partial_i^{\Delta^{n-1}}$ is identical to $\partial^\downarrow_{i+1}$ (the boundary map in the M\"obius chain complex). That is, $H_{i+1}^\downarrow I_{[a,b)}(b) \cong H_i (\Delta^{n-1})$. Since $\Delta^{n-1}$ is contractible, we obtain $H_1^\downarrow I_{[a,b)} (b)\cong \mathbb{R}$ and $H_i^\downarrow I_{[a,b)}(b) =0$ for $i \geq 2$ for $b \neq \infty$.

Now, suppose $c \neq a, b$ and consider first the case $a < c < b$. Let $m = \# \{x \in T \; | \; a \leq x < c\}$. Using a similar argument to above, the chain space $C_\bullet^\downarrow I_{[a,b)}(c)$ is:
\[
0 \xrightarrow{} \mathbb{R}^{{m}\choose{m}}\xrightarrow{\partial_{m}^\downarrow} \cdots \xrightarrow{} \mathbb{R}^{{m}\choose{2}}\xrightarrow{\partial_{2}^\downarrow} \mathbb{R}^{{m}\choose{1}}\xrightarrow{\partial_{1}^\downarrow} \mathbb{R} \xrightarrow{} 0.
\]
First, observe that $\partial_1^\downarrow$ is the all-ones $m$-vector. Indeed, each 1-chain $x < c$ for $x \in [a,c^-]$ induces the identity map $I_{[a,b)}(x) \to I_{[a,b)}(c)$. Thus, $H_0^\downarrow I_{[a,b)}(c) = 0$. Then, again relating chains in $[a,c^-]$ with the standard $(m-1)$-simplex, we have that $H_{i}^\downarrow I_{[a,b)}(c) \cong \tilde{H}_{i-1} (\Delta^{m-1})$ (reduced homology of the standard $(m-1)$-simplex) for all $i \geq 1$. That is, $H_i^\downarrow I_{[a,b)}(c) = 0$ for all $i$. 

When $c < a$, it is easy to see that $H_i^\downarrow I_{[a,b)}(c) \cong 0$ for all $i \geq 0$. For~$c > b$, we get the same chain complex as above but the indices are shifted up. That is, combining this paragraph with the preceding one, $H_i^\downarrow I_{[a,b)}(c) = 0$ whenever $c \neq a, b$. 
\end{proof}
\end{lem}

\begin{rmk}
When applying the lemma above to a persistence interval of the form $[a,\infty)$, the M\"obius homology space  associated to the interval module $I_{[a,\infty)}: T \to \Vect_\mathbb{R}$ at $c \in T$ is given by
\[
H_0^\downarrow I_{[a,\infty)}(c) \cong \begin{cases}
\mathbb{R} & \text{if } c = a,\\
0 & \text{else}.
\end{cases}
\]
Furthermore, $H_i^\downarrow I_{[a,\infty)}(c) = 0$ for all $i \geq 0$. Indeed, for $c \in T$, it is never true that $c = \infty$. Thus, we will never have a non-zero degree-one M\"obius homology space for interval modules of the form $I_{[a,\infty)}$.
\end{rmk}

By combining the lemma above with the structure theorem for one-parameter persistence modules, we can compute M\"obius homology spaces for \textit{any} one-parameter module $M$ at \textit{any} index $c \in T$. 

\begin{lem}\label{lem:distributes}
Let $T$ be a finite totally ordered set and let $M: T \to \Vect_\mathbb{R}$ be a one-parameter persistence module with interval decomposition
\[
M \cong \bigoplus_{[a,b) \in\Dgm(T)} I_{[a,b)}^{\oplus m_\mathcal{B}([a,b))}.
\]
Then, for each $c \in T$ and each $i \geq 0$, 
\[
H_i^\downarrow M(c) \cong \bigoplus_{[a,b) \in \Dgm(T)} \left(H_i^\downarrow I_{[a,b)}(c)\right)^{\oplus m_\mathcal{B}([a,b))}.
\]
\begin{proof}
Fix an isomorphism $\varphi: M \Rightarrow \bigoplus_{[a,b) \in\Dgm(T)} I_{[a,b)}^{\oplus m_\mathcal{B}([a,b))}.$ M\"obius homology respects isomorphisms and distributes over direct sums. Thus, we conclude. 
\end{proof}
\end{lem}

\begin{thm}\label{thm:mobius_int_decomp} Let $T$ be a finite totally ordered set and let $M: T \to \Vect_\mathbb{R}$ be a one-parameter persistence module with interval decomposition
\[
M \cong \bigoplus_{[a,b) \in\Dgm(T)} I_{[a,b)}^{\oplus m_\mathcal{B}([a,b))}.
\]
Then, for each $c \in T$, 
\[
\dim H_0^\downarrow M(c)  = \sum_{\substack{[a,b) \in \Dgm(T) \\ a = c}} m_\mathcal{B}([a,b)), \quad \dim H_1^\downarrow M(c)  = \sum_{\substack{[a,b) \in \Dgm(T) \\ b = c}} m_\mathcal{B}([a,b)),
\]
and $H_i^\downarrow M(c) = 0$ for all $i > 1$.
\begin{proof}
The result follows from Lemmas~\ref{lem:mobius_interval_module} and~\ref{lem:distributes}. 
\end{proof}
\end{thm}

In particular, degree-zero M\"obius homology records the birth indices of bars $[a,b)$ in the barcode $m_\mathcal{B}$ whereas degree-one M\"obius homology records the \textit{finite} death indices ($b \neq \infty$). 

\begin{ex}\label{ex:simple-mod}
Let $M: T \to \Vect_\mathbb{R}$ be the module over the totally ordered set $1 < 2 < 3 <4 $ from Example~\ref{ex:simple-fil-barcode}:
$$\mathbb{R} \xrightarrow{\begin{pmatrix}1 \\0\end{pmatrix}} \mathbb{R}^2 \xrightarrow{\begin{pmatrix}0 & 1\end{pmatrix}} \mathbb{R} \xrightarrow{\begin{pmatrix}0 \end{pmatrix}} 0.$$ 
We saw in Example~\ref{ex:simple-fil-barcode} that the interval decomposition of $M$ is given by $M \cong I_{[1,3)} \oplus I_{[2,4)}$. The M\"obius chain complex $C_\bullet^\downarrow M(2)$ is given by 
\[
 \cdots\xrightarrow{\partial_3^\downarrow}  0 \xrightarrow{\partial_2^\downarrow} \mathbb{R}\xrightarrow{\partial_1^\downarrow = \begin{pmatrix} 1 \\ 0
 \end{pmatrix}} \mathbb{R}^2 \xrightarrow{\partial_0^\downarrow} 0. 
\]
A quick calculation confirms that $H_0^\downarrow M(2) \cong \mathbb{R}$ and $H_i^\downarrow M(2) = 0$ for all $i > 0$. 

On the other hand, the M\"obius chain complex $C_\bullet^\downarrow M(3)$ is given by 
\[
 \cdots 0 \xrightarrow{\partial_3^\downarrow}  \mathbb{R} \xrightarrow{\partial_2^\downarrow = \begin{pmatrix}
     -1 \\ 1 \\ 0
 \end{pmatrix}} \mathbb{R}^3\xrightarrow{\partial_1^\downarrow = \begin{pmatrix} 0 & 0 &  1
 \end{pmatrix}} \mathbb{R}\xrightarrow{\partial_0^\downarrow} 0. 
\]
A quick calculation confirms that $H_1^\downarrow M(3) \cong \mathbb{R}$ and $H_i^\downarrow M(3) = 0$ for all $i \neq 1$.
\end{ex}

Thus, as expected, M\"obius homology captures the well-known structure of one-parameter persistence modules $M: T \to \Vect_\mathbb{R}$. By placing an inner product on each space $M(c)$, we may use the M\"obius homology Laplacian to study analytic properties of $M$. Again, we expect that the analytic and spectral properties of $M$ reflect its interval decomposition.

\subsection{M\"obius Homology Laplacian for One-Parameter Modules}
\label{sec:interval-decomp}
 
We rely heavily on the interval decomposition of $M$ to define an inner product structure on $M(c)$. For elements $c \in T$ such that $a \leq c < b$, endow the summand $I_{[a,b)}(c) \cong \mathbb{R}$ with an orthonormal basis (that is, 1 with the standard inner product). For all other values of $c \in T$, $I_{[a,b)}(c) = 0$. Then, $M(c) \cong \bigoplus_{[a,b) \in \Dgm(T)} I_{[a,b)}(c)^{m_\mathcal{B}([a,b))}$, by construction, inherits an orthonormal basis where each non-zero summand $I_{[a,b)}(c)$ determines an orthonormal basis vector in $M(c)$.

With the choice of inner product above, each one-parameter module is now valued in the category of Hilbert spaces: $M: T \to \Hilb_\mathbb{R}$.

Recall that $C_i^\downarrow M(c)$ is given by $\bigoplus_{\text{i-chains } \sigma} M(\min \sigma)$. Thus, we obtain an inner product on $C_i^\downarrow M(c)$ by declaring that vectors in distinct summands are orthogonal and extending linearly. With an inner product structure on each chain space $C_i^\downarrow M(c)$, we can compute its associated M\"obius homology Laplacian $L_i^\downarrow M(c)$. In this section, we will focus specifically on the degree-one M\"obius homology Laplacian $L_1^\downarrow M(c)$. Indeed, recall from Theorem~\ref{thm:mobius_int_decomp} that degree-one M\"obius homology $H_1^\downarrow M(c)$ has dimension given by the number of bars $[a,c)$ in the barcode of $M$ (counted with multiplicity).

The spectrum of the M\"obius homology Laplacian $L_1^\downarrow M(c)$ will depend heavily on the \textit{length} of persistence intervals $[a,b)$ in $\Dgm(T)$.

\begin{defn}\label{def:length-of-bar}
Let $P$ be a finite poset, and let $[a,b)$ be a symbolic persistence interval. Define the \textbf{length} of $[a,b)$ as \[\mathrm{len}([a,b)) = \# \{c \in P \; | \; a \leq c < b\}.\] We emphasize that the death index $b$ is \textit{not} counted in the length calculation. In particular, \[
\mathrm{len}([a,a)) = 0 \quad \text{and} \quad \mathrm{len}([a,\infty)) = \#\{c \in P \; | \; a \leq c \}.
\]
\end{defn}

The proof of Lemma~\ref{lem:mobius_interval_module} immediately shines a light on the degree-one M\"obius homology Laplacian of an interval module $I_{[a,b)}$ at the (finite) element $b\in T$, $L_1^\downarrow I_{[a,b)}(b).$

\begin{cor}\label{cor:laplacian_is_graphlaplacian}
Let $I_{[a,b)}: T \to \Hilb_\mathbb{R}$ be an interval module with $b \neq \infty$ and $n = \mathrm{len}([a,b))$. The degree-one M\"obius homology Laplacian $L_1^\downarrow I_{[a,b)}(b): C_1^\downarrow I_{[a,b)}(b) \to C_1^\downarrow I_{[a,b)}(b)$ is the (unnormalized) graph Laplacian on~$K_n$, the complete graph with $n$ vertices.
\begin{proof}
From the proof of Lemma~\ref{lem:mobius_interval_module}, recall that $C_\bullet^{\downarrow}I_{[a,b)}(b)$ is isomorphic to $C_\bullet(\Delta^{n-1})$. Recall that the degree-zero simplicial Hodge Laplacian $\partial_1^{\Delta^{n-1}} (\partial_1^{\Delta^{n-1}})^T$ on $C_0(\Delta^{n-1})$ is precisely the graph Laplacian on the 1-skeleton of $\Delta^{n-1}$, which is the complete graph on $n$ vertices. Furthermore, with our choice of inner product on $I_{[a,b)}(b)$, adjoint operators are given by the transpose. Since $\partial_1^{\Delta^{n-1}} = \partial_1^\downarrow$ (again, see the proof of Lemma~\ref{lem:mobius_interval_module}), we conclude.
\end{proof}
\end{cor}

\begin{cor}\label{cor:spectrum_is_n}
Let $I_{[a,b)}: T \to \Hilb_\mathbb{R}$ be an interval module with $b \neq \infty$ and $n = \mathrm{len}([a,b))$. Then, the eigenvalues of $L_1^\downarrow I_{[a,b)}(b)$ are 0 (with multiplicity 1) and $n$ (with multiplicity $n-1$).
\begin{proof}
This is a direct consequence of Corollary~\ref{cor:laplacian_is_graphlaplacian} and Proposition~\ref{prop:spectrum_of_L0Kn}.
\end{proof}
\end{cor}

Corollary~\ref{cor:spectrum_is_n} agrees with Lemma~\ref{lem:mobius_interval_module}: for $I_{[a,b)}$ with $b \neq \infty$, $L_1^\downarrow I_{[a,b)}(b)$ having eigenvalue 0 with multiplicity one agrees with the isomorphism $\ker L_1^\downarrow I_{[a,b)}(b) \cong H_1^\downarrow I_{[a,b)}(b) \cong \mathbb{R}$. 

Furthermore, we will utilize the corollary above to describe the eigenvalues of $L_1^\downarrow M(c)$ for \textit{any} one-parameter module $M \cong\bigoplus_{[a,b) \in \Dgm(T)}I_{[a,b)}^{m_\mathcal{B}([a,b))}$ and \textit{any} index $c \in T$. In particular, the spectrum of~$L_1^\downarrow M(c)$ will contain information about the \textit{length} and \textit{location} of bars $[a,b)$ in the barcode of $M$ in relation to the chosen element $c \in T$. We will make this statement precise in the following lemmas. First, we will utilize the following result from linear algebra.

\begin{lem}\label{prop:adding_identity}
    Let $A$ be an $n\times n$ matrix with eigenvalues $\lambda_1, \lambda_2, \ldots, \lambda_n$ and corresponding eigenvectors $\vec{x}_1, \vec{x}_2, \ldots, \vec{x}_n$. Then, the eigenvalues of $A + aI_{n}$ are $a + \lambda_1, a + \lambda_2, \ldots, a+ \lambda_n$ with corresponding eigenvectors $\vec{x}_1, \vec{x}_2, \ldots, \vec{x}_n$.
    \begin{proof} $(A + aI_{n})\vec{x}_i = A\vec{x}_i + aI_{n}\vec{x}_i = \lambda_i\vec{x}_i + a\vec{x}_i = (\lambda_i + a)\vec{x}_i.$
    \end{proof}
\end{lem}

We will now use Corollary~\ref{cor:spectrum_is_n} to show that the spectrum of $L_1^\downarrow M(c)$ reflects the bars $[a,b)$ in the barcode of $M$ ``below" $c$: that is, where $b \leq c$.

\begin{lem}\label{lem:spectrum_of_L1b_aboveb}
Let $M \cong \bigoplus_{[a,b) \in \Dgm (T)}I_{[a,b)}^{m_\mathcal{B}([a,b))}$ be a module over the totally ordered set $T$ and let~$c \in T$. Then, for each $[a,b) \in \Dgm(T)$ so that $m_\mathcal{B}([a,b)) \neq 0$ and $b \leq c$, $L_1^\downarrow M(c)$ has eigenvalues
\begin{align*}
     \mathrm{len}([b,c)) & \text{ with multiplicity }m_\mathcal{B}([a,b)), \\
     \mathrm{len}([a,c)) & \text{ with multiplicity }m_\mathcal{B}([a,b)) \cdot \left(\mathrm{len}([a,b)) - 1 \right).
\end{align*}
\begin{proof}
First, since $b \leq c$ and $c \in T$, we necessarily have that $b \neq \infty$.

Since the inner product of $M(c)$ is such that distinct interval modules are orthogonal, we can write $L_1^\downarrow M(c)$ as direct sum of $L_1^\downarrow M(c)\big|_{I_{[a,b)}}$ for all persistence intervals $[a,b)$ counted with multiplicity $m_\mathcal{B}([a,b))$. Thus, it suffices to look at one bar in the barcode at a time: for this lemma, we will focus on $[a,b)$ so that $b \leq c$. 

We claim that $L_1^\downarrow M(c)\big|_{I_{[a,b)}}$ can be expressed as $L_0^{K_{\mathrm{len}([a,b))}} + \mathrm{len}([b,c))I_{\mathrm{len}([a,b))}$. Indeed, view $L_1^\downarrow M(c)\big|_{I_{[a,b)}}$ as the (unnormalized) graph Laplacian on the complete graph on $\mathrm{len}([a,b))$ vertices with an extra edge attached to each vertex for every $x$ such that $b \leq x < c$. These edges connect an ``empty" vertex corresponding to $x$ to each element $d$ such $a \leq d < b$. Formally, the ``empty" vertices come from our module structure: $I_{[a,b)}(x) = 0$. This amounts to adding $\mathrm{len}([b,c))$ to the degree of each vertex: that is, $$L_1^\downarrow M(c)\big|_{I_{[a,b)}} =  L_0^{K_{\mathrm{len}([a,b))}} + \mathrm{len}([b,c))I_{\mathrm{len}([a,b))}.$$ 

By combining Lemma~\ref{prop:adding_identity} with Corollary~\ref{cor:spectrum_is_n}, the eigenvalues of $L_1^\downarrow M(c)\big|_{I_{[a,b)}}$ are given by $0 + \mathrm{len}([b,c))$ with multiplicity 1 and $\mathrm{len}([a,b)) + \mathrm{len}([b,c)) = \mathrm{len}([a,c))$ with multiplicity~$\mathrm{len}([a,b))-1$. Since the eigenvalues of a direct sum of matrices is the disjoint union of the eigenvalues of each summand, we conclude. 
\end{proof}
\end{lem}

Observe that the above result agrees with Corollary~\ref{cor:spectrum_is_n}. Indeed, consider a single interval module $I_{[a,b)}: T \to \Hilb_\mathbb{R}$ such that $b \neq \infty$. In this case, $m_\mathcal{B}([a,b)) = 1$. Then, $L_1^\downarrow I_{[a,b)}(b)$ has an eigenvalue of $\mathrm{len}([b,b)) = 0$ with multiplicity 1 and another eigenvalue of $\mathrm{len}([a,b))$ with multiplicity~$(\mathrm{len}([a,b)) - 1)$. 

Now, we will see that the spectrum of $L_1^\downarrow M(c)$ also contains information on the bars $[a,b)$ in the barcode of $M$ that ``intersect" $c$: that is, $a < c < b$.

\begin{lem}\label{lem:spectrum_of_L1b_belowb}
Let $M \cong \bigoplus_{[a,b) \in \Dgm (T)}I_{[a,b)}^{m_\mathcal{B}([a,b))}$ be a module over the totally ordered set $T$ and let~$c \in T$. Then, for each $[a,b) \in \Dgm(T)$ so that $m_\mathcal{B}([a,b)) \neq 0$ and $a < c < b$, $L_1^\downarrow M(c)$ has an eigenvalue
\begin{align*}
\mathrm{len}([a,c)) \text{ with multiplicity } m_\mathcal{B}([a,b))\cdot \mathrm{len}([a,c)).
\end{align*}
\begin{proof}
Again, we will consider one summand $L_1^\downarrow M(c)\big|_{I_{[a,b)}}$ of $L_1^\downarrow M(c)$ at a time. Now, we will focus on bars $[a,b)$ in the barcode so that $a < c< b$. 

The Laplacian $L_1^\downarrow M(c)\Big|_{I_{[a,b)}}$ is given by $L_0^{K_{\mathrm{len}([a,c))}} + (\partial_1^\downarrow)^* \partial_1^\downarrow$. Since $c < b$, the boundary map~$\partial_1^\downarrow$ is the $\mathrm{len}([a,c))$-vector of all 1's and, thus, $(\partial_1^\downarrow)^* \partial_1^\downarrow$ is the $\mathrm{len}([a,c))\times \mathrm{len}([a,c))$ matrix of all 1's, denoted $J$. Thus, $L_1^\downarrow M(c)\Big|_{I_{[a,b)}}= L_0^{K_{\mathrm{len}([a,c))}} + J = \mathrm{len}([a,c))I_{\mathrm{len}([a,c))}$. Since the eigenvalues of a diagonal matrix are precisely the diagonal entries, we conclude. 
\end{proof}
\end{lem}

Intuitively, the lemma above states that the spectrum of $L_1^\downarrow M(c)$ contains the length of the ``truncation" of the bar $[a,b)$ at $c$, with multiplicity also given by this length.

\begin{rmk}
In particular, Lemma~\ref{lem:spectrum_of_L1b_belowb} applies to persistence intervals of the form $[a,\infty)$. That is, let $M \cong \bigoplus_{[a,b) \in \Dgm (T)}I_{[a,b)}^{m_\mathcal{B}([a,b))}$ be a module over the totally ordered set $T$ and fix $c \in T$. Then, for each $[a,\infty) \in \Dgm(T)$ so that $m_\mathcal{B}([a,\infty) )\neq 0$ and $a < c$, then $L_1^\downarrow M(c)$ has eigenvalue $\mathrm{len}([a,c))$ with multiplicity $m_\mathcal{B}([a,\infty))\cdot \mathrm{len}([a,c))$
\end{rmk}

Combine the previous two lemmas to get a complete characterization of the eigenvalues of the degree-one Laplacian $L_1^\downarrow M(c)$ for any $c \in T$.

\begin{thm}\label{thm:spectrum_of_L1b}
Let $M \cong \bigoplus_{[a,b) \in \Dgm (T)}I_{[a,b)}^{m_\mathcal{B}([a,b))}$ be a module over the totally ordered set $T$ and let~$c \in T$. Then, the eigenvalues of $L_1^\downarrow M(c)$ are given by
\begin{align*}
\text{For } [a,b) \text{ such that }b \leq c: \;&  \mathrm{len}([b,c)) \text{ with multiplicity } m_\mathcal{B}([a,b)), \\
& \mathrm{len}([a,c)) \text{ with multiplicity } m_\mathcal{B}([a,b)) \cdot (\mathrm{len}([a,b)) - 1), \\ 
\text{For } [a,b) \text{ such that }a < c < b: \;&  \mathrm{len}([a,c)) \text{ with multiplicity } m_\mathcal{B}([a,b)). 
\end{align*}
\begin{proof}
    Since $L_1^\downarrow M(c)$ is diagonalizable (as it is self-adjoint), the number of eigenvalues (counted with multiplicity) must equal the dimension of $C_1^\downarrow M(c)$. By definition,~$C_1^\downarrow M(c)$ is given by the direct sum
    \[
    \bigoplus_{d < c} M(d) \cong \bigoplus_{d < c} \left(\bigoplus_{[a,b) \in \Dgm (T)}I_{[a,b)}(d)^{m_\mathcal{B}([a,b))}\right),
    \]
    and, thus 
    \begin{align*}
    \dim \bigoplus_{d < c} M(d) &= \dim \bigoplus_{d < c} \left(\bigoplus_{[a,b) \in \Dgm (T)}I_{[a,b)}(d)^{m_\mathcal{B}([a,b))}\right)\\
    &= \sum_{d < c} \dim \bigoplus_{[a,b) \in \Dgm (T)}I_{[a,b)}(d)^{m_\mathcal{B}([a,b))}.
    \end{align*}
    Observe that the chain space $C_1^\downarrow M(c)$ only includes summands of interval modules $I_{[a,b)}$ such that $a < c$. The corresponding bar $[a,b)$ either ``intersects" $c$ (meaning $a < c < b$) or does not (meaning $b \leq c$).
    
    If $a < c < b$, then $I_{[a,b)}$ appears $m_\mathcal{B}([a,b)) \cdot \mathrm{len}([a,c))$ many times in the direct sum $C_1^\downarrow M(c)$. In this case, Lemma~\ref{lem:spectrum_of_L1b_belowb} accounts for all possible eigenvalues. On the other hand, if a bar~$[a,b)$ is such that $b\leq c$, then there are $m_\mathcal{B}([a,b)) \cdot \mathrm{len}([a,b))$ many copies of $I_{[a,b)}$ in $C_1^\downarrow M(c)$. In this case, Lemma~\ref{lem:spectrum_of_L1b_aboveb} accounts for all possible eigenvalues.

    In conclusion, Lemmas~\ref{lem:spectrum_of_L1b_aboveb} and~\ref{lem:spectrum_of_L1b_belowb} provide the same number of eigenvalues (with multiplicity) as the dimension of $C_1^\downarrow M(c)$ and, thus, provide the entire spectrum.
\end{proof}
\end{thm}

Theorem~\ref{thm:spectrum_of_L1b} is best illustrated with an example.

\begin{figure}
% !TEX root = main.tex

\centering

\tikzset{every picture/.style={line width=0.75pt}} %set default line width to 0.75pt        

\begin{tikzpicture}[x=0.75pt,y=0.75pt,yscale=-1,xscale=1]
%uncomment if require: \path (0,112); %set diagram left start at 0, and has height of 112

%Straight Lines [id:da8093875541178013] 
\draw    (155,73.6) -- (517.2,72.61) ;
\draw [shift={(519.2,72.6)}, rotate = 179.84] [color={rgb, 255:red, 0; green, 0; blue, 0 }  ][line width=0.75]    (10.93,-3.29) .. controls (6.95,-1.4) and (3.31,-0.3) .. (0,0) .. controls (3.31,0.3) and (6.95,1.4) .. (10.93,3.29)   ;
%Straight Lines [id:da518556871065666] 
\draw [line width=4.5]    (187,61.6) -- (331,61.6) ;
%Straight Lines [id:da37474753004569095] 
\draw [line width=4.5]    (235,47.6) -- (475,47.6) ;
%Straight Lines [id:da9534973645619337] 
\draw [line width=4.5]    (332,33.6) -- (380,33.6) ;
%Straight Lines [id:da549296273927959] 
\draw [line width=1.5]    (187,61.6) -- (331,61.6) ;
\draw [shift={(331,61.6)}, rotate = 180] [color={rgb, 255:red, 0; green, 0; blue, 0 }  ][line width=1.5]      (6.71,-6.71) .. controls (3.01,-6.71) and (0,-3.7) .. (0,0) .. controls (0,3.7) and (3.01,6.71) .. (6.71,6.71) ;
\draw [shift={(187,61.6)}, rotate = 0] [color={rgb, 255:red, 0; green, 0; blue, 0 }  ][fill={rgb, 255:red, 0; green, 0; blue, 0 }  ][line width=1.5]      (0, 0) circle [x radius= 4.36, y radius= 4.36]   ;
%Straight Lines [id:da590407762122599] 
\draw [line width=1.5]    (235,47.6) -- (475,47.6) ;
\draw [shift={(475,47.6)}, rotate = 180] [color={rgb, 255:red, 0; green, 0; blue, 0 }  ][line width=1.5]      (6.71,-6.71) .. controls (3.01,-6.71) and (0,-3.7) .. (0,0) .. controls (0,3.7) and (3.01,6.71) .. (6.71,6.71) ;
\draw [shift={(235,47.6)}, rotate = 0] [color={rgb, 255:red, 0; green, 0; blue, 0 }  ][fill={rgb, 255:red, 0; green, 0; blue, 0 }  ][line width=1.5]      (0, 0) circle [x radius= 4.36, y radius= 4.36]   ;
%Straight Lines [id:da8684559995706247] 
\draw [line width=1.5]    (332,33.6) -- (380,33.6) ;
\draw [shift={(380,33.6)}, rotate = 180] [color={rgb, 255:red, 0; green, 0; blue, 0 }  ][line width=1.5]      (6.71,-6.71) .. controls (3.01,-6.71) and (0,-3.7) .. (0,0) .. controls (0,3.7) and (3.01,6.71) .. (6.71,6.71) ;
\draw [shift={(332,33.6)}, rotate = 0] [color={rgb, 255:red, 0; green, 0; blue, 0 }  ][fill={rgb, 255:red, 0; green, 0; blue, 0 }  ][line width=1.5]      (0, 0) circle [x radius= 4.36, y radius= 4.36]   ;

% Text Node
\draw (181,77.4) node [anchor=north west][inner sep=0.75pt]    {$1$};
% Text Node
\draw (230,77.4) node [anchor=north west][inner sep=0.75pt]    {$2$};
% Text Node
\draw (278,77.4) node [anchor=north west][inner sep=0.75pt]    {$3$};
% Text Node
\draw (325,77.4) node [anchor=north west][inner sep=0.75pt]    {$4$};
% Text Node
\draw (373,77.4) node [anchor=north west][inner sep=0.75pt]    {$5$};
% Text Node
\draw (421,77.4) node [anchor=north west][inner sep=0.75pt]    {$6$};
% Text Node
\draw (470,77.4) node [anchor=north west][inner sep=0.75pt]    {$7$};

\end{tikzpicture}
\vspace{-0.2in}
    \caption{The barcode $m_\mathcal{B}$ for module with interval decomposition $I_{[1,4)} \oplus I_{[2,7)} \oplus I_{[4,5)}$.}
    \label{fig:interval_example}
\end{figure}
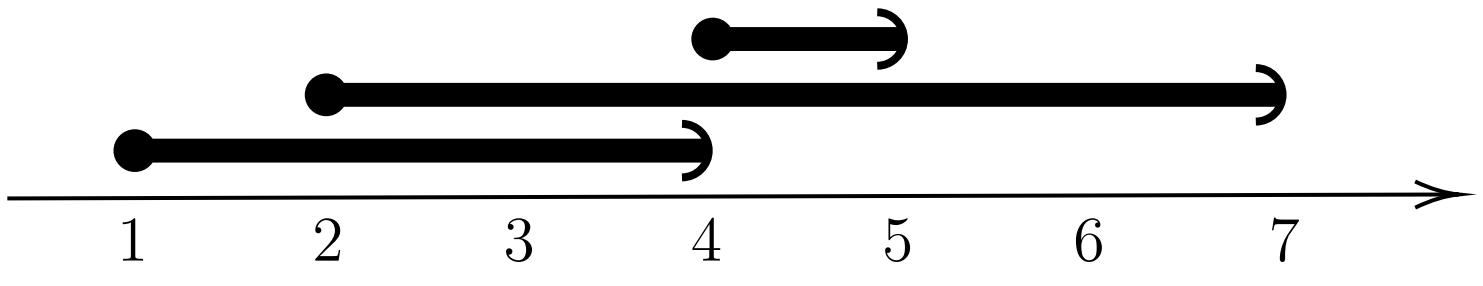

\begin{ex}\label{ex:L_1^4}
    Let $M: T \to \Hilb_\mathbb{R}$ be a module over the totally ordered set $T$ given by all integers in the interval $[1,7]$, and suppose $M \cong I_{[1, 4)} \oplus I_{[2,7)} \oplus  I_{[4,5)}$. See Figure~\ref{fig:interval_example} for a visualization of the barcode $m_\mathcal{B}$.

    We will compute $L_1^\downarrow M(4)$. 
   There are ${3\choose2} = 3$ two-chains in $T$ that end at 4 and ${3\choose1} = 3$ one-chains in $T$ that end at 4. Order these chains lexicographically. The M\"obius chain complex $C_\bullet^\downarrow M(4)$ is given by
   \[
   \cdots \xrightarrow{} \mathbb{R} \oplus \mathbb{R} \oplus \mathbb{R}^2 \xrightarrow{\partial_2^\downarrow} \mathbb{R}\oplus\mathbb{R}^2\oplus\mathbb{R}^2\xrightarrow{\partial_1^\downarrow} \mathbb{R}^2 \xrightarrow{}0
   \]
   where the boundary maps are 
   \[
   \partial_2^\downarrow = \begin{pmatrix} -1 & -1 & 0 & 0 \\
   1 & 0 & -1 & 0 \\
   0 & 0 & 0 & -1 \\
   0 & 1 & 1 & 0\\ 
   0 & 0 & 0 & 1
   \end{pmatrix} \quad \text{and} \quad  \partial_1^\downarrow = \begin{pmatrix} 0 & 0 & 1 & 0 & 1 \\
   0 & 0 & 0 & 0 & 0
   \end{pmatrix}.
   \]
   First, observe that $H_0^\downarrow M(4) = \mathbb{R}^2 / \im \partial_1^\downarrow \cong \mathbb{R}$, which coincides with Theorem~\ref{thm:mobius_int_decomp}. Indeed,~4 is the birth index of the bar $[4, 5)$. The Laplacian $L_1^\downarrow M(4)$ is the matrix given by:
   \[
   L_1^\downarrow M(4) = \partial_2^\downarrow(\partial_2^\downarrow)^T + (\partial_1^\downarrow)^T\partial_1^\downarrow  = \begin{pmatrix}
   2 & -1 & 0 & -1 & 0\\
   -1 & 2 & 0 & -1 & 0\\
   0 & 0 & 2 & 0 & 0 \\
   -1 & -1 & 0 & 2 & 0\\
   0 & 0 & 0 & 0 & 2
   \end{pmatrix}.
   \]
   Its eigenvalues are $\{0, 2, 2, 3, 3\}.$
   
   First, recall Theorem~\ref{thm:mobius_int_decomp}: $\dim\ker L_1^\downarrow M(4) = 1$ reflects the fact that 4 is the death index of the bar~$[1,4)$. Next, recall Lemma~\ref{lem:spectrum_of_L1b_belowb}: for this example, the only bar~$[a,b)$ with $b \leq 4$ is $[1,4)$. The bar~$[1,4)$ has length 3, and its death index is 0 steps away from 4: thus, the spectrum contains~$\{0, 3, 3\}$. Finally, recall Lemma~\ref{lem:spectrum_of_L1b_aboveb}:~$[2,7)$ is the only bar in the barcode that intersects~4. The ``truncated" bar $[2,4)$ has length 2, meaning the spectrum contains $\{2, 2\}$. By Theorem~\ref{thm:spectrum_of_L1b}, these two cases provide the entire spectrum: $\{0, 2, 2,3,3\}$.
 \end{ex}
 
In conclusion, we know precisely what the eigenvalues of $L_1^\downarrow M(c)$ are for \text{any} index $c \in T$. In particular, if $c = \max T$, then we can obtain even more information from the spectral decomposition of $L_1^\downarrow M(c)$. In this case, its \textit{eigenvectors} reflect the \textit{finite} bars in the barcode. We make this statement precise in the following corollary. 

\begin{cor}\label{cor:b_equal_infty}
    Let $M \cong \bigoplus_{[a,b) \in \Dgm (T)} I_{[a,b)}^{m_\mathcal{B}([a,b))}$ be a module over a finite totally ordered set $T$, and endow each $M(x)$ with the basis given by the interval decomposition. Fix a persistence interval $[a,b) \in \Dgm (T)$ with $b \neq \infty$, and let $n = \mathrm{len}([a,b))$.  Then, for each of the $m = m_\mathcal{B}([a,b))$ distinct summands
    \[
    I_{[a,b)}^{(1)}, \cdots, I_{[a,b)}^{(m)}
    \]
    in the interval decomposition of $M$, $L_1^\downarrow M(\max T)$ has an eigenvector $\vec{x}_j$ with eigenvalue $\mathrm{len}([b,\max T))$ so that:
    \begin{enumerate}
    \item $\vec{x}_j$ has exactly $n$ non-zero entries,
    \item each non-zero entry of $\vec{x}_j$ is equal to 1, and
    \item the support of $\vec{x}_j$ consists precisely of the basis elements $
    I_{[a,b)}^{(j)}(x) \in C_1^\downarrow M(\max T)$ for $a \leq x < b.
    $
    \end{enumerate}
    In particular, if $m > 1$, the $m$ eigenvectors $\vec{x}_1, \ldots, \vec{x}_m$ have pairwise disjoint supports.
    \begin{proof}
    With the basis for each $M(x)$ given by the interval decomposition of $M$,  $L_1^\downarrow M(\max T)$ is a block diagonal matrix. In particular, there is a block $L_1^\downarrow M(\max T)\Big|_{I_{[a,b)}^{(j)}}$ for each copy of $I_{[a,b)}$ in the interval decomposition of $M$. Thus, we focus on one such block for the remainder of the proof. 

By Lemma~\ref{lem:spectrum_of_L1b_aboveb}, we have that $L_1^\downarrow M(c)\Big|_{I_{[a,b)}^{(j)}} = L_0^{K_n} + \mathrm{len}([b,\max T))I$, where~$L_0^{K_n}$ denotes the graph Laplacian for the complete graph on $n$ vertices.
By Lemma~\ref{prop:adding_identity}, the eigenvectors of $L_1^\downarrow M(\max T)\Big|_{I_{[a,b)}^{(j)}}$ are the same as the eigenvectors of $L_0^{K_n}$. In this case, the all-ones vector is an eigenvector with eigenvalue $\mathrm{len}([b,\max T))$ (see the proof of Proposition~\ref{prop:spectrum_of_L0Kn}). 

Extending to the entire space $C_1^\downarrow M(\max T)$, we obtain an eigenvector of $L_1^\downarrow M(\max T)$ with exactly $n$ non-zero entries, each equal to 1. These non-zero entries correspond precisely to the basis elements $I_{[a,b)}^{(j)}(x) \in C_1^\downarrow M(\max T)$ for $a \leq x < b$. Since distinct copies $I_{[a,b)}^{(i)}$ determine distinct blocks in $L_1^\downarrow M(\max T)$, the resulting eigenvectors have distinct support.
    \end{proof}
\end{cor}

\begin{ex}\label{ex:L_1^7}
Recall Example~\ref{ex:L_1^4} and Figure~\ref{fig:interval_example}. Since there are no bars in the barcode of the form $[a,\infty)$, the spectral decomposition of $L_1^\downarrow M(7)$ provides a complete picture of the interval decomposition of the module $M \cong I_{[1,4)} \oplus I_{[2,7)} \oplus I_{[4,5)}$. Order the summands in $C_1^\downarrow M(7)$ according to the lexicographic order. That is, write $C_1^\downarrow M(c)$ as follows:
\[
M(1) \oplus M(2) \oplus M(3) \oplus M(4) \oplus M(5) \oplus M(6)
\]
or, in terms of interval modules,
\[
I_{[1,4)}(1) \oplus I_{[1,4)}(2) \oplus I_{[2,7)}(2) \oplus I_{[1,4)}(3) \oplus I_{[2,7)}(3) \oplus I_{[2,7)}(4) \oplus I_{[4,5)}(4) \oplus I_{[2,7)}(5) \oplus I_{[2,7)}(6).
\]
Since the boundary maps $\partial_1^\downarrow$ and $\partial_0^\downarrow$ are much larger than in Example~\ref{ex:L_1^4}, we omit them here. Computing $\partial_1^\downarrow (\partial_1^\downarrow)^* + (\partial_0^\downarrow)^* \partial_0^\downarrow$ yields
\[L_1^\downarrow M(7) = \begin{pmatrix}
5 & -1 & 0 & -1 & 0 & 0 & 0 & 0 & 0\\
-1 & 5 & 0 & -1 & 0 & 0 & 0 & 0 & 0\\
0 & 0 & 4 & 0 & -1 & -1 & 0 & -1 & -1\\
-1 & -1 & 0 & 5 & 0 & 0 & 0 & 0 & 0\\
0 & 0 & -1 & 0 & 4 & -1 & 0 & -1 & -1\\
0 & 0 & -1 & 0 & -1 & 4 & 0 & -1 & -1\\
0 & 0 & 0 & 0 & 0 & 0 & 2 & 0 & 0\\
0 & 0 & -1 & 0 & -1 & -1 & 0 & 4 & -1\\
0 & 0 & -1 & 0 & -1 & -1 & 0 & -1 & 4\\
\end{pmatrix}.
\]
The eigenvalues of $L_1^\downarrow M(7)$ are given by $\{0, 2, 3, 5, 5, 5, 5, 6, 6\}$ with eigenvectors given by:
\[
E = \begin{pmatrix}
    0 & 0 & 1 & 0 & 0 & 0 & 0 & -1 & -1 \\
    0 & 0 & 1 & 0 & 0 & 0 & 0 & 1 & 0\\
    1 & 0 & 0 & -1 & -1 & -1 & -1 & 0 & 0\\
    0 & 0 & 1 & 0 & 0 & 0 & 0 & 0 & 1 \\
    1 & 0 & 0 & 1 & 0 & 0 & 0 & 0 & 0\\
    1 & 0 & 0 & 0 & 1 & 0 & 0 & 0 & 0\\
    0 & 1 & 0 & 0 & 0 & 0 & 0 & 0 & 0\\
    1 & 0 & 0 & 0 & 0 & 1 & 0 & 0 & 0\\
    1 & 0 & 0 & 0 & 0 & 0 & 1 & 0 & 0
\end{pmatrix}.
\]
Here, the eigenvector in column $i$ of $E$ corresponds to the $i$-th smallest eigenvalue. 

Notice the first three eigenvectors: the first column of $E$ with eigenvalue $7-7=0$ has five non-zero entries, identifying the summands $I_{[2,7)}(x)$ for $2 \leq x < 7$. The second column with eigenvalue $7-5=2$ has one non-zero entry, identifying the summand $I_{[4,5)}(4)$. The third column with eigenvalue $7-4=3$ has three non-zero entries, identifying the summands $I_{[1,4)}(x)$ for $1 \leq x < 4$.
\end{ex}

In summary, an eigenvector of $L_1^\downarrow M(\max T)$ consisting of only 1's and 0's identifies a single summand $I_{[a,b)}$  (such that $b \neq \infty$) in the interval decomposition $M \cong \bigoplus_{[a,b) \in \Dgm (T)} I_{[a,b)}^{m_\mathcal{B}([a,b))}$ of a one-parameter module $M: T \to \Hilb_\mathbb{R}$. 

It is, thus, tempting to claim that we can construct the barcode of $M$ using the spectral decomposition of $L_1^\downarrow M(\max T)$. However, there is an important subtlety hiding here. Recall that we considered an inner product on each $M(c)$ so that its interval decomposition has an orthonormal basis. Provided an isomorphism $\varphi: M \Rightarrow \bigoplus_{[a,b) \in \Dgm (T)} I_{[a,b)}^{m_\mathcal{B}([a,b))}$, we can define such an inner product on each $M(c)$ as follows. Denote by $\langle\cdot, \cdot\rangle_\oplus$ the inner product on $M(c)$ so that $\bigoplus_{[a,b) \in \Dgm (T)} I_{[a,b)}^{m_\mathcal{B}([a,b))}$ has an orthonormal basis. Then, for any $x, y \in M(c)$, define 
\[
\langle x, y \rangle_{M(c)} = \langle \phi(x), \phi(y) \rangle_\oplus.
\]
However, having an isomorphism means that we already know the interval decomposition of $M$. To be able to truly construct this interval decomposition from $L_1^\downarrow M(\max T)$, we need to answer the following inverse problem: \textit{which inner products on $M$ guarantee that its interval decomposition has an orthonormal basis?} 

\subsection{Discussion and Directions for Future Work}
\label{subsec:oneparam-open-qs}
In Section~\ref{sec:mob-laplacian}, we introduced the \textit{M\"obius homology Laplacian}. This Laplace operator is self-adjoint and positive semidefinite. Moreover, its kernel is isomorphic to the associated M\"obius homology space and provides cycle representatives for homology classes with minimal norm squared. 

In Section~\ref{sec:interval-decomp}, we provide a complete characterization of the spectrum of $L_1^\downarrow M(c)$ for a one-parameter module $M: T \to \Hilb_\mathbb{R}$ at a fixed index $c \in T$, where $M$ has interval decompositon $\bigoplus_{[a,b) \in \Dgm (T)} I_{[a,b)}^{m_\mathcal{B}([a,b))}$. By Theorem~\ref{thm:spectrum_of_L1b}, the spectrum of $L_1^\downarrow M(c)$ captures the combinatorial structure of the barcode $m_\mathcal{B}: \Dgm(T) \to \mathbb{Z}_{\geq 0}$ for bars $[a,b)$ such that $b \leq c$. Moreover, by Corollary~\ref{cor:b_equal_infty}, certain eigenvalues of $L_1^\downarrow M(\max T)$ identify summands $I_{[a,b)}$ with finite death index, $b \neq \infty$. 

Interesting future work in this direction is as follows. As mentioned at the end of the previous section, these results rely heavily on a choice of inner product on each $M(c)$ that is compatible with its interval decomposition. Determining which inner product structure(s) on $M$ guarantee that $\bigoplus_{[a,b) \in \Dgm (T)} I_{[a,b)}^{m_\mathcal{B}([a,b))}$ admits an orthonormal basis is an interesting and subtle, though potentially difficult, problem. 
A solution to this problem may lead to a new method for computing the interval decomposition of $M$. 

Additionally, the results presented in this chapter rely on the interval decomposition of one-parameter modules. It would be interesting to study whether similar methods could be applied to the spectrum of the M\"obius homology Laplacian for zigzag modules or other multiparameter modules that are known to decompose into indecomposable summands \cite{CadS10}. These questions for future work are a direct continuation of the one-parameter case.

Finally, recall that persistence homology is motivated by the study of filtered simplicial complexes (see Section~\ref{sec:filtrations}). The evolving topology of such filtered complexes give rise to persistence modules: in particular, a \textit{birth-death module}. Furthermore, the \textit{geometry} of a filtered complex leads to a natural choice of an inner product structure on its associated birth-death spaces and hence M\"obius chain spaces. Thus, we are interested in studying \textit{analytical} properties of modules arising from filtrations using the geometry of the underlying simplical complex. In the following section, we will analyze a M\"obius homology Laplacian arising from filtrations (not necessarily one-parameter) of a simplicial complex.

%%%%%%%%%%%%%%%%%%%%%%%%%%%%%%%%%%%%%%%%%%%%%%
%%%%%%%%%%%%%%%%%%%%%%%%%%%%%%%%%%%%%%%%%%%%%%
\section{M\"obius Homology Laplacian for Simplicial Filtrations}
\label{sec:birth-death-laplacian}
%%%%%%%%%%%%%%%%%%%%%%%%%%%%%%%%%%%%%%%%%%%%%%
%%%%%%%%%%%%%%%%%%%%%%%%%%%%%%%%%%%%%%%%%%%%%%
In the previous sections, we defined a Laplace operator on the M\"obius chain complex, called the \textit{M\"obius homology Laplacian}. This Laplacian depends on a module $M: P \to \Hilb_\mathbb{R}$, an index $a \in P$, and dimension $i \geq 0$. When working specifically with one-parameter modules, we saw in Section~\ref{sec:the1dcase} that the spectrum of the M\"obius homology Laplacian captured the structure of the persistence diagram. Now, we will investigate the M\"obius homology Laplacian for the \textit{birth-death module} of a filtered simplicial complex. Unlike Section~\ref{sec:the1dcase}, we are no longer constrained to the one-parameter setting. That is, we consider modules over finite posets $P$. 

The extra data of a filtered simplicial complex will allow for an inner product structure on birth-death spaces that reflects the geometry of the underlying simplicial complex. In particular, this inner product structure will allow us to study analytical properties of the birth-death module by computing its associated M\"obius homology Laplacian. We will see that the kernel of the degree-zero M\"obius homology Laplacian has nice geometric properties. That is, it provides a method for obtaining cycle representatives for M\"obius homology that emphasizes ``newly added" simplices. Furthermore, the spectrum of the M\"obius homology Laplacian for the birth-death module will, again, reflect the combinatorial structure of the persistence diagram. 

\subsection{Birth-Death Module of a Simplicial Filtration}\label{sec:filtrations}

In this section, we will describe the persistent homology pipeline via M\"obius inversion of the birth-death function. We make two key choices. First, we return to the original motivation of persistent homology: simplicial filtrations. Second, we focus on the birth-death function as opposed to the rank function. Indeed, the birth-death function has a natural categorification, whereas the rank function does not (see Appendix~\ref{appendix:rank}). 

\begin{figure}
% !TEX root = main.tex

\centering

\tikzset{every picture/.style={line width=0.75pt}} %set default line width to 0.75pt        

\begin{tikzpicture}[x=0.75pt,y=0.75pt,yscale=-1,xscale=1]
%uncomment if require: \path (0,112); %set diagram left start at 0, and has height of 112

%Shape: Right Triangle [id:dp6008987595516877] 
\draw  [fill={rgb, 255:red, 155; green, 155; blue, 155 }  ,fill opacity=1 ] (360,35.59) -- (430,75.59) -- (360,75.59) -- cycle ;
%Shape: Right Triangle [id:dp9147113879045898] 
\draw   (430,75.59) -- (360,35.59) -- (430,35.59) -- cycle ;
%Shape: Circle [id:dp3932162631971121] 
\draw  [fill={rgb, 255:red, 0; green, 0; blue, 0 }  ,fill opacity=1 ] (426.59,75.59) .. controls (426.59,73.71) and (428.12,72.19) .. (430,72.19) .. controls (431.88,72.19) and (433.41,73.71) .. (433.41,75.59) .. controls (433.41,77.47) and (431.88,79) .. (430,79) .. controls (428.12,79) and (426.59,77.47) .. (426.59,75.59) -- cycle ;
%Shape: Circle [id:dp1126388992833316] 
\draw  [fill={rgb, 255:red, 0; green, 0; blue, 0 }  ,fill opacity=1 ] (426.59,35.59) .. controls (426.59,33.71) and (428.12,32.19) .. (430,32.19) .. controls (431.88,32.19) and (433.41,33.71) .. (433.41,35.59) .. controls (433.41,37.47) and (431.88,39) .. (430,39) .. controls (428.12,39) and (426.59,37.47) .. (426.59,35.59) -- cycle ;
%Shape: Circle [id:dp051351601886223674] 
\draw  [fill={rgb, 255:red, 0; green, 0; blue, 0 }  ,fill opacity=1 ] (356.59,35.59) .. controls (356.59,33.71) and (358.12,32.19) .. (360,32.19) .. controls (361.88,32.19) and (363.41,33.71) .. (363.41,35.59) .. controls (363.41,37.47) and (361.88,39) .. (360,39) .. controls (358.12,39) and (356.59,37.47) .. (356.59,35.59) -- cycle ;
%Shape: Circle [id:dp6497528216109199] 
\draw  [fill={rgb, 255:red, 0; green, 0; blue, 0 }  ,fill opacity=1 ] (356.59,75.59) .. controls (356.59,73.71) and (358.12,72.19) .. (360,72.19) .. controls (361.88,72.19) and (363.41,73.71) .. (363.41,75.59) .. controls (363.41,77.47) and (361.88,79) .. (360,79) .. controls (358.12,79) and (356.59,77.47) .. (356.59,75.59) -- cycle ;
%Shape: Right Triangle [id:dp920777351033554] 
\draw   (231,35.03) -- (301,75.03) -- (231,75.03) -- cycle ;
%Shape: Right Triangle [id:dp8340667003856201] 
\draw   (301,75.03) -- (231,35.03) -- (301,35.03) -- cycle ;
%Shape: Circle [id:dp19205225216726163] 
\draw  [fill={rgb, 255:red, 0; green, 0; blue, 0 }  ,fill opacity=1 ] (297.59,75.03) .. controls (297.59,73.15) and (299.12,71.62) .. (301,71.62) .. controls (302.88,71.62) and (304.41,73.15) .. (304.41,75.03) .. controls (304.41,76.91) and (302.88,78.44) .. (301,78.44) .. controls (299.12,78.44) and (297.59,76.91) .. (297.59,75.03) -- cycle ;
%Shape: Circle [id:dp49129625497124363] 
\draw  [fill={rgb, 255:red, 0; green, 0; blue, 0 }  ,fill opacity=1 ] (297.59,35.03) .. controls (297.59,33.15) and (299.12,31.62) .. (301,31.62) .. controls (302.88,31.62) and (304.41,33.15) .. (304.41,35.03) .. controls (304.41,36.91) and (302.88,38.44) .. (301,38.44) .. controls (299.12,38.44) and (297.59,36.91) .. (297.59,35.03) -- cycle ;
%Shape: Circle [id:dp5475414435711605] 
\draw  [fill={rgb, 255:red, 0; green, 0; blue, 0 }  ,fill opacity=1 ] (227.59,35.03) .. controls (227.59,33.15) and (229.12,31.62) .. (231,31.62) .. controls (232.88,31.62) and (234.41,33.15) .. (234.41,35.03) .. controls (234.41,36.91) and (232.88,38.44) .. (231,38.44) .. controls (229.12,38.44) and (227.59,36.91) .. (227.59,35.03) -- cycle ;
%Shape: Circle [id:dp6481260452525535] 
\draw  [fill={rgb, 255:red, 0; green, 0; blue, 0 }  ,fill opacity=1 ] (227.59,75.03) .. controls (227.59,73.15) and (229.12,71.62) .. (231,71.62) .. controls (232.88,71.62) and (234.41,73.15) .. (234.41,75.03) .. controls (234.41,76.91) and (232.88,78.44) .. (231,78.44) .. controls (229.12,78.44) and (227.59,76.91) .. (227.59,75.03) -- cycle ;
%Shape: Circle [id:dp1254926965928389] 
\draw  [fill={rgb, 255:red, 0; green, 0; blue, 0 }  ,fill opacity=1 ] (166.59,75) .. controls (166.59,73.12) and (168.12,71.59) .. (170,71.59) .. controls (171.88,71.59) and (173.41,73.12) .. (173.41,75) .. controls (173.41,76.88) and (171.88,78.41) .. (170,78.41) .. controls (168.12,78.41) and (166.59,76.88) .. (166.59,75) -- cycle ;
%Shape: Circle [id:dp4722626196120079] 
\draw  [fill={rgb, 255:red, 0; green, 0; blue, 0 }  ,fill opacity=1 ] (96.59,35) .. controls (96.59,33.12) and (98.12,31.59) .. (100,31.59) .. controls (101.88,31.59) and (103.41,33.12) .. (103.41,35) .. controls (103.41,36.88) and (101.88,38.41) .. (100,38.41) .. controls (98.12,38.41) and (96.59,36.88) .. (96.59,35) -- cycle ;
%Shape: Circle [id:dp20536841061239564] 
\draw  [fill={rgb, 255:red, 0; green, 0; blue, 0 }  ,fill opacity=1 ] (96.59,75) .. controls (96.59,73.12) and (98.12,71.59) .. (100,71.59) .. controls (101.88,71.59) and (103.41,73.12) .. (103.41,75) .. controls (103.41,76.88) and (101.88,78.41) .. (100,78.41) .. controls (98.12,78.41) and (96.59,76.88) .. (96.59,75) -- cycle ;
%Straight Lines [id:da019366533377991546] 
\draw    (187,61) -- (218,61.41) ;
\draw [shift={(220,61.44)}, rotate = 180.76] [color={rgb, 255:red, 0; green, 0; blue, 0 }  ][line width=0.75]    (10.93,-3.29) .. controls (6.95,-1.4) and (3.31,-0.3) .. (0,0) .. controls (3.31,0.3) and (6.95,1.4) .. (10.93,3.29)   ;
\draw [shift={(187,61)}, rotate = 0.76] [color={rgb, 255:red, 0; green, 0; blue, 0 }  ][line width=0.75]      (0,-11.18) .. controls (-3.09,-11.18) and (-5.59,-8.68) .. (-5.59,-5.59) .. controls (-5.59,-2.5) and (-3.09,0) .. (0,0) ;
%Straight Lines [id:da7828029083497321] 
\draw    (317,61) -- (348,61.41) ;
\draw [shift={(350,61.44)}, rotate = 180.76] [color={rgb, 255:red, 0; green, 0; blue, 0 }  ][line width=0.75]    (10.93,-3.29) .. controls (6.95,-1.4) and (3.31,-0.3) .. (0,0) .. controls (3.31,0.3) and (6.95,1.4) .. (10.93,3.29)   ;
\draw [shift={(317,61)}, rotate = 0.76] [color={rgb, 255:red, 0; green, 0; blue, 0 }  ][line width=0.75]      (0,-11.18) .. controls (-3.09,-11.18) and (-5.59,-8.68) .. (-5.59,-5.59) .. controls (-5.59,-2.5) and (-3.09,0) .. (0,0) ;
%Straight Lines [id:da566335005881706] 
\draw    (447,61) -- (478,61.41) ;
\draw [shift={(480,61.44)}, rotate = 180.76] [color={rgb, 255:red, 0; green, 0; blue, 0 }  ][line width=0.75]    (10.93,-3.29) .. controls (6.95,-1.4) and (3.31,-0.3) .. (0,0) .. controls (3.31,0.3) and (6.95,1.4) .. (10.93,3.29)   ;
\draw [shift={(447,61)}, rotate = 0.76] [color={rgb, 255:red, 0; green, 0; blue, 0 }  ][line width=0.75]      (0,-11.18) .. controls (-3.09,-11.18) and (-5.59,-8.68) .. (-5.59,-5.59) .. controls (-5.59,-2.5) and (-3.09,0) .. (0,0) ;
%Shape: Right Triangle [id:dp4787065887494397] 
\draw   (100,35) -- (170,75) -- (100,75) -- cycle ;
%Shape: Right Triangle [id:dp8666808913231936] 
\draw  [fill={rgb, 255:red, 155; green, 155; blue, 155 }  ,fill opacity=1 ] (492,35.41) -- (562,75.41) -- (492,75.41) -- cycle ;
%Shape: Right Triangle [id:dp5935749262611805] 
\draw  [fill={rgb, 255:red, 155; green, 155; blue, 155 }  ,fill opacity=1 ] (562,75.41) -- (492,35.41) -- (562,35.41) -- cycle ;
%Shape: Circle [id:dp4046375964466844] 
\draw  [fill={rgb, 255:red, 0; green, 0; blue, 0 }  ,fill opacity=1 ] (558.59,75.41) .. controls (558.59,73.53) and (560.12,72) .. (562,72) .. controls (563.88,72) and (565.41,73.53) .. (565.41,75.41) .. controls (565.41,77.29) and (563.88,78.81) .. (562,78.81) .. controls (560.12,78.81) and (558.59,77.29) .. (558.59,75.41) -- cycle ;
%Shape: Circle [id:dp7945260084318644] 
\draw  [fill={rgb, 255:red, 0; green, 0; blue, 0 }  ,fill opacity=1 ] (558.59,35.41) .. controls (558.59,33.53) and (560.12,32) .. (562,32) .. controls (563.88,32) and (565.41,33.53) .. (565.41,35.41) .. controls (565.41,37.29) and (563.88,38.81) .. (562,38.81) .. controls (560.12,38.81) and (558.59,37.29) .. (558.59,35.41) -- cycle ;
%Shape: Circle [id:dp47084114908835617] 
\draw  [fill={rgb, 255:red, 0; green, 0; blue, 0 }  ,fill opacity=1 ] (488.59,35.41) .. controls (488.59,33.53) and (490.12,32) .. (492,32) .. controls (493.88,32) and (495.41,33.53) .. (495.41,35.41) .. controls (495.41,37.29) and (493.88,38.81) .. (492,38.81) .. controls (490.12,38.81) and (488.59,37.29) .. (488.59,35.41) -- cycle ;
%Shape: Circle [id:dp47443506615638587] 
\draw  [fill={rgb, 255:red, 0; green, 0; blue, 0 }  ,fill opacity=1 ] (488.59,75.41) .. controls (488.59,73.53) and (490.12,72) .. (492,72) .. controls (493.88,72) and (495.41,73.53) .. (495.41,75.41) .. controls (495.41,77.29) and (493.88,78.81) .. (492,78.81) .. controls (490.12,78.81) and (488.59,77.29) .. (488.59,75.41) -- cycle ;

% Text Node
\draw (472,75) node [anchor=north west][inner sep=0.75pt]    {$v_{1}$};
% Text Node
\draw (470,25) node [anchor=north west][inner sep=0.75pt]    {$v_{2}$};
% Text Node
\draw (566,25) node [anchor=north west][inner sep=0.75pt]    {$v_{3}$};
% Text Node
\draw (566,75) node [anchor=north west][inner sep=0.75pt]    {$v_{4}$};

\end{tikzpicture}
\vspace{-0.2in}
\caption{A filtration $F: T \to \Delta K$ for $T$ given by $1 < 2 < 3 < 4$.}
\label{fig:simple-filt}
\end{figure}
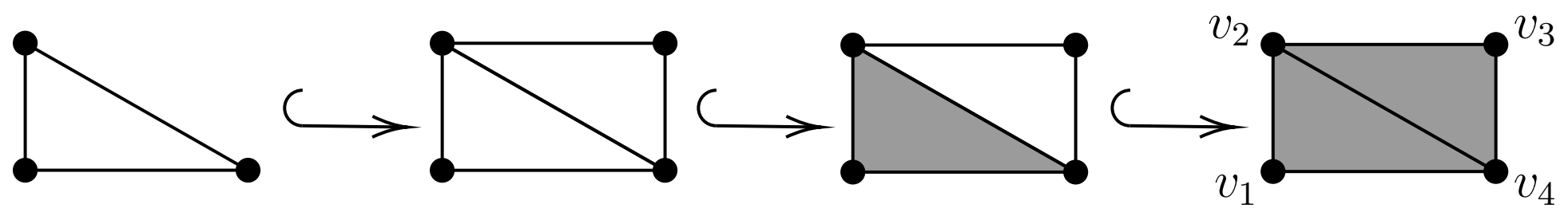

Recall that a \textit{subcomplex} of a simplicial complex $K$ is a subset $L \subseteq K$ such that $L$ is a simplicial complex itself. Given a fixed simplicial complex $K$, we denote by $\Delta K$ the category of subcomplexes of $K$.
\begin{defn}
A \textbf{filtration} of a finite simplicial complex $K$ is a functor from a poset $P$ to the category of subcomplexes of $K$, denoted $F: P \to \Delta K$. 
\end{defn}
In this section, we will consider filtrations over a \textit{finite} poset $P$. See Figure~\ref{fig:simple-filt} for an example of a filtration.

We wish to understand the evolving topology of the simplicial complex $K$ as it moves through a filtration. That is, we wish to answer the following question: when do cycles appear in $K$, and when do they disappear (become boundaries)? 

In order to answer the questions above, we will define a module and a function defined over the collection of symbolic persistence intervals. Thus, we begin by placing an ordering on $\Dgm(P)$.

\begin{defn}\label{def:product-order-general}
Let $P$ be a finite poset with $\overline{P} = P \cup \{\infty\}$. Denote by $(\Dgm(P), \preceq)$ the poset of persistence intervals with the \textbf{product order} $\preceq$. That is, we define 
\[
[a,b) \preceq [c,d) \quad \text{if}  \quad a \leq c \text{ and } b \leq d.\]
\end{defn}

We now define the \textit{birth-death module}.
This module is constructed using the birth-death spaces first introduced in \cite{saecular_persistence, McCPa20-edit-distance}.

\begin{defn}\label{def:birth_death_module_filt}
Let $K$ be a finite simplicial complex, $P$ be a finite poset, and $F: P \to \Delta K$ be a filtration. Fix $p \geq 0$. The \textbf{birth-death module} associated to $F$ is the module $\BD_p F: (\Dgm (P),  \preceq) \to \Vect_\mathbb{R}$ defined on objects by 
\[
\BD_p F([a,b)) = \begin{cases} Z_p F(a) \cap B_p F(b), & \text{if } b \in P, \\
Z_p F(a), & \text{if } b = \infty,
\end{cases}
\]
where the intersection is taken in $C_p(K)$. 

For each relation $[c, d) \preceq [a,b)$ in $(\Dgm (P), \preceq)$, the birth-death module assigns an injective linear map $\BD_p F\left([c, d) \preceq [a,b) \right): \BD_p F([c,d)) \hookrightarrow \BD_p F([a,b))$ induced by the inclusions of cycle and boundary spaces.
\end{defn}

That is, the birth-death module applied to a persistence interval $[a,b)$ consists of cycles at step $a$ of the filtration $F$ that are boundaries at step $b$. 

\begin{defn}[{\cite[Definition 5.6]{McCPa20-edit-distance}}]\label{def:bd-pers-fromfilt}
Let $K$ be a finite simplicial complex and $P$ be a finite poset. Denote by $(\Dgm(P), \preceq)$ the poset of symbolic persistence intervals with the product order. Let $F: P \to \Delta K$ be a filtration, and fix $p \geq 0$. The \textbf{birth-death function} associated to $F$ is the integer-valued map \[
\bd_p F: (\Dgm (P), \preceq) \to \mathbb{Z}_{\geq 0}
\]
given by the dimension function of the birth-death module $\BD_p F: (\Dgm(P), \preceq) \to \Vect_\mathbb{R}$. That is,
\[
\bd_p F([a,b)) = \dim \left( Z_p F(a) \cap B_p F(b) \right) = \dim \BD_p F([a,b)).
\]
\end{defn}

Because the structure maps in the birth-death module are inclusions, the birth-death function is a weakly increasing function.

Intuitively, the birth-death function counts the number of \textit{cycles} in the subcomplex $F(a)$ that are boundaries in the larger subcomplex $F(b)$. In other words, $\bd_p F([a,b))$ counts the number of cycles born at or before $a$ that become boundaries at or before $b$. 

\begin{defn}[{\cite[Definition 8.1]{McCPa20-edit-distance}}]\label{def:pd-from-filt}
Let $K$ be a finite simplicial complex and $P$ be a finite poset. Let $F: P \to \Delta K$ be a filtration, and fix $p \geq 0$. The \textbf{persistence diagram} of $F$ is the M\"obius inversion of the birth-death function
\[
\partial \bd_p F: (\Dgm (P), \preceq) \to \mathbb{Z}.
\]
That is, it is the unique function $\partial \bd_p F$ satisfying 
\[
\bd_p F([a,b)) = \sum_{[c,d) \preceq [a,b)} \partial \bd_p F([c,d)).
\]
\end{defn}

Intuitively, the persistence diagram captures changes in the dimension of the birth-death spaces as we move through the poset $(\Dgm(P), \preceq)$.

\begin{ex}\label{ex:birthdeath-of-filt}
Let $T$ be the totally ordered set $1 < 2 < 3 < 4$, and let $F: T \to \Delta K$ be the filtration shown in Figure~\ref{fig:simple-filt}. The birth-death module $\BD_1 F: (\Dgm(T), \preceq) \to \Vect_\mathbb{R}$ and the resulting persistence diagram $\partial \bd_1 F: (\Dgm(T), \preceq) \to \mathbb{Z}$ are shown in Figure~\ref{fig:birthdeath-of-filt}. 
\end{ex}

\begin{figure}
\include{figures-tex/birthdeath-of-filt}
\vspace{-0.2in}
\caption{From left to right: the birth-death module $\BD_1 F: (\Dgm (T), \preceq) \to \Vect_\mathbb{R}$ associated to the filtration $F$ in Figure~\ref{fig:simple-filt} and the resulting persistence diagram $\partial \bd_1 F: (\Dgm(T), \preceq) \to \mathbb{Z}$. }
\label{fig:birthdeath-of-filt}
\end{figure}

\begin{rmk}
In the one parameter setting, the two approaches to persistent homology described in this paper are known to coincide off the diagonal. Indeed, consider a one-parameter filtration $F: T \to \Delta K$ for $T$ a finite totally ordered set. Define the \textit{persistent homology module} as the composition $H_p F: P \xrightarrow{F} \Delta K \xrightarrow{H_p(-)} \Vect_\mathbb{R}$, where homology is computed for every subcomplex $F(a)$. Let $H_pF \cong \bigoplus I_{[a,b)}^{\oplus m_\mathcal{B}([a,b))}$ be the interval decomposition of $H_pF$. Then, it is known that $m_\mathcal{B}([a,b)) = \bd_pF [a,b)$ for $a \neq b$: this is obtained by combining \cite[Section 9.1]{McCPa20-edit-distance} with \cite[Theorem 12]{CaZo09}.\footnote{Theorem 12 in \cite{CaZo09} actually states that the rank invariant is a complete invariant in the one-parameter setting. However, in the proof of this result, Carlsson and Zomorodian also prove a bijection between possible interval decompositions and possible rank functions.}
\end{rmk}

\subsection{Hilbert Space Structure and the M\"obius Homology Laplacian}\label{sec:defining-L0ab-on-bd}

In order to utilize the M\"obius homology Laplacian to analyze birth-death modules, we must first endow birth-death spaces with an inner product structure. Recall that, for every interval $[a, b) \in (\Dgm (P), \preceq)$, the M\"obius chain complex $C_\bullet^\downarrow \BD_p F([a, b))$ is constructed as follows (see Definition~\ref{def:mobius_chain_complex}):
\[
\cdots \xrightarrow{\partial_{i+1}^\downarrow}
\bigoplus_{\substack{i\text{-chains }c \\ \max c = [a,b)}} \BD_p F(\min c)\xrightarrow{\partial_i^\downarrow} \bigoplus_{\substack{(i-1)\text{-chains }c \\ \max c = [a,b)}} \BD_p F(\min c) \xrightarrow{\partial_{i-1}^\downarrow} \cdots
\]
We equip each birth-death space with the following inner product and use it to induce an inner product on each chain space $C_i^\downarrow\BD_p F([a,b))$. Let $K^{[p]}$ denote the collection of $p$-simplices in~$K$. Given two $p$-cycles $z = \sum_{\sigma \in K^{[p]}} a_\sigma \sigma$ and $w = \sum_{\sigma \in K^{[p]}} b_\sigma \sigma$ in $\BD_p F([a,b))$, define their inner product as follows:
\[
\langle z, w \rangle_{\BD_p F([a,b))} = \sum_{\sigma \in K^{[p]}} a_\sigma b_\sigma.
\]
That is, the inner product on $\BD_p F([a,b))$ is inherited from the standard inner product on $C_p(K)$ as discussed in Section~\ref{subsec:combo-lap}. Indeed, this is the inner product structure commonly used in the literature: see~\cite{basu2022harmonic, Memloi2022Laplacian, lieuter_talk, gulen2025grassmannian, gulen2025grassmannian1d, wang20persistentspectralgraphs}.

Thus, we now view the birth-death module as a module valued in the category of Hilbert spaces: $\BD_p F: (\Dgm (P), \preceq) \to \Hilb_\mathbb{R}$. Thus, we obtain an inner product structure on the M\"obius chain spaces $C_i^\downarrow \BD_p F([a,b))$. The Hilbert space structure allows us to define the M\"obius homology Laplacian for birth-death modules as in Definition~\ref{def:mob-laplacian}. That is, the degree-$i$ M\"obius homology Laplacian for $\BD_p F([a,b))$ is the operator
\[
L_i^\downarrow \BD_p F([a,b)) = (\partial_i^\downarrow)^*\partial_i^\downarrow + \partial_{i+1}^\downarrow(\partial_{i+1}^\downarrow)^*
\]
acting on $C_i^\downarrow\BD_p F([a,b)).$

For the remainder of this section, we will focus specifically on $L_0^\downarrow \BD_p F([a,b)) = \partial_{1}^\downarrow(\partial_{1}^\downarrow)^*$. That is, we will focus on the end of the Möbius chain complex for $\BD_p F([a,b))$, which is given in terms of birth-death spaces below:
\[
\cdots \xrightarrow{\partial_2^\downarrow} \bigoplus_{[c, d) \prec [a,b)} \BD_p F([c,d))\xrightarrow{\partial_1^\downarrow}  \BD_p F([a,b)) \xrightarrow{}  0.
\]

By Theorems~\ref{prop:harmonic_cosheaf_reps} and~\ref{prop:projection-section}, we know that $\ker L_0^\downarrow \BD_p F([a,b)) \cong H_0^\downarrow \BD_p F([a,b))$ and that elements in $\ker L_0^\downarrow \BD_p F([a,b))$ provide canonical cycle representatives for $H_0^\downarrow \BD_p F([a,b))$. Thus, we focus on the degree-zero M\"obius homology Laplacian as it enables an analytic study of the birth-death spaces themselves. 

Since $L_0^\downarrow \BD_p F([a,b))$ is self-adjoint, the spectral theorem yields 
\[
C_0^\downarrow\BD_p F([a,b)) = \BD_p F([a,b)) = \im L_0^\downarrow \BD_p F([a,b)) \oplus \ker L_0^\downarrow \BD_p F([a,b)).
\]
In this section, we will study both the kernel and the image of $L_0^\downarrow \BD_p F([a,b))$. Indeed, we will see that the entire spectral decomposition of $L_0^\downarrow \BD_p F([a,b))$ reflects geometric and combinatorial properties of the persistence diagram $\partial \bd_p F$. First, however, we briefly describe the matrix representation of $L_0^\downarrow \BD_p F([a,b))$ as it will be useful for both theory and computations.

\subsection{Matrix Representation}
\label{sec:matrix_l0ab}
Let $P$ be a finite poset with $\overline{P} = P \cup \{\infty\}$. Denote by $(\Dgm(P), \preceq)$ the poset of symbolic persistence intervals $[a,b)$ where $a \in P$ and $b \in \overline{P}$. Let $F: P \to \Delta K$ be a filtration of a finite simplicial complex $K$, and let $\BD_p F: (\Dgm (P), \preceq) \to \Hilb_\mathbb{R}$ be the associated birth-death module. Fix a persistence interval $[a,b)$. In this section, we will describe the matrix representation of $L_0^\downarrow \BD_p F([a,b))$ according to given bases.

Fix a basis for the birth-death space $\BD_p F([a,b))$, and let $M_{a,b}$ denote the inner product matrix according to this basis.

Let $\{[c_i, d_i)\}_{i = 1}^N$ be the persistence intervals strictly less than $[a,b)$ in $(\Dgm (P), \preceq)$ and fix a basis for each $\BD_pF[c_i, d_i)$. Thus, we obtain a basis for $C_1^\downarrow \BD_p F([a,b)) = \bigoplus_{i = 1}^N \BD_p F([c_i, d_i))$. Let $M$ be the inner product matrix according to this basis.

The boundary map $\bigoplus_{i = 1}^N \BD_p F([c_i, d_i)) \xrightarrow{\partial_1^\downarrow} \BD_p F([a,b))$ can be written as the block matrix 
\[
\partial_1^\downarrow = \begin{pmatrix}
    A_1 & A_2 & \cdots & A_N
\end{pmatrix}
\]
where each $A_i$ is the matrix representation of the inclusion map $\BD_p F[c_i, d_i) \to \BD_p F([a,b))$.

Putting it all together, the Laplacian $L_0^\downarrow \BD_p F([a,b))$ can be written as:
\[
    L_0^\downarrow \BD_p F([a,b)) = \partial_1^\downarrow (\partial_1^\downarrow)^*= \partial_1^\downarrow M^{-1} (\partial_1^\downarrow)^T M_{a,b}.
\]
We will use this matrix representation and its block structure to help us understand and compute both the kernel and the image of $L_0^\downarrow \BD_p F([a,b))$.

\subsection{Harmonic Representatives and Newly Added Simplices}
\label{subsec:harmonic-reps}

In this section, we will explore the kernel of $L_0^\downarrow \BD_p F([a,b))$. Again, by Theorem~\ref{prop:projection-section}, we know that $\ker L_0^\downarrow \BD_p F([a,b))$ provides cycle representatives for $H_0^\downarrow \BD_p F([a,b))$. Moreover, the representatives in $\ker L_0^\downarrow \BD_p F([a,b))$ have a nice geometric interpretation: by Theorem~\ref{prop:minimal-energy}, elements in $\ker L_0^\downarrow \BD_p F([a,b))$ have minimal energy (that is, they have minimal squared norm). In the case of the birth-death module for filtrations, we can say more about the geometric properties of elements in $\ker L_0^\downarrow \BD_p F([a,b))$. Inspired by Basu and Cox's work on harmonic persistent homology \cite{basu2022harmonic}, our choice of an inner product structure on each $\BD_p F([a,b))$ results in harmonic representatives that maximize the contribution of ``newly added" simplices. In the definitions and propositions to follow, we make this claim precise. 

\begin{defn}[{\cite[Definition 4.2]{basu2022harmonic}}]\label{def:support-of-chain}
The \textbf{support} of a $p$-chain $z = \sum_{j = 0}^n a_j \sigma_j \in C_p(K)$ is given by the set $\text{supp}(z) = \left\{ \sigma_j \in K^{[i]} \; | \; a_j \neq 0 \right\}$. 
\end{defn}
That is, the support of a chain $z$ consists of all $p$-simplices with non-zero coefficient in $z$.

\begin{defn}\label{def:relative-new-content}
Let $\sum_{[c,d) \prec [a,b)}\BD_p F([c, d))$ denote the internal sum of the birth-death spaces $\BD_p F([c,d))$, and let $A = \BD_p F([a,b)) \backslash \sum_{[c,d) \prec [a,b)}\BD_p F([c, d))$. Then, given a symbolic persistence interval $[a,b) \in \Dgm(P)$, define
\[
\nu[a,b) = \bigcap_{z \in A} \text{supp}(z).
\]
We call $\nu[a,b)$ the \textbf{set of newly added simplices} associated to the persistence interval $[a,b)$. 
Then, the \textbf{relative new content} of a chain $z = \sum_{\sigma \in K^{[p]}} a_\sigma \sigma \in \BD_p F([a,b))$ is defined as
\[
\mathrm{New}(z) = \left(\frac{\sum_{\sigma \in \nu[a,b)} a_\sigma^2}{\sum_{\sigma \in K^{[p]}} a_\sigma^2}\right)^{1/2}.
\]
\end{defn}
The definitions of the set of newly added simplices associated to an interval and the relative new content of a chain are adapted from Basu and Cox's set of essential simplices and relative essential content, respectively. See \cite[Definitions 4.4 and 4.5]{basu2022harmonic}.

Intuitively, the relative new content of a chain $z$ in $\BD_p F([a,b))$ quantifies the proportion of simplices in $z$ that are ``newly added:" that is, they do not appear in previous birth-death spaces. Observe that, if $\BD_p F([a,b)) \backslash \sum_{[c,d) \prec [a,b)}\BD_p F([c, d)) = \emptyset$, then $\nu[a,b) = \emptyset$ as well. In this case, the relative new content of \textit{any} chain $z \in \BD_p F([a,b))$ will be zero. 

Recall Definition~\ref{def:space-of-reps}: given the natural surjection $\varphi: \BD_p F([a,b)) \to H_0^\downarrow \BD_p F([a,b))$, the \textit{space of cycle representatives} for an equivalence class $[x] \in H_0^\downarrow \BD_p F([a,b))$ is given by $\varphi^{-1}([x])$. 

\begin{thm}\label{prop:relative-new-content}
    Let $F: P \to \Delta K$ be a filtration and fix $p \geq 0$. Consider $[a,b) \in \Dgm(P)$ so that~$\BD_p F([a,b))$ has the property $\dim \BD_p F([a,b)) - \dim \sum_{[c,d) \prec [a,b)} \BD_p F([c,d)) = 1.$ Let $z_0 \in \ker L_0^\downarrow \BD_p F([a,b))$, and let $z \in \varphi^{-1}([z_0])$ be any other cycle representative. Then, 
    \[
    \mathrm{New}(z_0) \geq \mathrm{New}(z).
    \]
    \begin{proof}
        The argument is nearly identical to that of \cite[Theorem 4.8]{basu2022harmonic}, as their theorem was the inspiration for this result. We provide the proof for our specific context here.
        
       Recall that $\partial_1^\downarrow: C_1^\downarrow \BD _pF([a,b)) \to C_0^\downarrow \BD_p F([a,b))$ is the direct sum of the inclusion maps $\BD_p F([c, d) \preceq [a,b))$ for $[c,d) \prec [a,b)$. Thus, $\im \partial_1^\downarrow = \sum_{[c,d) \prec [a,b)} \BD_p F([c,d))$. 

        Let $z_0 \in \ker L_0^\downarrow \BD_p F([a,b))$, and let $z  = \sum_{\sigma \in K^{[p]}} a_\sigma \sigma \in \varphi^{-1}([z_0])$ be a cycle representative of the same homology class $[z_0] \in H_0^\downarrow \BD_p F([a,b))$. Define $z_1 = \sum_{\sigma \in \nu[a,b)} a_\sigma \sigma$ to be only the contribution of ``newly added" simplices, and define $z_2 = z - z_1$. 

        Restated in our context, \cite[Lemma 4.10]{basu2022harmonic} states that $\sigma \in \nu[a,b)$ if and only if $\sigma \in (\im \partial_1^\downarrow)^\perp$. Thus, each $\sigma \in \mathrm{supp}(z_1)$ is in $(\im \partial_1^\downarrow )^\perp$, and therefore $z_1 \in (\im \partial_1^\downarrow )^\perp$. Thus
\begin{align*}
z_0 & = \mathrm{proj}_{(\im\partial_1^\downarrow)^\perp}(z) & \text{(by Theorem~\ref{prop:projection-section})} \\
& = \mathrm{proj}_{(\im\partial_1^\downarrow)^\perp}(z_1) + \mathrm{proj}_{(\im\partial_1^\downarrow)^\perp}(z_2) &\\
& = z_1 + \mathrm{proj}_{(\im\partial_1^\downarrow)^\perp}(z_2) & \text{(since } z_1 \in (\im \partial_1^\downarrow )^\perp\text{).}
\end{align*}

        Now, we will compute $\mathrm{proj}_{(\im \partial_1^\downarrow )^\perp}(z_2)$. For $\sigma \in \nu[a,b)$, $\sigma \notin \mathrm{supp}(z_2)$, and hence \cite[Lemma 4.10]{basu2022harmonic} yields $\sigma \notin \mathrm{supp}(\mathrm{proj}_{(\im\partial_1^\downarrow )^\perp}(z_2))$. In particular, 
        \begin{equation}
        \mathrm{supp}(z_1) \cap \mathrm{supp}(\mathrm{proj}_{(\im \partial_1^\downarrow )^\perp}(z_2)) = \emptyset \implies \langle z_1,
\mathrm{proj}_{(\im\partial_1^\downarrow)^\perp}(z_2)
\rangle
= 0.\label{eq:summands-are-orthogonal}\end{equation}
We will now compare relative new content. 
\begin{align*}
\mathrm{New}(z)^2 & = \frac{\| z_1\|^2}{\|z\|^2} \\
& = \frac{\| z_1\|^2}{\|z_1\|^2 + \| z_2\|^2} & \text{(since } \langle z_1, z_2\rangle = 0\text{ by construction)} \\
& \leq \frac{\|z_1\|^2}
{\|z_1\|^2 + \|\mathrm{proj}_{(\im\partial_1^\downarrow)^\perp}(z_2)\|^2} & \text{(orthogonal projection does not increase norm)} \\
& = \frac{\|z_1\|^2}{\|z_0\|^2} & \text{(by Equation~\ref{eq:summands-are-orthogonal}).}
\end{align*}
Since $\sigma \in \nu[a,b) \implies \sigma \notin \mathrm{supp}(\mathrm{proj}_{(\im\partial_1^\downarrow )^\perp}(z_2))$, and since $z_0 = z_1 + \mathrm{proj}_{(\im\partial_1^\downarrow )^\perp}(z_2)$, the only simplices in $\nu[a,b)$ appearing in $z_0$ are those in $z_1$. That is, $\mathrm{New}(z_0) = \frac{\|z_1\|^2}{\|z_0\|^2}$ and we obtain 
\[
\mathrm{New} (z)^2 \leq  \mathrm{New} (z_0)^2, 
\]
as desired. \end{proof}
\end{thm}

\begin{ex}\label{ex:oneparam-example1}

Consider the filtration $F: T \to \Delta K$ in Figure~\ref{fig:simple-filt} and its birth-death module $\BD_1 F$ as shown in Figure~\ref{fig:birthdeath-of-filt}.

Let the birth-death space $\BD_1F([2,4))$ be generated by $z_1 = [v_1 v_2] + [v_2 v_4] - [v_1 v_4]$ and $z_2 = [v_2 v_3] + [v_3 v_4] - [v_2 v_4]$. The M\"obius chain space $C_1^\downarrow \BD_1 F ([2,4))$, then, is given by the direct sum $$\BD_1 F([1, 3)) \oplus \BD_1 F([1, 4)) \oplus \BD_1 F([2, 3))$$ as these are the only persistence intervals strictly less than~$[2,4)$ in the product order so that $\BD_1 F([c,d)) \neq 0$. All three spaces are generated by~$z_1$. Thus, the internal sum $\BD_1 F([1, 3)) + \BD_1 F([1, 4)) + \BD_1 F([2, 3))$ has dimension 1 and $\BD_1 F([2,4))$ satisfies the requirements of Theorem~\ref{prop:relative-new-content}. 

The end of the M\"obius chain complex $C_\bullet^\downarrow \BD_1F([2,4))$ is given by
\[
\cdots \xrightarrow{\partial_2^\downarrow} \mathbb{R} \oplus \mathbb{R} \oplus \mathbb{R} \xrightarrow{\partial_1^\downarrow = \begin{pmatrix}
    1 & 1 & 1 \\ 0 & 0 & 0
\end{pmatrix}} \mathbb{R}^2 \rightarrow 0.
\]
    We now compute the Laplacian $L_0^\downarrow \BD_1 F([2,4))$. The inner product matrix $M_{2,4}$ is given by $\begin{pmatrix}
3 & -1 \\ -1 & 3
\end{pmatrix}$. The inner product matrix on $C_1^\downarrow \BD_1 F ([2,4))$ is given by $M = 3 \cdot I_3$. Thus, the degree-zero M\"obius homology Laplacian for $\BD_1F([2, 4))$ is 
\[
L_0^\downarrow \BD_1 F([2,4)) = \partial_1^\downarrow M^{-1} (\partial_1^\downarrow)^T M_{2,4} = \begin{pmatrix} 3 & -1 \\ 0 & 0
\end{pmatrix}.
\]
It is easy to check that $\ker L_0^\downarrow \BD_1 F([2,4))$ is generated by $\begin{pmatrix}
1/3 \\ 1
\end{pmatrix}$, which represents the cycle 
\begin{align*}
z &= \frac{1}{3} \left([v_1 v_2] + [v_2 v_4] - [v_1 v_4] \right) + \left( [v_2 v_3] + [v_3 v_4] - [v_2 v_4] \right) \\
&= \frac{1}{3}[v_1 v_2] - \frac{1}{3}[v_1 v_4] + [v_2 v_3] -  \frac{2}{3}[v_2 v_4] + [v_3 v_4].
\end{align*}
The energy of $z$ is $\| z\|^2 = 8/3$. Observe that $z_2 = z - \frac{1}{3} z_1$ is a different cycle representative for the point $[2,4)$ in the persistence diagram. The energy of $z_2$ is given by $\| z_2\|^2 = 3$. 

We also compute the relative new content of $z$ and $z_2$. First, observe that the set of newly added simplices for $[2,4)$ is given by $\nu[2,4) = \{ [v_2 v_3], [v_3v_4]\}$. Then,
\begin{align*}
\mathrm{New}(z) &= \left( \frac{1^2 + 1^2}{8/3}\right)^{1/2} = \sqrt{3/4}, \text{ and} \\
 \mathrm{New}(z_2) & = \left( \frac{1^2 + 1^2}{3}\right)^{1/2} = \sqrt{2/3}.
\end{align*}
Indeed, $\mathrm{New}(z) \geq \mathrm{New}(z_2)$. Although the supports of both $z$ and $z_2$ contain $\nu[2,4)$, the relative new simplices make up a larger proportion of the support of $z \in \ker L_0^\downarrow \BD_1 F([2,4))$.
\end{ex}

Observe that the results of Theorem~\ref{prop:relative-new-content} do not require a one-parameter filtration. See Appendix~\ref{appendix:multiparam} for an example of using the M\"obius homology Laplacian to obtain cycle representatives in the multiparameter setting.

\begin{rmk}
Let $F: P \to \Delta K$ be a filtration, and let $[a,b) \in (\Dgm(P), \preceq)$ be such that $\dim\BD_p F([a,b)) = 1$ and $\dim \BD_p F ([c, d)) = 0$ for all $[c, d) \prec [a, b)$. Since $\BD_p F([a,b))$ satisfies all requirements of Theorem~\ref{prop:relative-new-content}, elements in $\ker L_0^\downarrow \BD_p F([a,b))$ will provide cycle representatives that emphasize newly added content. In this case, there are some subtleties worth pointing out. Computing $C_\bullet^\downarrow \BD_p F([a,b))$ yields
\[
\cdots \xrightarrow{} 0 \xrightarrow{\partial_1^\downarrow} \BD_p F([a,b)) \xrightarrow{} 0.
\]
Thus, $L_0^\downarrow \BD_p F([a,b)) = 0$. Choose a generator $x \in \BD_p F([a,b))$. Then,~$x \in \ker L_0^\downarrow \BD_p F([a,b))$. It is vacuously true that $x$ maximizes relative new content. Since there are \textit{no} cycles present in previous birth-death spaces, \textit{all} simplices present in $x$ are new. Thus, $\mathrm{New}(x) = 1$. 
\end{rmk}

\subsection{Combinatorial Spectral Properties}
\label{subsec:spectral-props}

Let $F: P \to \Delta K$ be a filtration, and let $\BD_p F: (\Dgm (P), \preceq) \to \Hilb_\mathbb{R}$ be its associated birth-death module. In the previous section, we investigated the \textit{kernel} of the degree-zero M\"obius homology Laplacian $L_0^\downarrow \BD_p F([a,b))$. In particular, Theorem~\ref{prop:relative-new-content} states that, under certain conditions, the elements in $\ker L_0^\downarrow \BD_p F([a,b))$ have nice geometric properties. Namely, they have minimal energy, which in this case means they emphasize \textit{relative new content}.

The M\"obius homology Laplacian, however, provides additional structure \textit{beyond} its kernel. Much like classical combinatorial Hodge Laplace operators, we expect that the non-zero spectrum of $L_0^\downarrow \BD_p F([a,b))$ will reflect the combinatorial structure of the persistence diagram associated to $F$, the underlying object in question here.

However, observe that the results from Section~\ref{sec:interval-decomp} do not apply here. Indeed, even when~$P$ \textit{is} totally ordered, $\Dgm (P)$ is often not.\footnote{For $T$ a totally ordered set, the interval poset $(\Dgm (T), \preceq)$ is only totally ordered when $T$ has two or fewer elements.} Thus, new methods are needed to investigate the spectrum of~$L_0^\downarrow \BD_p F([a,b))$. 

Recall that, to compute the eigenvalues of $L_0^\downarrow \BD_p F([a,b))$, we solve for the roots of the characteristic polynomial $$\det(\lambda I_{m} - L_0^\downarrow \BD_p F([a,b))) = \sum_{i = 0}^m a_i \lambda^{m - i}, $$ where $m = \bd_p F([a,b))$. It is known that $\mathrm{tr\,}L_0^\downarrow \BD_p F([a,b))$ (the sum of its eigenvalues) is given by $-a_1$ and that $\det L_0^\downarrow \BD_p F([a,b))$ (the product of its eigenvalues) is given by $(-1)^m a_m$. 

Instead of computing the eigenvalues of $L_0^\downarrow \BD_p F([a,b))$ explicitly, we will relate its trace and the determinant to the persistence diagram $\partial \bd_p F([a,b))$. In fact, we can compute the trace exactly. 

\begin{thm}\label{trace-of-L0ab}
Let $F: P \to \Delta K$ be a filtration, and let $\BD_p F: (\Dgm(P), \preceq) \to \Hilb_\mathbb{R}$ be its associated birth-death module for a fixed $p \geq 0$. Fix $[a,b) \in \Dgm(P)$, and let $L_0^\downarrow \BD_p F([a,b))$ denote the associated M\"obius homology Laplacian. Then, 
\[\mathrm{tr\,}L_0^\downarrow \BD_p F([a,b)) = \sum_{[c,d) \prec [a,b)} \bd_pF([c,d)).\]
\begin{proof}
We begin by changing basis so that, for every interval $[c,d) \preceq [a,b)$, $\BD_p F([c,d))$ has an orthonormal basis. In particular, this means the direct sum $\bigoplus_{[c,d) \prec [a,b)}\BD_p F([c,d))$ also has an orthonormal basis. Denote by $U^{-1}$ and $V^{-1}$ the change of basis matrices for $\BD_p F([a,b))$ and~$\bigoplus_{[c,d) \prec [a,b)}\BD_p F([c,d))$, respectively. Changing basis changes the matrix representation of $L_0^\downarrow \BD_p F([a,b))$ as follows. Denote by $\widetilde{L_0}^\downarrow \BD_p F([a,b))$ the new representation given by 
\begin{align*}
\widetilde{L_0}^\downarrow \BD_p F([a,b)) & = U^{-1}\partial_1^\downarrow V V^{-1} (\partial_1^\downarrow)^* U \\
& = U^{-1}\partial_1^\downarrow V (U^{-1}\partial_1^\downarrow V)^*,
\end{align*}
where the second equality follows from the fact that $U$ and $V$ are orthogonal matrices. Furthermore, since we have changed to orthonormal bases, 
\[
\widetilde{L_0}^\downarrow \BD_p F([a,b)) = U^{-1}\partial_1^\downarrow V (U^{-1}\partial_1^\downarrow V)^T.
\]
Let $\widetilde{\partial_1^\downarrow} = U^{-1}\partial_1^\downarrow V$.
Recall that changing basis does not change eigenvalues.
In particular, $\mathrm{tr\,}L_0^\downarrow \BD_p F([a,b)) = \mathrm{tr\,}\widetilde{L_0}^\downarrow \BD_p F([a,b))$. 

Now, $\widetilde{\partial_1^\downarrow}$ has block structure given by $\begin{pmatrix}
    \widetilde{A_1} & \widetilde{A_2} & \cdots & \widetilde{A_N}
\end{pmatrix}$ (see Section~\ref{sec:matrix_l0ab}). Thus, we can write $\widetilde{L_0}^\downarrow \BD_p F([a,b))=  \widetilde{A_1} \widetilde{A_1}^T +  \widetilde{A_2} \widetilde{A_2}^T + \cdots + \widetilde{A_N} \widetilde{A_N}^T$. Each summand $\widetilde{A_i} \widetilde{A_i}^T$ is a projection matrix that projects onto the subspace $\BD_p F([c_i, d_i)) \subseteq \BD_p F([a,b))$. Thus, $\mathrm{tr\,}\widetilde{A_i} \widetilde{A_i}^T = \dim \BD_p F([c_i, d_i))$ and, because the trace distributes over matrix addition, $$\mathrm{tr\,}L_0^\downarrow \BD_p F([a,b)) = \sum_{[c,d) \prec [a,b)} \dim \BD_p F([c,d))$$
as desired.
\end{proof}
\end{thm}

Furthermore, the methods utilized in the proof of Theorem~\ref{trace-of-L0ab} will also provide insight into the determinant of $L_0^\downarrow \BD_p F([a,b))$. 

\begin{thm}\label{prop:lower-bound-L0ab}
Let $F: P \to \Delta K$ be a filtration, and let $\BD_p F: (\Dgm(P), \preceq) \to \Hilb_\mathbb{R}$ be its associated birth-death module for a fixed $p \geq 0$. Fix $[a,b) \in \Dgm(P)$, and let $L_0^\downarrow \BD_p F([a,b))$ denote the associated M\"obius homology Laplacian. 
Let $m = \bd_p F([a,b))$. Then, $$(\det L_0^\downarrow \BD_p F([a,b)))^{1/m} \geq \# \left\{ [c,d) \; | \; [c,d) \prec [a,b) \text{ and } \bd_pF([c,d)) = \bd_pF([a,b))\right\}.$$
\begin{proof}
First, change basis as in the proof of Theorem~\ref{trace-of-L0ab}. Because changing basis does not change the determinant, $$\det L_0^\downarrow \BD_p F([a,b)) = \det\widetilde{L_0}^\downarrow \BD_p F([a,b)) = \det(\widetilde{A_1} \widetilde{A_1}^T +  \widetilde{A_2} \widetilde{A_2}^T + \cdots + \widetilde{A_N} \widetilde{A_N}^T).$$ By the Minkowski determinant inequality (see \cite[Part II Section 4.1.8]{MarcusMinc1964}), 
\begin{equation}
\begin{split}
&(\det(\widetilde{A_1} \widetilde{A_1}^T +  \widetilde{A_2} \widetilde{A_2}^T + \cdots + \widetilde{A_N} \widetilde{A_N}^T))^{1/m} \\
\geq& (\det(\widetilde{A_1} \widetilde{A_1}^T))^{1/m}  +  (\det(\widetilde{A_2} \widetilde{A_2}^T))^{1/m}  + \cdots + (\det(\widetilde{A_N} \widetilde{A_N}^T))^{1/m} .
\end{split}\label{eq:inductive-minkowski}
\end{equation}
Since each $\widetilde{A_i} \widetilde{A_i}^T$ is a projection matrix, $\det(\widetilde{A_i} \widetilde{A_i}^T) = 1$ if $\bd_p F([c_i,d_i)) = \bd_p F([a,b))$ whereas $\det(\widetilde{A_i} \widetilde{A_i}^T) = 0$ if $\bd_p F([c_i,d_i)) < \bd_p F([a,b))$. Thus, we conclude. 
\end{proof}
\end{thm}

Intuitively, the $1/\bd_pF([a,b))$-th power of the determinant of $L_0^\downarrow \BD_p F([a,b))$ is bounded from below by the number of intervals $[c, d) \prec [a,b)$ where $\bd_pF([c,d)) = \bd_p F([a,b))$. 

The lower bound for $(\det L_0^\downarrow \BD_p F([a,b)))^{1/m} $ given by Theorem~\ref{prop:lower-bound-L0ab} is less informative as the number of intervals $[c,d) \prec [a,b)$ grows. Indeed, Equation~\ref{eq:inductive-minkowski} follows inductively:
\begin{align*}
& (\det(\widetilde{A_1} \widetilde{A_1}^T +  \widetilde{A_2} \widetilde{A_2}^T + \cdots + \widetilde{A_N} \widetilde{A_N}^T))^{1/m} \\  \geq &(\det(\widetilde{A_1} \widetilde{A_1}^T))^{1/m}  +  (\det(\widetilde{A_2} \widetilde{A_2}^T + \cdots + \widetilde{A_N} \widetilde{A_N}^T))^{1/m}  \\
\geq &(\det(\widetilde{A_1} \widetilde{A_1}^T))^{1/m}  +  (\det(\widetilde{A_2} \widetilde{A_2}^T))^{1/m}  + (\det(\cdots + \widetilde{A_N} \widetilde{A_N}^T))^{1/m}  \\
\geq &\cdots \\
\geq &(\det(\widetilde{A_1} \widetilde{A_1}^T))^{1/m}  +  (\det(\widetilde{A_2} \widetilde{A_2}^T))^{1/m}  + \cdots + (\det(\widetilde{A_N} \widetilde{A_N}^T))^{1/m} .
\end{align*}
That is, each summand $\widetilde{A_i} \widetilde{A_i}^T$ adds an inequality.

\begin{ex}
Consider the filtration $F: T \to \Delta K$ shown in Figure~\ref{fig:simple-filt} and its birth-death module $\BD_1 F$ as shown in Figure~\ref{fig:birthdeath-of-filt}. In Example~\ref{ex:oneparam-example1}, we computed the M\"obius homology Laplacian $L_0^\downarrow \BD_1 F([2,4)) = \begin{pmatrix}
    3 & -1 \\ 0 & 0 
\end{pmatrix}$. The characteristic polynomial, then, is given by $$\det(\lambda I_2 - L_0^\downarrow \BD_1 F([2,4)))=
\lambda^2 - 3\lambda.$$
First, observe that $\mathrm{tr\,}L_0^\downarrow \BD_1 F([2,4)) = 3 = \bd_1 F[1,3) + \bd_1 F[1,4) + \bd_1 F[2,3).$ Then, because there are zero birth-death spaces below $\BD_1 F([2,4))$ with the same dimension, $\det L_0^\downarrow \BD_p F([a,b)) = 0$.
\end{ex}

In particular, the two theorems presented in this section relate the combinatorial structure of the birth-death function to the spectrum of $L_0^\downarrow \BD_p F([a,b))$ through linear algebraic invariants (namely, the trace and the determinant). Again, the results of Theorems~\ref{trace-of-L0ab} and~\ref{prop:lower-bound-L0ab} do not require a one-parameter filtration. See Appendix~\ref{appendix:multiparam} for multiparameter examples illustrating Theorems~\ref{trace-of-L0ab} and~\ref{prop:lower-bound-L0ab}.

To conclude this section, we provide an example of a one-parameter filtration $F: T \to \Delta K$ with a degree-zero Laplacian $L_0^\downarrow \BD_p F([a,b))$ whose eigenvalues are not rational.

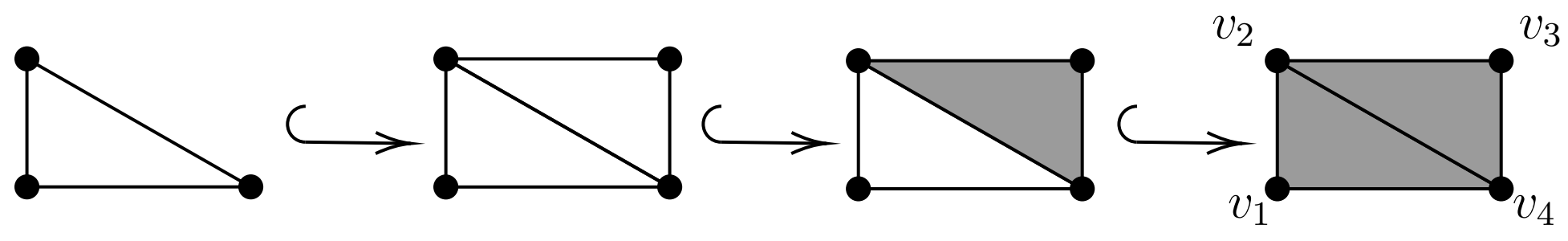
\begin{figure}
    \centering
\tikzset{every picture/.style={line width=0.75pt}} %set default line width to 0.75pt        

\begin{tikzpicture}[x=0.75pt,y=0.75pt,yscale=-1,xscale=1]
%uncomment if require: \path (0,138); %set diagram left start at 0, and has height of 138

%Shape: Right Triangle [id:dp8192669136188677] 
\draw  [fill={rgb, 255:red, 155; green, 155; blue, 155 }  ,fill opacity=1 ] (511,39.59) -- (581,79.59) -- (511,79.59) -- cycle ;
%Shape: Right Triangle [id:dp0643833213536551] 
\draw  [fill={rgb, 255:red, 155; green, 155; blue, 155 }  ,fill opacity=1 ] (581,79.59) -- (511,39.59) -- (581,39.59) -- cycle ;
%Shape: Circle [id:dp8303678135026318] 
\draw  [fill={rgb, 255:red, 0; green, 0; blue, 0 }  ,fill opacity=1 ] (577.59,79.59) .. controls (577.59,77.71) and (579.12,76.19) .. (581,76.19) .. controls (582.88,76.19) and (584.41,77.71) .. (584.41,79.59) .. controls (584.41,81.47) and (582.88,83) .. (581,83) .. controls (579.12,83) and (577.59,81.47) .. (577.59,79.59) -- cycle ;
%Shape: Circle [id:dp019612402619962888] 
\draw  [fill={rgb, 255:red, 0; green, 0; blue, 0 }  ,fill opacity=1 ] (577.59,39.59) .. controls (577.59,37.71) and (579.12,36.19) .. (581,36.19) .. controls (582.88,36.19) and (584.41,37.71) .. (584.41,39.59) .. controls (584.41,41.47) and (582.88,43) .. (581,43) .. controls (579.12,43) and (577.59,41.47) .. (577.59,39.59) -- cycle ;
%Shape: Circle [id:dp42571815728314344] 
\draw  [fill={rgb, 255:red, 0; green, 0; blue, 0 }  ,fill opacity=1 ] (507.59,39.59) .. controls (507.59,37.71) and (509.12,36.19) .. (511,36.19) .. controls (512.88,36.19) and (514.41,37.71) .. (514.41,39.59) .. controls (514.41,41.47) and (512.88,43) .. (511,43) .. controls (509.12,43) and (507.59,41.47) .. (507.59,39.59) -- cycle ;
%Shape: Circle [id:dp584224699270691] 
\draw  [fill={rgb, 255:red, 0; green, 0; blue, 0 }  ,fill opacity=1 ] (507.59,79.59) .. controls (507.59,77.71) and (509.12,76.19) .. (511,76.19) .. controls (512.88,76.19) and (514.41,77.71) .. (514.41,79.59) .. controls (514.41,81.47) and (512.88,83) .. (511,83) .. controls (509.12,83) and (507.59,81.47) .. (507.59,79.59) -- cycle ;
%Shape: Right Triangle [id:dp36245258419819004] 
\draw   (380,39.59) -- (450,79.59) -- (380,79.59) -- cycle ;
%Shape: Right Triangle [id:dp45384535072912136] 
\draw  [fill={rgb, 255:red, 155; green, 155; blue, 155 }  ,fill opacity=1 ] (450,79.59) -- (380,39.59) -- (450,39.59) -- cycle ;
%Shape: Circle [id:dp7199450212806473] 
\draw  [fill={rgb, 255:red, 0; green, 0; blue, 0 }  ,fill opacity=1 ] (446.59,79.59) .. controls (446.59,77.71) and (448.12,76.19) .. (450,76.19) .. controls (451.88,76.19) and (453.41,77.71) .. (453.41,79.59) .. controls (453.41,81.47) and (451.88,83) .. (450,83) .. controls (448.12,83) and (446.59,81.47) .. (446.59,79.59) -- cycle ;
%Shape: Circle [id:dp37191815084774316] 
\draw  [fill={rgb, 255:red, 0; green, 0; blue, 0 }  ,fill opacity=1 ] (446.59,39.59) .. controls (446.59,37.71) and (448.12,36.19) .. (450,36.19) .. controls (451.88,36.19) and (453.41,37.71) .. (453.41,39.59) .. controls (453.41,41.47) and (451.88,43) .. (450,43) .. controls (448.12,43) and (446.59,41.47) .. (446.59,39.59) -- cycle ;
%Shape: Circle [id:dp787525794882526] 
\draw  [fill={rgb, 255:red, 0; green, 0; blue, 0 }  ,fill opacity=1 ] (376.59,39.59) .. controls (376.59,37.71) and (378.12,36.19) .. (380,36.19) .. controls (381.88,36.19) and (383.41,37.71) .. (383.41,39.59) .. controls (383.41,41.47) and (381.88,43) .. (380,43) .. controls (378.12,43) and (376.59,41.47) .. (376.59,39.59) -- cycle ;
%Shape: Circle [id:dp5165101376315151] 
\draw  [fill={rgb, 255:red, 0; green, 0; blue, 0 }  ,fill opacity=1 ] (376.59,79.59) .. controls (376.59,77.71) and (378.12,76.19) .. (380,76.19) .. controls (381.88,76.19) and (383.41,77.71) .. (383.41,79.59) .. controls (383.41,81.47) and (381.88,83) .. (380,83) .. controls (378.12,83) and (376.59,81.47) .. (376.59,79.59) -- cycle ;
%Shape: Right Triangle [id:dp8939804627926395] 
\draw   (251,39.03) -- (321,79.03) -- (251,79.03) -- cycle ;
%Shape: Right Triangle [id:dp7290217700665276] 
\draw   (321,79.03) -- (251,39.03) -- (321,39.03) -- cycle ;
%Shape: Circle [id:dp011212153658897828] 
\draw  [fill={rgb, 255:red, 0; green, 0; blue, 0 }  ,fill opacity=1 ] (317.59,79.03) .. controls (317.59,77.15) and (319.12,75.62) .. (321,75.62) .. controls (322.88,75.62) and (324.41,77.15) .. (324.41,79.03) .. controls (324.41,80.91) and (322.88,82.44) .. (321,82.44) .. controls (319.12,82.44) and (317.59,80.91) .. (317.59,79.03) -- cycle ;
%Shape: Circle [id:dp36709974963290626] 
\draw  [fill={rgb, 255:red, 0; green, 0; blue, 0 }  ,fill opacity=1 ] (317.59,39.03) .. controls (317.59,37.15) and (319.12,35.62) .. (321,35.62) .. controls (322.88,35.62) and (324.41,37.15) .. (324.41,39.03) .. controls (324.41,40.91) and (322.88,42.44) .. (321,42.44) .. controls (319.12,42.44) and (317.59,40.91) .. (317.59,39.03) -- cycle ;
%Shape: Circle [id:dp9098976514284345] 
\draw  [fill={rgb, 255:red, 0; green, 0; blue, 0 }  ,fill opacity=1 ] (247.59,39.03) .. controls (247.59,37.15) and (249.12,35.62) .. (251,35.62) .. controls (252.88,35.62) and (254.41,37.15) .. (254.41,39.03) .. controls (254.41,40.91) and (252.88,42.44) .. (251,42.44) .. controls (249.12,42.44) and (247.59,40.91) .. (247.59,39.03) -- cycle ;
%Shape: Circle [id:dp35643922452244525] 
\draw  [fill={rgb, 255:red, 0; green, 0; blue, 0 }  ,fill opacity=1 ] (247.59,79.03) .. controls (247.59,77.15) and (249.12,75.62) .. (251,75.62) .. controls (252.88,75.62) and (254.41,77.15) .. (254.41,79.03) .. controls (254.41,80.91) and (252.88,82.44) .. (251,82.44) .. controls (249.12,82.44) and (247.59,80.91) .. (247.59,79.03) -- cycle ;
%Shape: Circle [id:dp05275425236595899] 
\draw  [fill={rgb, 255:red, 0; green, 0; blue, 0 }  ,fill opacity=1 ] (186.59,79) .. controls (186.59,77.12) and (188.12,75.59) .. (190,75.59) .. controls (191.88,75.59) and (193.41,77.12) .. (193.41,79) .. controls (193.41,80.88) and (191.88,82.41) .. (190,82.41) .. controls (188.12,82.41) and (186.59,80.88) .. (186.59,79) -- cycle ;
%Shape: Circle [id:dp12083200007355832] 
\draw  [fill={rgb, 255:red, 0; green, 0; blue, 0 }  ,fill opacity=1 ] (116.59,39) .. controls (116.59,37.12) and (118.12,35.59) .. (120,35.59) .. controls (121.88,35.59) and (123.41,37.12) .. (123.41,39) .. controls (123.41,40.88) and (121.88,42.41) .. (120,42.41) .. controls (118.12,42.41) and (116.59,40.88) .. (116.59,39) -- cycle ;
%Shape: Circle [id:dp5791735086443972] 
\draw  [fill={rgb, 255:red, 0; green, 0; blue, 0 }  ,fill opacity=1 ] (116.59,79) .. controls (116.59,77.12) and (118.12,75.59) .. (120,75.59) .. controls (121.88,75.59) and (123.41,77.12) .. (123.41,79) .. controls (123.41,80.88) and (121.88,82.41) .. (120,82.41) .. controls (118.12,82.41) and (116.59,80.88) .. (116.59,79) -- cycle ;
%Straight Lines [id:da29801628313885953] 
\draw    (207,65) -- (238,65.41) ;
\draw [shift={(240,65.44)}, rotate = 180.76] [color={rgb, 255:red, 0; green, 0; blue, 0 }  ][line width=0.75]    (10.93,-3.29) .. controls (6.95,-1.4) and (3.31,-0.3) .. (0,0) .. controls (3.31,0.3) and (6.95,1.4) .. (10.93,3.29)   ;
\draw [shift={(207,65)}, rotate = 0.76] [color={rgb, 255:red, 0; green, 0; blue, 0 }  ][line width=0.75]      (0,-11.18) .. controls (-3.09,-11.18) and (-5.59,-8.68) .. (-5.59,-5.59) .. controls (-5.59,-2.5) and (-3.09,0) .. (0,0) ;
%Straight Lines [id:da7585548284673012] 
\draw    (337,65) -- (368,65.41) ;
\draw [shift={(370,65.44)}, rotate = 180.76] [color={rgb, 255:red, 0; green, 0; blue, 0 }  ][line width=0.75]    (10.93,-3.29) .. controls (6.95,-1.4) and (3.31,-0.3) .. (0,0) .. controls (3.31,0.3) and (6.95,1.4) .. (10.93,3.29)   ;
\draw [shift={(337,65)}, rotate = 0.76] [color={rgb, 255:red, 0; green, 0; blue, 0 }  ][line width=0.75]      (0,-11.18) .. controls (-3.09,-11.18) and (-5.59,-8.68) .. (-5.59,-5.59) .. controls (-5.59,-2.5) and (-3.09,0) .. (0,0) ;
%Straight Lines [id:da8279612500035542] 
\draw    (467,65) -- (498,65.41) ;
\draw [shift={(500,65.44)}, rotate = 180.76] [color={rgb, 255:red, 0; green, 0; blue, 0 }  ][line width=0.75]    (10.93,-3.29) .. controls (6.95,-1.4) and (3.31,-0.3) .. (0,0) .. controls (3.31,0.3) and (6.95,1.4) .. (10.93,3.29)   ;
\draw [shift={(467,65)}, rotate = 0.76] [color={rgb, 255:red, 0; green, 0; blue, 0 }  ][line width=0.75]      (0,-11.18) .. controls (-3.09,-11.18) and (-5.59,-8.68) .. (-5.59,-5.59) .. controls (-5.59,-2.5) and (-3.09,0) .. (0,0) ;
%Shape: Right Triangle [id:dp05416561481279114] 
\draw   (120,39) -- (190,79) -- (120,79) -- cycle ;

% Text Node
\draw (494,80.21) node [anchor=north west][inner sep=0.75pt]    {$v_{1}$};
% Text Node
\draw (489,24.21) node [anchor=north west][inner sep=0.75pt]    {$v_{2}$};
% Text Node
\draw (585,24.21) node [anchor=north west][inner sep=0.75pt]    {$v_{3}$};
% Text Node
\draw (583,80.21) node [anchor=north west][inner sep=0.75pt]    {$v_{4}$};

\end{tikzpicture}
\vspace{-0.2in}
    \caption{A filtration $F: T \to \Delta K$ for $T$ given by $1 < 2 < 3 < 4$.}
    \label{fig:nonint-eigenvalues}
\end{figure}

\begin{figure}
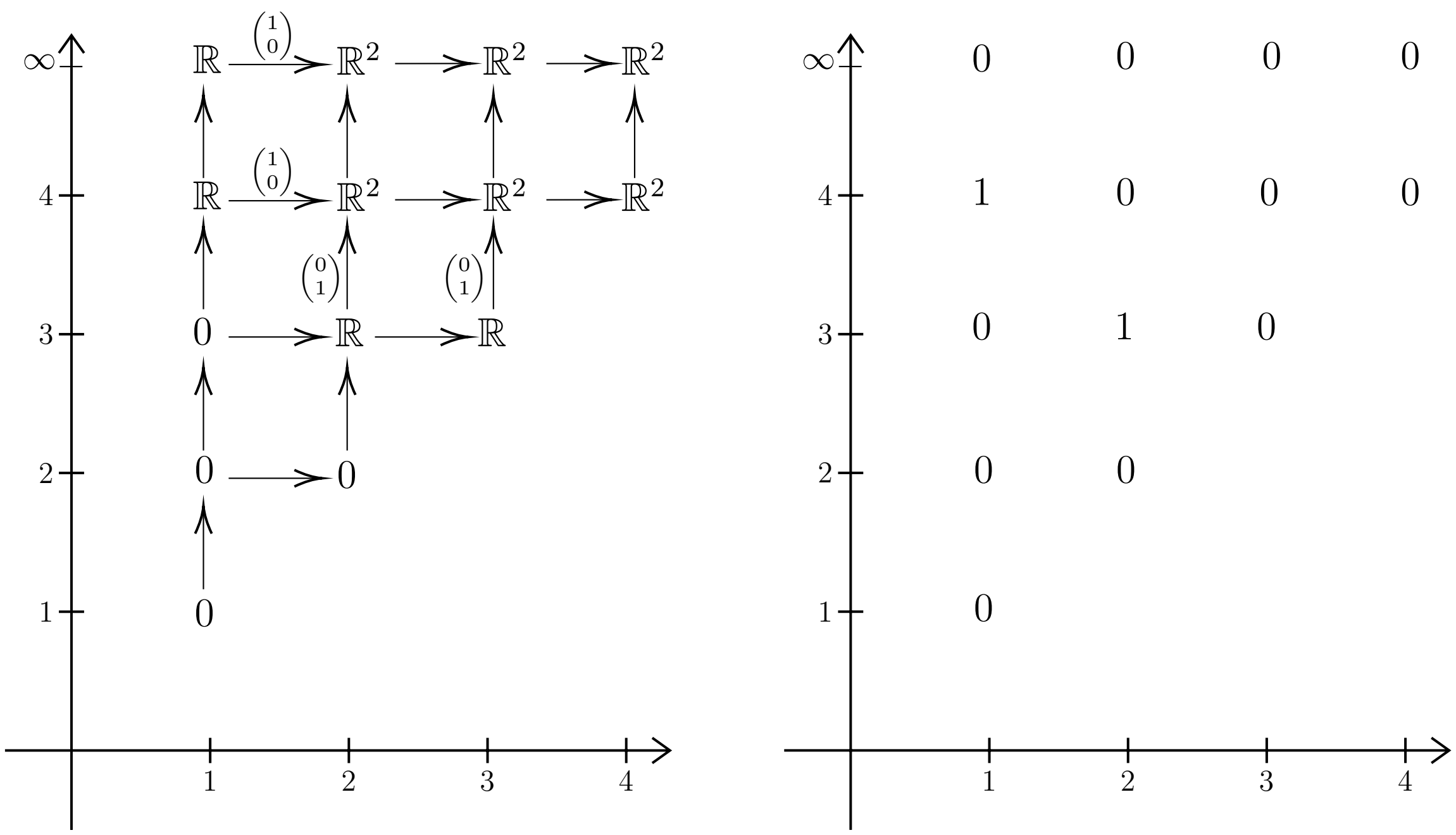

    \centering
    \include{figures-tex/bd-2}
    \vspace{-0.2in}
    \caption{From left to right: the birth-death module $\BD_1 F: (\Dgm(T), \preceq) \to \Vect_\mathbb{R}$ associated to the filtration $F$ in Figure~\ref{fig:nonint-eigenvalues} and the resulting persistence diagram $\partial \bd_1 F:(\Dgm(T), \preceq) \to \mathbb{Z}$.}
    \label{fig:nonint-eigenvalues-bd}
\end{figure}

\begin{ex}
Consider the filtration $F: T \to \Delta K$ given in Figure~\ref{fig:nonint-eigenvalues}. Its birth-death module $\BD_1 F$ and resulting persistence diagram $\partial \bd_1 F$ are shown in Figure~\ref{fig:nonint-eigenvalues-bd}. 

Focus on the interval $[2,4)$. The M\"obius chain space $C_1^\downarrow \BD_1 F([2,4))$ is given by the direct sum $\BD_1 F([1,4)) \oplus \BD_1 F([2,3))$.
The birth-death space $\BD_1 F[2,4)$ is generated by the cycles $z_1 = [v_2 v_3] + [v_3 v_4] - [v_2 v_4]$ and $z_2 = [v_1 v_2] + [v_2 v_4] - [v_1 v_4]$. First, observe that $\BD_1 F([1,3)) \cong 0$. Then, let $\BD_1 F([1,4))$ be generated by $z_2$ and $\BD_1 F([2,3))$ be generated by $z_1$.

The M\"obius homology Laplacian $L_0^\downarrow \BD_1 F([2,4))$ is given by
\[
L_0^\downarrow \BD_1 F([2,4)) = \begin{pmatrix}
    1/2 & 1/2 \\
    1/6 & 3/2
\end{pmatrix}.
\]
The eigenvalues of $L_0^\downarrow \BD_1 F([2,4))$ are $1 \pm \frac{1}{\sqrt{3}}$, and it has characteristic polynomial
\[
\det \left(\lambda I_2 - L_0^\downarrow \BD_1 F([2,4)) \right) = \lambda^2 - 2\lambda + \frac{2}{3}.
\]
In this case, $$\mathrm{tr\,} L_0^\downarrow \BD_1 F([2,4)) = 2 = \bd_1 F([1,4)) + \bd_1 F([2,3))$$ and $$\det \left( L_0^\downarrow \BD_1 F([2,4))\right)^{1/2} = (2/3)^{1/2} \geq 0.$$

\end{ex}

\subsection{Discussion and Future Directions}
\label{subsec:bd-open-qs}

In this section, we considered the degree-zero M\"obius homology Laplacian for the birth-death module of a simplicial filtration. Under certain conditions, we found that $\ker L_0^\downarrow \BD_p F([a,b))$ provided cycle representatives for the M\"obius homology space $H_0^\downarrow \BD_p F([a,b))$. These representatives maximize relative new content in the same way that Basu and Cox's harmonic representatives maximize relative essential content \cite{basu2022harmonic}. Furthermore, we provided a first spectral analysis of $L_0^\downarrow \BD_p F([a,b))$. We computed the trace of $L_0^\downarrow \BD_p F([a,b))$ exactly and obtained a lower bound on the determinant of $L_0^\downarrow \BD_p F([a,b))$. Both quantities were given in terms of the birth-death function $\bd_p F$. 

There are many avenues for future work that remain to be investigated. For one, further investigation into the characteristic polynomial of $L_0^\downarrow \BD_p F([a,b))$ would be interesting. Indeed, when $\dim \BD_p F([a,b)) > 2$, there are more coefficients in $\det(\lambda I_m - L_0^\downarrow \BD_p F([a,b)))$ than only the trace and determinant. It would be interesting to investigate whether the coefficients of this characteristic polynomial admit combinatorial or topological interpretations, as in the case of the graph and combinatorial Laplacians \cite{Bernardi16, Contreras24, Biggs93AlgebraicGraphTheory}. A better understanding of the characteristic polynomial of $L_0^\downarrow \BD_p F([a,b))$ will aid in better understanding of its spectrum as a whole.

Furthermore, in this work, we investigated the \textit{degree-zero} M\"obius homology Laplacian $L_0^\downarrow \BD_p F([a,b))$. The appearance of non-diagonal persistence intervals $[a,b)$ such that $\partial \bd_p F ([a,b)) < 0$, in part, motivated the introduction of M\"obius homology. Indeed, a \textit{negative} value for $\partial \bd_p F ([a,b))$ suggests the presence of odd degree M\"obius homology spaces $H_{2i + 1}^\downarrow \BD_p F([a,b))$ (see Equation~\ref{eq:mob-eulerchar}). Investigating the degree-one M\"obius homology Laplacian, for example, may provide useful insight.

Work by G\"ulen and McCleary provides a method for computing a birth-death module and the associated persistence diagram for a persistence module purely algebraically, utilizing free presentations rather than requiring the data of a simplicial filtration \cite{GulenMcCleary}. Future work would investigate properties of the M\"obius homology Laplacian for their birth-death module and compare to the work presented in this paper.

Finally, the M\"obius homology Laplacian is isomorphic to existing constructions of cycle representatives: see Appendix~\ref{appendix:lit-review} for an explicit relationship between our work and Gülen, Mémoli, and Wan's Grassmannian persistence diagram \cite{gulen2025grassmannian, gulen2025grassmannian1d}. Deeper connections between the two constructions remain to be investigated.

\section{Conclusion}
\label{sec:conclusion}

In this paper, we introduced the M\"obius homology Laplacian as a tool for studying the structure of persistence modules. Our results fell under two main categories: first, we provided explicit eigenvalue computations for the M\"obius homology Laplacian of one-parameter modules. Second, we used the kernel of the M\"obius homology Laplacian for the birth-death module to obtain cycle representatives, then we related linear algebraic invariants (the trace and the determinant) to the structure of the corresponding persistence diagram. However, there is much research to be done to better understand both the kernel and the non-zero spectrum of the M\"obius homology Laplacian in persistence theory and beyond.

%%%%%%%%%%%%%%%%%%%%%%%%%%%%%%%%%%%%%%%%%%%%%%
%%%%%%%%%%%%%%%%%%%%%%%%%%%%%%%%%%%%%%%%%%%%%%
%\section{Statements and Declarations}
%%%%%%%%%%%%%%%%%%%%%%%%%%%%%%%%%%%%%%%%%%%%%%
%%%%%%%%%%%%%%%%%%%%%%%%%%%%%%%%%%%%%%%%%%%%%%
%\paragraph{Conflict of Interest}
%\todo{Add statement on AI use?}
%\newpage
%%%%%%%%%%%%%%%%%%%%%%%%%%%%%%%%%%%%%%%%%%%%%%
\bibliographystyle{plain}
\bibliography{references}
%%%%%%%%%%%%%%%%%%%%%%%%%%%%%%%%%%%%%%%%%%%%%%

\appendix
\newpage

\section{Relationship to Existing Work}
\label{appendix:lit-review}
As discussion in Section~\ref{sec:related-work}, we are not the first to use harmonic methods to study simplicial filtrations or persistence diagrams. 
In our work, particularly in Section~\ref{sec:birth-death-laplacian}, we emphasize the M\"obius inversion approach to persistent homology and obtain higher order structure via the birth-death module and M\"obius homology. 

A similar approach to obtaining higher-order structure for persistence using the theory of M\"obius inversion appears in a pair of papers by G\"ulen, M\'emoli, and Wan in which they introduce \textit{orthogonal M\"obius inversion} and explore various properties \cite{gulen2025grassmannian, gulen2025grassmannian1d}. In general, their approach applies a generalization of M\"obius inversion to Hilbert spaces directly, rather than passing to integer-valued functions. In particular, they apply this inversion to the birth-death module of a simplicial filtration $F: P \to \Delta K$. We summarize their construction in the definition below, combining \cite[Definition 3.4]{gulen2025grassmannian} and \cite[Definition 4.1]{gulen2025grassmannian}.

\begin{defn}[\cite{gulen2025grassmannian}]
Let $F: P \to \Delta K$ be a filtration, let $p \geq 0$, and consider the birth-death module $\BD_p F: (\Dgm(P), \preceq) \to \Hilb_\mathbb{R}$.\footnote{Use the inner product structure on birth-death spaces $\BD_pF([a,b))$ defined in Section~\ref{sec:defining-L0ab-on-bd}.} The \textbf{Grassmannian persistence diagram} of $\BD_pF$ is the function $\mathrm{OI} \left( \BD_p F\right): (\Dgm(P), \preceq) \to \Hilb_\mathbb{R}$ where 
\[
\mathrm{OI} \left( \BD_p F\right) ([a,b)) = \BD_p F([a,b)) \cap \left( \sum_{[c,d) \prec [a,b)} \BD_p F([c,d))\right)^\perp.
\]
Recall that $ \sum_{[c,d) \prec [a,b)} \BD_p F([c,d))$ denotes the internal sum of the birth-death spaces, and here we take the internal sum within the larger birth-death space $\BD_p F([a,b))$.
\end{defn}

We refer interested readers to their paper for more details and important properties, but of particular interest to this work is the following isomorphism.

\begin{thm}[{\cite[Theorem 9]{gulen2025grassmannian}}]
\label{thm:grassmannianPD-mob-iso}
Let $F: P \to \Delta K$ be a filtration with birth-death module $\BD_p F: (\Dgm(P), \preceq) \to \Hilb_\mathbb{R}$. Fix an interval $[a,b) \in \Dgm(P)$. Then, 
\[
\mathrm{OI} \left( \BD_p F\right) ([a,b)) \cong H_0^\downarrow \BD_p F([a,b)).
\]
\end{thm}

Thus, we immediately obtain the following connection between our construction and the Grassmannian persistence diagram.

\begin{cor}
Let $F: P \to \Delta K$ be a filtration with $\BD_p F: (\Dgm(P), \preceq) \to \Hilb_\mathbb{R}$ the birth-death module for $p \geq 0$. Fix an interval $[a,b) \in \Dgm(P)$. Then, 
\[
\ker L_0^\downarrow \BD_p F([a,b)) \cong \mathrm{OI} \left( \BD_p F\right) ([a,b)).
\]
\begin{proof}
    Combine Theorem~\ref{thm:grassmannianPD-mob-iso} with Theorem~\ref{prop:harmonic_cosheaf_reps}.%, which states that $\ker L_0^\downarrow \BD_p F([a,b)) \cong H_0^\downarrow \BD_p F([a,b))$.
\end{proof}
\end{cor}

Interesting future work includes investigating deeper connections between the two constructions.

\section{Proofs from Section~\ref{sec:mob-laplacian}}
\label{appendix:proofs}
In this appendix, we provide the proofs for statements in Section~\ref{sec:mob-laplacian}. Because these proofs are standard in discrete Hodge theory, we provide them here as opposed to the main body of the text. 

\subsection{Proof of Theorem~\ref{prop:projection-section}}
\label{appendix:firstproof}

Recall the statement of the theorem: Let $P$ be a finite poset, let $c \in P$, and let $M:P \to \Hilb_\mathbb{R}$ be a persistence module.
The inclusion $s: \ker L_i^\downarrow M(c) \hookrightarrow \ker \partial_i^\downarrow$ induces a section of the quotient map $\varphi: \ker \partial_i^\downarrow \to H_i^\downarrow M(c)$.
\begin{proof}
Let $[x] = [y] \in H_i^\downarrow M(c)$, meaning $x = y + \partial_{i+1}^\downarrow z$ for some $z \in C_{i+1}^\downarrow M(c)$. Since
\[\mathrm{proj}_{(\mathrm{im}\partial_{i+1}^\downarrow)^\perp}(x)
= \mathrm{proj}_{(\mathrm{im}\partial_{i+1}^\downarrow)^\perp}(y + \partial_{i+1}^\downarrow z)
= \mathrm{proj}_{(\mathrm{im}\partial_{i+1}^\downarrow)^\perp}(y),\]
the map $s$ is well-defined and, thus, yields a canonical representative for $[x] \in H_i^\downarrow M(c)$ contained in $\ker L_i^\downarrow M(c)$.
\end{proof}

\subsection{Proof of Theorem~\ref{prop:minimal-energy}}
\label{appendix:secondproof}

Recall the statement of the theorem: Let $P$ be a finite poset, let $c \in P$, and let $M:P \to \Hilb_\mathbb{R}$ be a module. Consider the degree-$i$ M\"obius homology Laplacian $L_i^\downarrow M(c)$, and let $H_i^\downarrow M(c)$ denote the respective M\"obius homology space.
Let $x \in \ker L_i^\downarrow M(c)$ and let $y = x + \partial_{i+1}^\downarrow z$ be another cycle representative of the same homology class. Then,
\[
\| x\|^2 \leq \|y \|^2.
\]
That is, $x \in \ker L_i^\downarrow M(c)$ is the unique representative with minimal energy.
\begin{proof}
Since $y \in \varphi^{-1}([x])$, $[y] = [x]$ and we can write $y = x + \partial_{i+1}^\downarrow z$ for some $z \in C_{i+1}^\downarrow M(c)$. Using the inner product structure on $C_i^\downarrow M(c)$, 
\begin{align*}
\| y\|^2 = \langle y,y \rangle &= \langle y, x \rangle + \langle y, \partial_{i+1}^\downarrow z \rangle\\
& = \langle x, x \rangle + \langle \partial_{i+1}^\downarrow z, x \rangle + \langle x, \partial_{i+1}^\downarrow z \rangle +  \langle \partial_{i+1}^\downarrow z, \partial_{i+1}^\downarrow z \rangle\\
&= \langle x, x \rangle +  \langle \partial_{i+1}^\downarrow z, \partial_{i+1}^\downarrow z \rangle\\
&\geq \| x\|^2.
\end{align*}
Note that equality holds if and only if $\partial_{i+1}^\downarrow z = 0 \implies y = x$.
\end{proof}

\section{On the Choice of the Birth-Death Approach}
\label{appendix:rank}

In this paper, we used the \textit{birth-death module} when describing the evolving topology of a simplicial complex $K$ through a filtration $F: P \to \Delta K$. The birth-death module categorifies the birth-death function, which is then used to compute the persistence diagram of $F$.

On the other hand, the rank function is another commonly used function in the literature for computing persistent homology \cite{ZoCa05, clause2024generalizedrankinvariantmobius, kim2018generalized}. The rank function, however, does not have a natural categorification, even in the one-parameter setting. Indeed, consider a one-parameter module $M: T \to \Vect_\mathbb{R}$. First, define a new ordering on $\Dgm(T)$.

\begin{defn} Given a finite totally ordered set $T$, the \textbf{reverse containment order} $\supseteq$ on $\Dgm(T)$ is given by setting $[a,b) \supseteq [c,d)$ when $a \leq c$ and $d \leq b$.
\end{defn}

\begin{defn} Given a one-parameter module $M: T \to \Vect_\mathbb{R}$, define the \textbf{rank function} $\rk_M: (\Dgm(T), \supseteq) \to \mathbb{Z}_{\geq 0}$ as follows: $$\rk_M([a,b)) = \dim(\im M(a\leq b^-)),$$ where $b^-$ denotes the greatest element in $T$ that is strictly less than $b$. 
\end{defn}

The natural choice for lifting the rank function to a vector space-valued module is the \textit{image}. However, taking images is \textit{not} functorial. To demonstrate this claim, we provide a counterexample below. 

\begin{ex}
Let $T$ be the totally ordered set $1 < 2 < 3 < 4 < 5$, and let $M: T \to \Vect_\mathbb{R}$ be the module given by
\[
\mathbb{R} \xrightarrow{\begin{pmatrix} 1 \\ 0
\end{pmatrix}} \mathbb{R}^2 \xrightarrow{I_2} \mathbb{R}^2 \xrightarrow{\begin{pmatrix} 1& 0 \\ 0 & 0 \\ 0  & 0
\end{pmatrix}} \mathbb{R}^3 \to 0
\]
where $I_2$ denotes the identity matrix.

The reverse containment order on $\Dgm(T)$ admits the relation $[1,5) \supseteq [2,4)$. We compute the images $\im M(1 \leq 5^-) = \im M(1 \leq 4)$ and $\im M(2 \leq 4^-) = \im M(2 \leq 3)$ below:
\[
\im M(1 \leq 4) = \left\langle\begin{pmatrix}
1 \\ 0 \\ 0
\end{pmatrix}\right\rangle \quad \text{and} \quad \im M(2 \leq 3) = \left\langle\begin{pmatrix}
1 \\ 0
\end{pmatrix}, \begin{pmatrix}
0 \\ 1
\end{pmatrix}\right\rangle.
\]
However, observe the data of the module $M$ does not provide a canonical map from the space $\im M(1 \leq 4)$ to $\im M(2 \leq 3)$. That is, the relation $[1,5) \supseteq [2,4)$ does not induce a map between the corresponding images. 
\end{ex}

Indeed, the example above demonstrates that images of linear maps, the natural choice for a categorification of the rank function, are \textit{not} functorial. Thus, in order to obtain higher-order structure with which to study simplicial filtartions, the birth-death module is prefered.

\section{Multiparameter Examples}
\label{appendix:multiparam}

In Section~\ref{sec:birth-death-laplacian}, we no longer restrict ourselves to one-parameter persistence. Indeed, many of the results in that section apply even for multiparameter filtrations. We now provide two examples of multiparameter filtrations and their associated M\"obius homology Laplacians.

\begin{figure}
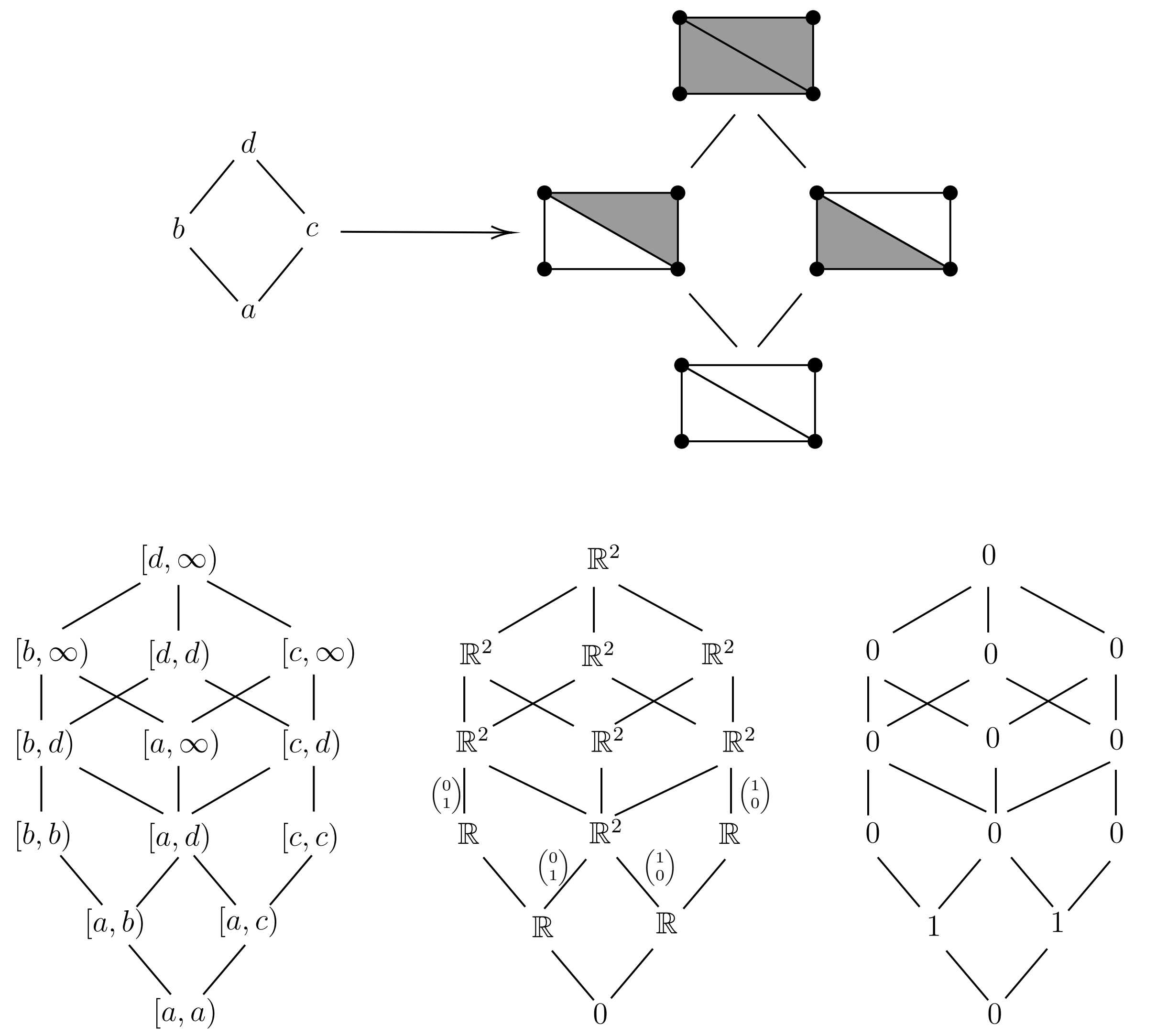

\include{figures-tex/multiparam-counterex}
\vspace{-0.2in}
\caption{Above: a filtration $F: P \to \Delta K$ where $P$ is not totally ordered. Below: $(\Dgm(P), \preceq)$ (left), the birth-death module $\BD_1 F$ (center), and the persistence diagram $\partial \bd_1F$ (right). In $\BD_1F$, all unlabeled maps are either the zero map or the identity.}
\label{fig:multiparam-counterex}
\end{figure}

\begin{ex}\label{ex:multiparam-counterex}
Consider the filtration $F: P \to \Delta K$ in Figure~\ref{fig:multiparam-counterex}. Obtain an order for the birth-death spaces in $C_1^\downarrow \BD_1F[d,d)$ by going from top to bottom, left to right in the Hasse diagram of $(\Dgm(P), \preceq)$. Let $\BD_1F([a,b))$ and $\BD_1 F[b,b)$ be generated by the cycle $z = [v_1 v_2] + [v_2 v_4] - [v_1 v_4]$, and let $\BD_1 F[a,c)$ and $\BD_1 F[c,c)$ be generated by the cycle $w = [v_2 v_3] + [v_3 v_4] - [v_2 v_4]$. The remaining birth-death spaces (for intervals $[a,d), [b,d),$ and $[c,d)$) are generated by both $z$ and $w$. Use these bases to obtain a basis for the direct sum $C_1^\downarrow \BD_1 F[d,d)$. Let $\BD_1 F[d,d)$, and thus $C_0^\downarrow \BD_1 F[d,d)$, also be generated by $z$ and $w$.

With the choice of basis given above, the M\"obius homology Laplacian for $\BD_1F[d,d)$ is the $2 \times 2$ matrix given by:
\[
L_0^\downarrow \BD_1F[d,d) = \begin{pmatrix}
5 & -2/3 \\
-2/3 & 5
\end{pmatrix}.
\]
The eigenvalues of $L_0^\downarrow \BD_1F[d,d)$ are $\frac{13}{3}$ and $\frac{17}{3}$. The characteristic polynomial is
\[
\det \left(\lambda I_2 - L_0^\downarrow \BD_1 F([d,d))\right) = \lambda^2 - 10\lambda + 24\frac{5}{9}.
\]
In this case, $$\mathrm{tr\,} L_0^\downarrow \BD_1 F([d,d)) = 10 = \sum_{[x,y) \prec [d,d)}\bd_1 F([x,y)),$$ illustrating Theorem~\ref{trace-of-L0ab}. Further, $$\det \left( L_0^\downarrow \BD_1 F([d,d))\right)^{1/2} = (221/9)^{1/2} \geq 3,$$
illustrating Theorem~\ref{prop:lower-bound-L0ab}.
\end{ex}

Next, we provide an example for obtaining cycle representatives using the M\"obius homology Laplacian given a multiparameter filtration.

 \begin{ex}\label{ex:multiparam-overcounting-cyclereps}
Consider the multiparameter filtration from Figure~\ref{fig:multiparam-overcounting}. Let $\BD_1F([c,d))$ be generated by the cycles $z_1 = [v_1 v_2] + [v_2 v_4] - [v_1 v_4]$ and $z_2 = [v_2 v_3] + [v_3 v_4] - [v_2 v_4]$, and let $\BD_1 F([a, d))$ be generated by $z_3 = z_1 + z_2$. The M\"obius homology Laplacian for $\BD_1 F([c,d))$ is given by \[
L_0^\downarrow \BD_1F([c,d)) = \begin{pmatrix}
    1/2 & 1/2 \\
    1/2 & 1/2
\end{pmatrix}.
\]
It is easy to check that $\ker L_0^\downarrow \BD_1F([c,d))$ is generated by $\begin{pmatrix}1 \\ -1 \end{pmatrix}$, which represents the cycle \begin{align*}
z &= ([v_1 v_2] + [v_2 v_4] - [v_1 v_4]) - ([v_2 v_3] + [v_3 v_4] - [v_2 v_4]) \\
&= [v_1v_2] +2 [v_2v_4] - [v_1 v_4] - [v_2 v_3] - [v_3v_4].
\end{align*}
Observe that $z_1$ is another representative for $\BD_1 F[c,d)$. We can write $z_1$ as the sum $\frac{1}{2}z + \frac{1}{2}(z_1 + z_2)$ where $z_1 + z_2$ is the image of $z_3$ under $\partial_1^\downarrow$. Up to rescaling, the harmonic representative~$z$ has minimal energy: $\| z_1\|^2 = 3$ whereas $\|\frac{1}{2}z\|^2 = \frac{1}{4}8 = 2$.
Furthermore,
\begin{align*}
    \mathrm{New}(z_1) &= \frac{1}{3},\\
    \mathrm{New}(z) &= \frac{4}{8} = \frac{1}{2}.
\end{align*}
That is, the representative $z$ in $\ker L_0^\downarrow \BD_1 F[c,d)$ maximizes relative new content, illustrating Theorem~\ref{prop:relative-new-content}.

The eigenvalues of $L_0^\downarrow \BD_1F([c,d))$ are 0 and 1. %Indeed, $1 = \mathrm{len}\left([[a,d), [c,d)) \right) = \# \{[a,d)\}$. 
The characteristic polynomial is
\[
\det \left(\lambda I_2 - L_0^\downarrow \BD_1 F([c,d))\right) = \lambda^2 - \lambda.
\]
In this case, $$\mathrm{tr\,} L_0^\downarrow \BD_1 F([c,d)) = 1 = \bd_1 F([a,d))$$ and $$\det \left( L_0^\downarrow \BD_1 F([c,d))\right)^{1/2} = 0^{1/2} \geq 0.$$
\end{ex}

\begin{figure}
\include{figures-tex/multiparam-overcounting}
\vspace{-0.2in}
\caption{Above: A filtration $F: P \to \Delta K$ where $P$ is not totally ordered. Below: $(\Dgm(P), \preceq)$ of $P$ (left), the birth-death module $\BD_1F$ (center), and the persistence diagram $\partial\bd_1F$ (right). In $\BD_1F$, all unlabeled maps are either the zero map or the identity.}
\label{fig:multiparam-overcounting}
\end{figure}

\end{document}